\documentclass[12pt]{article}
\usepackage{amsmath,amssymb,amsthm,amscd,latexsym,}
 
\usepackage{graphicx}

\makeatletter
\renewcommand{\mod}[1]{\allowbreak \if@display \mkern 8mu \else
\mkern 5mu\fi {\operator@font mod}\,\,#1}
\makeatother

\newcommand{\bz}{\mathbb Z}

\newcommand{\aaa}{\mathbb A}
\newcommand{\ddd}{\mathbb D}
\newcommand{\eee}{\mathbb E}

\DeclareMathOperator{\rk}{rk}

\newtheorem{theorem}{Theorem}

\newtheorem{corollary}[theorem]{Corollary}

\newtheorem{conjecture}[theorem]{Conjecture}

\begin{document}
\title{Degenerations of K\"ahlerian K3 surfaces with finite
symplectic automorphism groups}
\date{}
\author{Viacheslav V. Nikulin}
\maketitle

\begin{abstract} Using results of our preprint
"K\"ahlerian K3 surfaces and Niemeier
lattices" \newline arXiv:1109.2879 (and the corresponding
papers), we classify degenerations of K\"ahlerian
K3 surfaces with finite symplectic automorphism groups.
\end{abstract}

\centerline{Dedicated to the memory of Andrey Todorov}

\section{Introduction}
\label{introduction}

Using results of our preprint
"K\"ahlerian K3 surfaces and Niemeier
lattices" arXiv:1109.2879 \cite{Nik7} (and the corresponding
papers), we classify degenerations of K\"ahlerian
K3 surfaces with finite symplectic automorphism groups.

Here we mainly consider degenerations of codimension $1$.
We describe possible Dynkin diagrams and discriminant groups
of the degenerations. See Section \ref{sec:degALL}.
In Section \ref{sec:Markdegen}, one can find detailed
description of possible degenerations using markings by
Niemeier lattices; it uses classification of so called
KahK3 conjugacy classes of automorphism groups
of Niemeier lattices obtained in \cite{Nik7}, \cite{Nik8}
and Section \ref{sec:KahK3cl} here.
In Section \ref{sec:Conj} we formulate an important Conjecture.

In Section \ref{sec:Kummer}, we consider an example of
Kummer surfaces which can be considered as a degeneration of
codimension 1 of K\"ahlerian K3 surfaces with symplectic
automorphism group $(C_2)^4$, and the general theory.

In Section \ref{sec:primembbunim}, we remind to a reader results of \cite{Nik1}
about existence of a primitive embedding of a lattice
(i.e., an integral symmetric bilinear form) into one of even unimodular
lattices of the fixed signature. See Theorem \ref{th:primembb1}.
It is crucial for our classification since by the epimorphicity of
Torelli map for K\"ahlerian K3 surfaces proved by A. Todorov \cite{Tod}
and Y. Siu \cite{Siu}, it permits to describe possible Picard lattices of
K3 surfaces. All the time we use Programs 5 and 6 from Section \ref{sec:Appendix},
Appendix based on Theorem \ref{th:primembb1}.

We hope to give more exact results (e.g. see Conjecture \ref{conj}) 
and consider other cases (for example, degenerations of other 
codimensions) in further publications.

\section{Existence of a primitive embedding of an even lattice
into even unimodular lattices, according to \cite{Nik1}}
\label{sec:primembbunim}

In this paper, we use notations, definitions and results of
\cite{Nik1} about lattices
(that is non-degenerate integral (over $\bz$) symmetric bilinear
forms). In
particular, $\oplus$ denotes the orthogonal sum of lattices,
quadratic forms.
For a prime $p$, we denote by $\bz_p$ the ring of $p$-adic integers,
and by $\bz_p^\ast$ its group of invertible elements.

Let $S$ be a lattice. Let $A_S=S^\ast/S$ be its discriminant group,
and $q_S$ its discriminant
quadratic form on $A_S$ where we assume that the lattice $S$ is
even: that is
$x^2$ is even for any $x\in S$. We denote by $l(A_S)$ the minimal
number
of generators of the finite Abelian group $A_S$, and by $|A_S|$
its order.
For a prime $p$, we denote by ${q_{S}}_p=q_{S\otimes \bz_p}$ the
$p$-component of $q_S$ (equivalently,
the discriminant quadratic form of the $p$-adic lattice
$S\otimes \bz_p$). A quadratic form on a group of order $2$
is denoted by  $q_\theta^{(2)}(2)$. A $p$-adic lattice $K({q_{S}}_p)$
or the rank $l({A_S}_p)$ with the discriminant quadratic form
${q_{S}}_p$ is denoted by $K({q_{S}}_p)$. It is unique, up to
isomorphisms,
for $p\not=2$, and for $p=2$, if
${q_{S}}_2\not\cong q_\theta^{(2)}(2)\oplus
q^\prime$. We have the following result where an embedding
$S\subset L$
of lattices is called {\it primitive} if $L/S$ has no torsion.

\begin{theorem} (Theorem 1.12.2 in \cite{Nik1}).

Let $S$ be an even lattice of the signature $(t_{(+)},t_{(-)})$, and
$l_{(+)}$, $l_{(-)}$ are integers.

Then, there exists a primitive embedding of $S$ into one of
even unimodular lattices
of the signature $(l_{(+)},\,l_{(-)})$ if and only if the following
conditions satisfy:

(1) $l_{(+)}-l_{(-)}\equiv 0\mod 8$;

(2) $l_{(+)}-t_{(+)}\ge 0$, $l_{(-)}-t_{(-)}\ge 0$,
$l_{(+)}+l_{(-)}-t_{(+)}-t_{(-)}\ge l(A_S)$;

(3) $(-1)^{l_{(+)}-t_{(+)}}|A_S|\equiv \det {K(q_{S_p})}\mod ({\bz_p}^\ast)^2$
for each odd prime $p$ such that $l_{(+)}+l_{(-)}-t_{(+)}-t_{(-)}=l(A_{S_p})$;

(4) $|A_S|\equiv \pm \det {K(q_{S_2})}\mod ({\bz_2}^\ast)^2$, if
$l_{(+)}+l_{(-)}-t_{(+)}-t_{(-)}= l(A_{S_2})$ and
${q_{S}}_2\not\cong q_\theta^{(2)}(2)\oplus q^\prime$.
\label{th:primembb1}
\end{theorem}

Remark that if the last inequality in (2) is strict then one does not need the
conditions (3) and (4). If
${q_{S}}_2\cong q_\theta^{(2)}(2)\oplus q^\prime$, then one does not need
the condition (4).

\section{An example related to Kummer surfaces,\\ and the general theory.}
\label{sec:Kummer}

Let us consider K\"ahlerian K3 surfaces which are Kummer sufaces.
They are minimal desingularizations of the quotients of 2-dimensional complex tori by
the involution of the inversion. A general Kummer surface $X$ has the
Picard number 16 and has 16 non-singular rational curves which
define the Dynkin diagram $16\aaa_1$. Moreover, $X$ has the
symplectic automorphism group $(C_2)^4$ defined by translations
of order two elements of the tori. General K\"ahlerian K3 surfaces
with symplectic automorphism group $(C_2)^4$ have Picard number $15$.
See \cite{Nik-1}, \cite{Nik0}.
Thus Kummer surfaces can be considered as the {\it degeneration of codimension $1$}
of K\"ahlerian K3 surfaces with symplectic automorphism group $(C_2)^4$.

This classical example raises the following interesting questions.

Do K\"ahlerian K3 surfaces with finite symplectic automorphism group $(C_2)^4$
have other degenerations of codimension $1$ which are different from Kummer?

What are degenerations of codimension $1$ (and other codimensions) of K3 surfaces
with other finite symplectic automorphism groups?

Using results of our preprint \cite{Nik7}, we can answer this questions.

For example, we show that K\"ahlerian K3 surfaces with finite symplectic automorphism group
$(C_2)^4$ have another degeneration of codimension $1$ which has
the type $4\aaa_1$. It follows from Theorem \ref{theorem:degengeneral} below.
Thus, the Picard lattice of the K3 surface has the rank $16$,
the K3 surface has exactly $4$ non-singular
rational curves which have the Dynkin diagram $4\aaa_1$, and it has the
automorphism group $(C_2)^4$. By \cite{Bri}, \cite{Gra} and \cite{Mum}, it
can be considered as the minimal resolution of singularities of
the complex surface with the corresponding four Du Val singularities of
the type $\aaa_1$ and the automorphism group $(C_2)^4$.

\medskip

Using markings of K\"ahlerian K3 surfaces by Niemeier lattices which was developed in
\cite{Nik7}, this is reduced to the following question. We use notations and results
from \cite{Nik7}.

Let $N_j$ be one of $23$ Niemeier lattices which have non-trivial root systems.
We fix the bases $P(N_j)$ of the root system and consider the automorphism group
$A(N_j)$ of $N_j$ which preserves the bases $P(N_j)$. Let $H\subset A(N_j)$ be
a subgroup of $A(N_j)$. We denote by $(N_j)_H=(N_j^H)^\perp_{N_j} \subset N_j$
the coinvariant sublattice of $H$. We say that $H$ is {\it maximal}
if $H=Clos(H)$ is the maximal
subgroup of $A(N_j)$ with the same coinvariant sublattice $N_H$.
We can restrict considering the maximal case only. A subgroup $H$ is
{\it K\"ahlerian K3 automorphism group} if $N_H$ has a primitive embedding into
the lattice $L_{K3}$ which is an even unimodular lattice
of the signature $(3,19)$.
It is isomorphic to the cohomology lattice $H^2(X,\bz)$ of K\"ahlerian K3 surfaces.
Equivalently, $N_H$ satisfies conditions of Theorem \ref{th:primembb1} for
$l_{(+)}=3$, $l_{(-)}=19$.

Let $\alpha\in P(N_j)$ and $H(\alpha)$ is the orbit of
$H$ for $\alpha\in P(N_j)$. The primitive sublattice
$S=[N_H,\alpha]_{pr}\subset N_j$ contains $H(\alpha)$
and $H(\alpha)$ gives the basis of the system of roots
with square $(-2)$ in $S$.
If $S$ has a primitive embedding into $L_{K3}$ (equivalently,
if $S$ satisfies Theorem \ref{th:primembb1} for $l_{(+)}=3$,
$l_{(-)}=19$), then $S$ can be considered as Picard lattice of
general K\"ahlerian K3 surfaces with symplectic automorphism
group $H$ and classes of non-singular rational curves $H(\alpha)$.
Here we use epimorphicity of the period map for K\"ahlerian
K3 surfaces, see A. Todorov \cite{Tod} and Y. Siu \cite{Siu} and
global Torelli Theorem for K3 surfaces \cite{BR} (see also
\cite{Kul} and \cite{PS}). They permit to describe possible Picard
lattices of
Since $\rk S=\rk N_H + 1$, these K3 surfaces can be considered
as degeneration of codimension $1$ of K\"ahlerian K3 surfaces
with the automorphism group $H$.
The primitive embedding $S\subset N_j$ can be considered as
{\it marking} of the K\"ahlerian K3 surfaces by the Niemeier lattice $N_j$.
See \cite{Nik7} for details.

In \cite{Nik7} we described K\"ahlerian K3 conjugacy classes for
Niemeier lattices $N_j$, $j=1,\dots,23$. Thus, using these classification
and Theorem \ref{th:primembb1}, we can describe degenerations of
codimension $1$ of K\"ahlerian K3 surfaces with finite
symplectic automorphism groups. By considering several orbits of $H$,
we can describe
all degenerations.

\medskip

For example, the group $(C_2)^4$ corresponds to the
following K\"ahlerian K3 conjugacy classes for $N_{23}=N(24A_1)$.
This case can be marked by $N_{23}$ only. See \cite{Nik7}.
We use the basis $\alpha_1,\dots, \alpha_{24}$ with the Dynkin diagram $24\aaa_1$ for the
roots system $24A_1$ of $N_{23}$.

{\bf n=21,} $H\cong C_2^4$ ($|H|=16$, $i=14$):
$\rk N_H=15$, $(N_H)^\ast/N_H \cong
\bz/8\bz\times (\bz/2\bz)^6$.
$$
H_{21,1}=
[(\alpha_{2}\alpha_{20})(\alpha_{3}\alpha_{10})
(\alpha_{5}\alpha_{6})(\alpha_{8}\alpha_{11})
(\alpha_{9}\alpha_{21})(\alpha_{12}\alpha_{22})
(\alpha_{17}\alpha_{23})(\alpha_{19}\alpha_{24}),
$$
$$
(\alpha_{2}\alpha_{19})(\alpha_{3}\alpha_{5})
(\alpha_{6}\alpha_{10})(\alpha_{8}\alpha_{9})
(\alpha_{11}\alpha_{21})(\alpha_{12}\alpha_{23})
(\alpha_{17}\alpha_{22})(\alpha_{20}\alpha_{24}),
$$
$$
(\alpha_{1}\alpha_{16})(\alpha_{2}\alpha_{20})
(\alpha_{3}\alpha_{6})(\alpha_{5}\alpha_{10})
(\alpha_{12}\alpha_{23})(\alpha_{14}\alpha_{18})
(\alpha_{17}\alpha_{22})(\alpha_{19}\alpha_{24}),
$$
$$
(\alpha_{1}\alpha_{14})(\alpha_{2}\alpha_{24})
(\alpha_{3}\alpha_{5})(\alpha_{6}\alpha_{10})
(\alpha_{12}\alpha_{22})(\alpha_{16}\alpha_{18})
(\alpha_{17}\alpha_{23})(\alpha_{19}\alpha_{20})]
$$
with orbits
$
\{\alpha_{1},\alpha_{16},\alpha_{14},\alpha_{18}\},
\{\alpha_{2},\alpha_{20},\alpha_{19},\alpha_{24}\},
\{\alpha_{3},\alpha_{10},\alpha_{5},\alpha_{6}\},
\{\alpha_{8},\alpha_{11},\alpha_{9},\alpha_{21}\},
\newline
\{\alpha_{12},\alpha_{22},\alpha_{23},\alpha_{17}\}
$;

$$
H_{21,2}=
[(\alpha_{1}\alpha_{3})(\alpha_{2}\alpha_{23})
(\alpha_{5}\alpha_{14})(\alpha_{6}\alpha_{16})
(\alpha_{10}\alpha_{18})(\alpha_{12}\alpha_{20})
(\alpha_{17}\alpha_{24})(\alpha_{19}\alpha_{22}),
$$
$$
(\alpha_{1}\alpha_{2})(\alpha_{3}\alpha_{23})
(\alpha_{5}\alpha_{17})(\alpha_{6}\alpha_{12})
(\alpha_{10}\alpha_{22})(\alpha_{14}\alpha_{24})
(\alpha_{16}\alpha_{20})(\alpha_{18}\alpha_{19}),
$$
$$
(\alpha_{1}\alpha_{16})(\alpha_{2}\alpha_{20})
(\alpha_{3}\alpha_{6})(\alpha_{5}\alpha_{10})
(\alpha_{12}\alpha_{23})(\alpha_{14}\alpha_{18})
(\alpha_{17}\alpha_{22})(\alpha_{19}\alpha_{24}),
$$
$$
(\alpha_{1}\alpha_{14})(\alpha_{2}\alpha_{24})
(\alpha_{3}\alpha_{5})(\alpha_{6}\alpha_{10})
(\alpha_{12}\alpha_{22})(\alpha_{16}\alpha_{18})
(\alpha_{17}\alpha_{23})(\alpha_{19}\alpha_{20})]
$$
with orbits
$
\{\alpha_{1},\alpha_{3},\alpha_{2},\alpha_{16},\alpha_{14},
\alpha_{23},\alpha_{6},\alpha_{5},
\alpha_{20},\alpha_{24},\alpha_{18},\alpha_{12},\alpha_{17},
\alpha_{10},\alpha_{19},\alpha_{22}\}
$
(we give only orbits with more than one element).
\medskip

The group $H_{21,1}$ and its orbits with $4$ elements
(there are also four orbits with
$1$ element) give degenerations
of the type $4\aaa_1$ (by Theorem \ref{th:primembb1}, only for these orbits
the lattice $S$ has a primitive embedding into $L_{K3}$). The group $H_{21,2}$
and its orbit with $16$ elements gives the degeneration of the type $16\aaa_1$
which gives the case of Kummer surfaces (by Theorem \ref{th:primembb1},
only for this orbit the lattice $S$ has a primitive embedding into $L_{K3}$).
Thus K\"ahlerian K3 surfaces with symplectic automorphism group $H\cong (C_2)^4$
have codimension $1$ degenerations only of the types $4\aaa_1$ and $16\aaa_1$.
We also calculate the discriminant group $A_S=S^\ast /S$ of $S$ for both
cases. See the case $n=21$ of Theorem \ref{theorem:degengeneral}
in Section \ref{sec:degALL}.

Similarly below we describe degenerations of codimension $1$ which are
marked by all Niemeier lattices. As a result, we obtain classification
of types of degenerations of codimension $1$ of K\"ahlerian K3 surfaces.
We hope to consider more exact classification and degenerations of arbitrary codimension
in further variants of this paper and further publications.

\section{Classification of KahK3 conjugacy classes\\
for $N_j$, $j=1,\dots, 21$.}
\label{sec:KahK3cl}

We use the same notations as in \cite{Nik7} and \cite{Nik8}.

\begin{figure}
\begin{center}
\includegraphics[width=10cm]{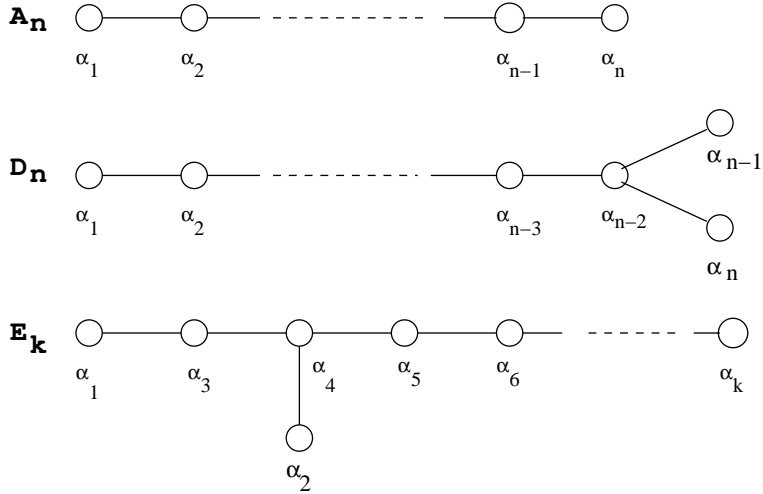}
\end{center}
\caption{Bases of Dynkin diagrams $\aaa_n$, $\ddd_n$, $\eee_k$. }
\label{fig:dyndiagr}
\end{figure}

Below, we use the basis of a root lattice $A_n$, $D_n$ or $E_k$,
$k=6,7,8$,
which is shown on Figure \ref{fig:dyndiagr}.

For $A_n$, $n\ge 1$, we denote
$\varepsilon_1=(\alpha_1+2\alpha_2+\cdots +n\alpha_n)/(n+1)$.
It gives the generator of the discriminant group
$A_n^\ast/A_n\cong \bz/(n+1)\bz$.

For $D_n$, $n\ge 4$ and $n\equiv 0\mod 2$, we denote
$\varepsilon_1=(\alpha_1+\alpha_3+\cdots +\alpha_{n-3}+\alpha_{n-1})/2$,
$\varepsilon_2=(\alpha_{n-1}+\alpha_n)/2$,
$\varepsilon_3=(\alpha_1+\alpha_3+\cdots +\alpha_{n-3}+\alpha_n)/2$. They
give all non-zero elements of the discriminant group
$D_n^\ast/D_n\cong (\bz/2\bz)^2$.

For $D_n$, $n\ge 4$ and $n\equiv 1\mod 2$, we denote
$\varepsilon_1=(\alpha_1+\alpha_3+\cdots +\alpha_{n-2})/2+
\alpha_{n-1}/4-\alpha_n/4$,
$\varepsilon_2=(\alpha_{n-1}+\alpha_n)/2$,
$\varepsilon_3=(\alpha_1+\alpha_3+\cdots +\alpha_{n-2})/2-\alpha_{n-1}/4+
\alpha_n/4$. They give all non-zero elements of
$D_n^\ast/D_n\cong \bz/4\bz$.

For $E_6$, we denote $\varepsilon_1=(\alpha_1-\alpha_3+\alpha_5-\alpha_6)/3$,
$\varepsilon_2=(-\alpha_1+\alpha_3-\alpha_5+\alpha_6)/3$.
They give all non-zero elements of
$E_6^\ast/E_6\cong \bz/3\bz$.

For $E_7$, we denote $\varepsilon_1=(\alpha_2+\alpha_5+\alpha_7)/2$.
It gives the non-zero element of $E_7^\ast/E_7\cong \bz/2\bz$.

If the Dynkin diagram of a root lattice has several connected
components, the second index of a basis numerates
the corresponding connected component.

\vskip1cm

{\bf Case 21.} For the Niemeier lattice
$$
N=N_{21}=N(8A_3)=[8A_3,[3(2001011)]]=
$$
$$
[8A_3,\ -\varepsilon_{1,1}+2\varepsilon_{1,2}+
\varepsilon_{1,5}+\varepsilon_{1,7}+\varepsilon_{1,8},\
-\varepsilon_{1,1}+\varepsilon_{1,2}+2\varepsilon_{1,3}+
\varepsilon_{1,6}+\varepsilon_{1,8},\,
$$
$$
-\varepsilon_{1,1}+\varepsilon_{1,2}+\varepsilon_{1,3}+
2\varepsilon_{1,4}+\varepsilon_{1,7},\
-\varepsilon_{1,1}+\varepsilon_{1,3}+\varepsilon_{1,4}+
2\varepsilon_{1,5}+\varepsilon_{1,8},\
$$
$$
-\varepsilon_{1,1}+\varepsilon_{1,2}+\varepsilon_{1,4}+
\varepsilon_{1,5}+2\varepsilon_{1,6},\
-\varepsilon_{1,1}+\varepsilon_{1,3}+\varepsilon_{1,5}+
\varepsilon_{1,6}+2\varepsilon_{1,7},\
$$
$$
-\varepsilon_{1,1}+\varepsilon_{1,4}+\varepsilon_{1,6}+
\varepsilon_{1,7}+2\varepsilon_{1,8}]
$$
the group $A(N_{21})$ has the
order $2688$, and it is generated by
$$
\widetilde{F1}_7=
(\alpha_{1,2}\alpha_{1,3}\alpha_{1,4}\alpha_{1,5}\alpha_{1,6}\alpha_{1,7}\alpha_{1,8})
(\alpha_{2,2}\alpha_{2,3}\alpha_{2,4}\alpha_{2,5}\alpha_{2,6}\alpha_{2,7}\alpha_{2,8})
$$
$$
(\alpha_{3,2}\alpha_{3,3}\alpha_{3,4}\alpha_{3,5}\alpha_{3,6}\alpha_{3,7}\alpha_{3,8}),
$$
$$
\widetilde{F1}_4=
(\alpha_{1,2}\alpha_{3,2})
(\alpha_{1,3}\alpha_{1,4}\alpha_{3,7}\alpha_{1,5})
(\alpha_{2,3}\alpha_{2,4}\alpha_{2,7}\alpha_{2,5})
(\alpha_{3,3}\alpha_{3,4}\alpha_{1,7}\alpha_{3,5})
$$
$$
(\alpha_{1,6}\alpha_{1,8}\alpha_{3,6}\alpha_{3,8})
(\alpha_{2,6}\alpha_{2,8}),
$$
$$
\widetilde{T}_{12}=(\alpha_{1,1}\alpha_{1,2})
(\alpha_{3,1}\alpha_{3,2})(\alpha_{1,3}\alpha_{1,7})
(\alpha_{3,3}\alpha_{3,7})(\alpha_{1,4}\alpha_{1,5})
(\alpha_{3,4}\alpha_{3,5})
$$
$$
(\alpha_{1,6}\alpha_{1,8})
(\alpha_{3,6}\alpha_{3,8})
$$
See \cite[Ch. 16]{CS} and \cite{Nik7}, \cite{Nik8}.

Using classification by Mukai \cite{Muk} and its
refinement by Xiao \cite{Xiao}
of abstract finite symplectic automorphism groups of K3,
and GAP Program \cite{GAP}
(similarly to Hashimoto in \cite{Hash}),
like for $N_{23}$ and $N_{22}$ in \cite{Nik7} and \cite{Nik8},
we obtain the following classification.
All the time we use Theorem \ref{th:primembb1} (we use Program 5 in
Section \ref{sec:Appendix} based on this Theorem).

The same we do for all other cases below.
This style of classification will be more convenient for studying
of degenerations of K3 than our classification in \cite{Nik7} and \cite{Nik8}.

\vskip1cm

\centerline {\bf Classification of KahK3 conjugacy classes for $A(N_{21})$:}

\vskip1cm

{\bf n=74,} $H\cong L_2(7)$ ($|H|=168$, $i=42$):
$\rk N_H=19$ and $(N_H)^\ast/N_H \cong \bz/28\bz \times \bz/7\bz$.
$$
H_{74,1}=
[(\alpha_{1,2}\alpha_{1,3}\alpha_{1,5})(\alpha_{2,2}\alpha_{2,3}\alpha_{2,5})
(\alpha_{3,2}\alpha_{3,3}\alpha_{3,5})
(\alpha_{1,4}\alpha_{1,7}\alpha_{1,6})
(\alpha_{2,4}\alpha_{2,7}\alpha_{2,6})
(\alpha_{3,4}\alpha_{3,7}\alpha_{3,6}),
$$
$$
(\alpha_{1,2}\alpha_{3,2})(\alpha_{1,3}\alpha_{3,6})
(\alpha_{2,3}\alpha_{2,6})(\alpha_{3,3}\alpha_{1,6})
(\alpha_{1,4}\alpha_{3,4})(\alpha_{1,7}\alpha_{1,8})
(\alpha_{2,7}\alpha_{2,8})
(\alpha_{3,7}\alpha_{3,8})]
$$
with orbits
$
\{\alpha_{1,2}, \alpha_{1,3},\alpha_{3,2}, \alpha_{1,5},\alpha_{3,6}, \alpha_{3,3},
\alpha_{3,4}, \alpha_{3,5}, \alpha_{1,6}, \alpha_{3,7}, \alpha_{1,4}, \alpha_{3,8},
\alpha_{1,7},\alpha_{1,8}\},\
\linebreak
\{\alpha_{2,2}, \alpha_{2,3}, \alpha_{2,5}, \alpha_{2,6}, \alpha_{2,4},
\alpha_{2,7}, \alpha_{2,8}\}
$.

\medskip


{\bf n=51,} $H\cong C_2\times {\mathfrak S}_4$ ($|H|=48$, $i=48$):
$\rk N_H=18$, $(N_H)^\ast/N_H\cong
(\bz/12\bz)^2\times (\bz/2\bz)^2$.
$$
H_{51,1}=[(\alpha_{1,2}\alpha_{3,2})(\alpha_{1,3}\alpha_{3,4})
(\alpha_{2,3}\alpha_{2,4})(\alpha_{3,3}\alpha_{1,4})
(\alpha_{1,5}\alpha_{1,7})(\alpha_{2,5}\alpha_{2,7})
(\alpha_{3,5}\alpha_{3,7})(\alpha_{1,6}\alpha_{3,6}),
$$
$$
 (\alpha_{1,1}\alpha_{3,1})(\alpha_{1,3}\alpha_{3,7}
 \alpha_{3,3}\alpha_{1,7})(\alpha_{2,3}\alpha_{2,7})
(\alpha_{1,4}\alpha_{3,6}\alpha_{3,5}\alpha_{3,8})
(\alpha_{2,4}\alpha_{2,6}\alpha_{2,5}\alpha_{2,8})
(\alpha_{3,4}\alpha_{1,6}\alpha_{1,5}\alpha_{1,8})]
$$
with orbits
$
\{\alpha_{1,1}, \alpha_{3,1} \},\  \{\alpha_{1,2}, \alpha_{3,2}\},\
\{\alpha_{1,3}, \alpha_{3,4}, \alpha_{3,7}, \alpha_{1,6}, \alpha_{3,5},
\alpha_{3,3},
\alpha_{3,6}, \alpha_{1,5}, \alpha_{3,8}, \alpha_{1,4},
\linebreak
\alpha_{1,7}, \alpha_{1,8}\},\
\{\alpha_{2,3}, \alpha_{2,4}, \alpha_{2,7}, \alpha_{2,6},
\alpha_{2,5}, \alpha_{2,8}\}
$;
$$
H_{51,2}=[
(\alpha_{1,2}\alpha_{1,3})
(\alpha_{2,2}\alpha_{2,3})
(\alpha_{3,2}\alpha_{3,3})
(\alpha_{1,4}\alpha_{3,4})
(\alpha_{1,5}\alpha_{3,8})
(\alpha_{2,5}\alpha_{2,8})
(\alpha_{3,5}\alpha_{1,8})
(\alpha_{1,6}\alpha_{3,6}),
$$
$$
(\alpha_{1,1}\alpha_{3,1})
(\alpha_{1,3}\alpha_{3,7}\alpha_{3,3}\alpha_{1,7})
(\alpha_{2,3}\alpha_{2,7})
(\alpha_{1,4}\alpha_{3,6}\alpha_{3,5}\alpha_{3,8})
(\alpha_{2,4}\alpha_{2,6}\alpha_{2,5}\alpha_{2,8})
(\alpha_{3,4}\alpha_{1,6}\alpha_{1,5}\alpha_{1,8})]
$$
with orbits
$
\{\alpha_{1,1},\alpha_{3,1}\},\
\{\alpha_{1,2},\alpha_{1,3},\alpha_{3,7},\alpha_{3,3},\alpha_{3,2},\alpha_{1,7} \},\
\{\alpha_{2,2},\alpha_{2,3},
\alpha_{2,7}\},\ \linebreak
\{\alpha_{1,4},\alpha_{3,4},\alpha_{3,6},\alpha_{1,6},\alpha_{3,5},
\alpha_{1,5},\alpha_{1,8},\alpha_{3,8}\},\
\{\alpha_{2,4},\alpha_{2,6},\alpha_{2,5},\alpha_{2,8}\}
$.

\medskip


{\bf n=35,} $H\cong C_2\times {\mathfrak A}_4$, ($|H|=24$, $i=13$):
$$
H_{35,1}=
[(\alpha_{1,1}\alpha_{3,1})
(\alpha_{1,2}\alpha_{3,2})
(\alpha_{1,3}\alpha_{3,4}\alpha_{1,6}\alpha_{3,3}\alpha_{1,4}\alpha_{3,6})
(\alpha_{2,3}\alpha_{2,4}\alpha_{2,6})
(\alpha_{1,5}\alpha_{3,8}\alpha_{1,7}\alpha_{3,5}\alpha_{1,8}\alpha_{3,7})
$$
$$
(\alpha_{2,5}\alpha_{2,8}\alpha_{2,7}),
$$
$$
(\alpha_{1,1}\alpha_{3,1})
(\alpha_{1,2}\alpha_{3,2})
(\alpha_{1,3}\alpha_{1,4}\alpha_{1,8}\alpha_{3,3}\alpha_{3,4}\alpha_{3,8})
(\alpha_{2,3}\alpha_{2,4}\alpha_{2,8})
(\alpha_{1,5}\alpha_{1,6}\alpha_{3,7}\alpha_{3,5}\alpha_{3,6}\alpha_{1,7})
$$
$$
(\alpha_{2,5}\alpha_{2,6}\alpha_{2,7})]
$$
with $Clos(H_{35,1})=H_{51,1}$ above;
$$
H_{35,2}=
[(\alpha_{1,1}\alpha_{3,1})
(\alpha_{1,2}\alpha_{1,3}\alpha_{3,7}\alpha_{3,2}\alpha_{3,3}\alpha_{1,7})
(\alpha_{2,2}\alpha_{2,3}\alpha_{2,7})
(\alpha_{1,4}\alpha_{1,6}\alpha_{3,8}\alpha_{3,4}\alpha_{3,6}\alpha_{1,8})
$$
$$
(\alpha_{2,4}\alpha_{2,6}\alpha_{2,8})
(\alpha_{1,5}\alpha_{3,5}),
$$
$$
(\alpha_{1,1}\alpha_{3,1})
(\alpha_{1,2}\alpha_{1,3}\alpha_{1,7}\alpha_{3,2}\alpha_{3,3}\alpha_{3,7})
(\alpha_{2,2}\alpha_{2,3}\alpha_{2,7})
(\alpha_{1,4}\alpha_{1,5}\alpha_{3,6}\alpha_{3,4}\alpha_{3,5}\alpha_{1,6})
$$
$$
(\alpha_{2,4}\alpha_{2,5}\alpha_{2,6})
(\alpha_{1,8}\alpha_{3,8})]
$$
with $Clos(H_{35,2})=H_{51,2}$ above.

\medskip

{\bf n=34}, $H\cong {\mathfrak S}_4$ ($|H|=24$, $i=12$):
$\rk N_H=17$ and $(N_H)^\ast/N_H\cong
(\bz/12\bz)^2\times \bz/4\bz$.
$$
H_{34,1}=[
(\alpha_{1,3}\alpha_{3,4}\alpha_{1,8})
(\alpha_{2,3}\alpha_{2,4}\alpha_{2,8})
(\alpha_{3,3}\alpha_{1,4}\alpha_{3,8})
(\alpha_{1,5}\alpha_{3,6}\alpha_{3,7})
$$
$$
(\alpha_{2,5}\alpha_{2,6}\alpha_{2,7})
(\alpha_{3,5}\alpha_{1,6}\alpha_{1,7}),
$$
$$
(\alpha_{1,1}\alpha_{3,1})
(\alpha_{1,3}\alpha_{3,7}\alpha_{3,3}\alpha_{1,7})
(\alpha_{2,3}\alpha_{2,7})
(\alpha_{1,4}\alpha_{3,6}\alpha_{3,5}\alpha_{3,8})
$$
$$
(\alpha_{2,4}\alpha_{2,6}\alpha_{2,5}\alpha_{2,8})
(\alpha_{3,4}\alpha_{1,6}\alpha_{1,5}\alpha_{1,8})]
$$
with orbits
$
\{\alpha_{1,1}, \alpha_{3,1} \},\  \{\alpha_{1,3}, \alpha_{3,4},
\alpha_{3,7}, \alpha_{1,8}, \alpha_{1,6}, \alpha_{1,5},
\alpha_{3,3}, \alpha_{1,7},
\alpha_{3,6},
\alpha_{1,4}, \alpha_{3,5}, \alpha_{3,8}\},\
\linebreak
 \{\alpha_{2,3}, \alpha_{2,4},
\alpha_{2,7}, \alpha_{2,8}, \alpha_{2,6}, \alpha_{2,5} \}
$;
$$
H_{34,2}=[
(\alpha_{1,2}\alpha_{1,3}\alpha_{3,7})
(\alpha_{2,2}\alpha_{2,3}\alpha_{2,7})
(\alpha_{3,2}\alpha_{3,3}\alpha_{1,7})
(\alpha_{1,4}\alpha_{3,8}\alpha_{3,5})
$$
$$
(\alpha_{2,4}\alpha_{2,8}\alpha_{2,5})
(\alpha_{3,4}\alpha_{1,8}\alpha_{1,5}),
$$
$$
 (\alpha_{1,1}\alpha_{3,1})
 (\alpha_{1,3}\alpha_{3,7}\alpha_{3,3}\alpha_{1,7})
 (\alpha_{2,3}\alpha_{2,7})
 (\alpha_{1,4}\alpha_{3,6}\alpha_{3,5}\alpha_{3,8})
 $$
 $$
 (\alpha_{2,4}\alpha_{2,6}\alpha_{2,5}\alpha_{2,8})
 (\alpha_{3,4}\alpha_{1,6}\alpha_{1,5}\alpha_{1,8})]
$$
with orbits
$
\{\alpha_{1,1},\alpha_{3,1}\},\
\{\alpha_{1,2},\alpha_{1,3},\alpha_{3,7},
\alpha_{3,3},\alpha_{1,7},\alpha_{3,2}\},\
\{\alpha_{2,2},\alpha_{2,3},
\alpha_{2,7}\},\
\{\alpha_{1,4},\alpha_{3,8},\alpha_{3,6},\linebreak
\alpha_{3,5}\},\
\{\alpha_{2,4},\alpha_{2,8},\alpha_{2,6},\alpha_{2,5}\},\
\{\alpha_{3,4},\alpha_{1,8},\alpha_{1,6},\alpha_{1,5}\};
$

\medskip

$$
H_{34,3}=[
(\alpha_{1,1}\alpha_{1,3}\alpha_{1,7})
(\alpha_{2,1}\alpha_{2,3}\alpha_{2,7})
(\alpha_{3,1}\alpha_{3,3}\alpha_{3,7})
(\alpha_{1,4}\alpha_{1,8}\alpha_{1,6})
$$
$$
(\alpha_{2,4}\alpha_{2,8}\alpha_{2,6})
(\alpha_{3,4}\alpha_{3,8}\alpha_{3,6}),
$$
$$
(\alpha_{1,1}\alpha_{3,1})
(\alpha_{1,3}\alpha_{3,7}\alpha_{3,3}\alpha_{1,7})
(\alpha_{2,3}\alpha_{2,7})
(\alpha_{1,4}\alpha_{3,6}\alpha_{3,5}\alpha_{3,8})
$$
$$
(\alpha_{2,4}\alpha_{2,6}\alpha_{2,5}\alpha_{2,8})
(\alpha_{3,4}\alpha_{1,6}\alpha_{1,5}\alpha_{1,8}) ]
$$
with orbits
$
\{\alpha_{1,1},\alpha_{1,3},\alpha_{3,1},\alpha_{1,7},
\alpha_{3,7},\alpha_{3,3}\},\
\{\alpha_{2,1},\alpha_{2,3},\alpha_{2,7}\},\
\{\alpha_{1,4},
\alpha_{1,8},\alpha_{3,6},\alpha_{1,6},
\alpha_{3,4},\alpha_{3,5},
\linebreak
\alpha_{1,5},\alpha_{3,8}\},\
\{\alpha_{2,4},\alpha_{2,8},\alpha_{2,6},\alpha_{2,5}\}
$.

\medskip


{\bf n=33,} $H\cong C_7\rtimes C_3$ ($|H|=21$, $i=1$):
$\rk N_H=18$ and $(N_H)^\ast/N_H
\cong (\bz/7\bz)^3$.
$$
H_{33,1}=
[(\alpha_{1,3}\alpha_{1,5}\alpha_{3,8})
(\alpha_{2,3}\alpha_{2,5}\alpha_{2,8})
(\alpha_{3,3}\alpha_{3,5}\alpha_{1,8})
(\alpha_{1,4}\alpha_{3,6}\alpha_{1,7})
$$
$$
(\alpha_{2,4}\alpha_{2,6}\alpha_{2,7})
(\alpha_{3,4}\alpha_{1,6}\alpha_{3,7}),
$$
$$
(\alpha_{1,2}\alpha_{3,4}\alpha_{1,6}\alpha_{3,5}\alpha_{3,7}\alpha_{3,3}\alpha_{1,8})
(\alpha_{2,2}\alpha_{2,4}\alpha_{2,6}\alpha_{2,5}\alpha_{2,7}\alpha_{2,3}\alpha_{2,8})
$$
$$
(\alpha_{3,2}\alpha_{1,4}\alpha_{3,6}\alpha_{1,5}\alpha_{1,7}\alpha_{1,3}\alpha_{3,8})]
$$
with orbits
$
\{\alpha_{1,2},\alpha_{3,4},\alpha_{1,6},\alpha_{3,7},
\alpha_{3,5},\alpha_{3,3},\alpha_{1,8}\},\
\{\alpha_{2,2},\alpha_{2,4},\alpha_{2,6},\alpha_{2,7},
\alpha_{2,5},\alpha_{2,3},\alpha_{2,8}\},\ \linebreak
\{\alpha_{3,2},\alpha_{1,4},\alpha_{3,6},\alpha_{1,7},
\alpha_{1,5},\alpha_{1,3},\alpha_{3,8}\}
$.

\medskip

{\bf n=22,} $H\cong C_2\times D_8$ ($|H|=16$, $i=11$):
$\rk N_H=16$, $(N_H)^\ast/N_H\cong
(\bz/4\bz)^4\times (\bz/2\bz)^2$.
$$
H_{22,1}=[
(\alpha_{1,3}\alpha_{3,3})(\alpha_{1,4}\alpha_{3,5})
(\alpha_{2,4}\alpha_{2,5})(\alpha_{3,4}\alpha_{1,5})
(\alpha_{1,6}\alpha_{1,8})
$$
$$
(\alpha_{2,6}\alpha_{2,8})
(\alpha_{3,6}\alpha_{3,8})(\alpha_{1,7}\alpha_{3,7}),
$$
$$
(\alpha_{1,2}\alpha_{3,2})(\alpha_{1,3}\alpha_{3,6})
(\alpha_{2,3}\alpha_{2,6})(\alpha_{3,3}\alpha_{1,6})
(\alpha_{1,4}\alpha_{3,4})(\alpha_{1,7}\alpha_{1,8})
$$
$$
(\alpha_{2,7}\alpha_{2,8})(\alpha_{3,7}\alpha_{3,8}),
$$
$$
(\alpha_{1,1}\alpha_{3,1})(\alpha_{1,2}\alpha_{3,2})
(\alpha_{1,3}\alpha_{3,3})(\alpha_{1,4}\alpha_{3,4})
(\alpha_{1,5}\alpha_{3,5})(\alpha_{1,6}\alpha_{3,6})
$$
$$
(\alpha_{1,7}\alpha_{3,7})(\alpha_{1,8}\alpha_{3,8}) ]
$$
with orbits
$
\{\alpha_{1,1},\alpha_{3,1}\},\  \{\alpha_{1,2},\alpha_{3,2}\},\
\{\alpha_{1,3},\alpha_{3,3},\alpha_{3,6},\alpha_{1,6},\alpha_{3,8},
\alpha_{1,8},
\alpha_{3,7},\alpha_{1,7}\},\ \linebreak
\{\alpha_{2,3},\alpha_{2,6},\alpha_{2,8},\alpha_{2,7}\},\
\{\alpha_{1,4},\alpha_{3,5},\alpha_{3,4},\alpha_{1,5}\},\
\{\alpha_{2,4},\alpha_{2,5}\}
$.

\medskip


{\bf n=18,} $H\cong D_{12}$ ($|H|=12$, $i=4$):
$\rk N_H=16$ and $(N_H)^\ast/N_H\cong (\bz/6\bz)^4$.
$$
H_{18,1}=[
(\alpha_{1,2}\alpha_{3,2})(\alpha_{1,3}\alpha_{3,3})
(\alpha_{1,4}\alpha_{3,6})(\alpha_{2,4}\alpha_{2,6})
(\alpha_{3,4}\alpha_{1,6})(\alpha_{1,5}\alpha_{1,8})
$$
$$
(\alpha_{2,5}\alpha_{2,8})(\alpha_{3,5}\alpha_{3,8}),
$$
$$
(\alpha_{1,1}\alpha_{3,1})(\alpha_{1,2}\alpha_{3,2})
(\alpha_{1,3}\alpha_{3,4}\alpha_{1,6}\alpha_{3,3}\alpha_{1,4}\alpha_{3,6})
(\alpha_{2,3}\alpha_{2,4}\alpha_{2,6})
$$
$$
(\alpha_{1,5}\alpha_{3,8}\alpha_{1,7}\alpha_{3,5}\alpha_{1,8}\alpha_{3,7})
(\alpha_{2,5}\alpha_{2,8}\alpha_{2,7})]
$$
with orbits
$
\{\alpha_{1,1},\alpha_{3,1}\},\  \{\alpha_{1,2},\alpha_{3,2}\},\
\{\alpha_{1,3},\alpha_{3,3},\alpha_{3,4},\alpha_{1,4},\alpha_{1,6},\alpha_{3,6}\},\
 \{\alpha_{2,3},\alpha_{2,4},\alpha_{2,6}\},\ \linebreak
\{\alpha_{1,5},\alpha_{1,8},\alpha_{3,8},\alpha_{3,7},\alpha_{3,5},\alpha_{1,7}\},\
\{\alpha_{2,5},\alpha_{2,8},\alpha_{2,7}\}$.

\medskip

{\bf n=17,} $H\cong {\mathfrak A}_4$ ($|H|=12$, $i=3$):
$\rk N_H=16$ and $(N_H)^\ast/N_H\cong (\bz/12\bz)^2\times (\bz/2\bz)^2$.
$$
H_{17,1}=[
(\alpha_{1,2}\alpha_{1,3}\alpha_{3,8})(\alpha_{2,2}\alpha_{2,3}\alpha_{2,8})
(\alpha_{3,2}\alpha_{3,3}\alpha_{1,8})(\alpha_{1,4}\alpha_{1,6}\alpha_{3,7})
$$
$$
(\alpha_{2,4}\alpha_{2,6}\alpha_{2,7})(\alpha_{3,4}\alpha_{3,6}\alpha_{1,7}),
$$
$$
(\alpha_{1,2}\alpha_{3,2})(\alpha_{1,3}\alpha_{3,6})(\alpha_{2,3}\alpha_{2,6})
(\alpha_{3,3}\alpha_{1,6})(\alpha_{1,4}\alpha_{3,4})(\alpha_{1,7}\alpha_{1,8})
$$
$$
(\alpha_{2,7}\alpha_{2,8})(\alpha_{3,7}\alpha_{3,8})]
$$
with orbits
$
\{\alpha_{1,2},\alpha_{1,3},\alpha_{3,2},\alpha_{3,8},\alpha_{3,6},\alpha_{3,3},
\alpha_{3,7}, \alpha_{1,7},\alpha_{1,8},\alpha_{1,6},\alpha_{1,4},
\alpha_{3,4}\},\
\{\alpha_{2,2},\alpha_{2,3},\alpha_{2,8}, \linebreak
\alpha_{2,6},\alpha_{2,7},\alpha_{2,4}\}
$;
$$
H_{17,2}=[
(\alpha_{1,2}\alpha_{1,4}\alpha_{3,5})(\alpha_{2,2}\alpha_{2,4}\alpha_{2,5})
(\alpha_{3,2}\alpha_{3,4}\alpha_{1,5})(\alpha_{1,6}\alpha_{3,7}\alpha_{3,8})
$$
$$
(\alpha_{2,6}\alpha_{2,7}\alpha_{2,8})(\alpha_{3,6}\alpha_{1,7}\alpha_{1,8}),
$$
$$
(\alpha_{1,2}\alpha_{3,2})(\alpha_{1,3}\alpha_{3,6})(\alpha_{2,3}\alpha_{2,6})
(\alpha_{3,3}\alpha_{1,6})(\alpha_{1,4}\alpha_{3,4})(\alpha_{1,7}\alpha_{1,8})
$$
$$
(\alpha_{2,7}\alpha_{2,8})(\alpha_{3,7}\alpha_{3,8})]
$$
with orbits
$
\{\alpha_{1,2},\alpha_{1,4},\alpha_{3,2},\alpha_{3,5},\alpha_{3,4},\alpha_{1,5}\},\
\{\alpha_{2,2},\alpha_{2,4},\alpha_{2,5}\},\  \{\alpha_{1,3},
\alpha_{3,6},\alpha_{1,7},\alpha_{1,8}\},\
\{\alpha_{2,3},\alpha_{2,6}, \linebreak
\alpha_{2,7},\alpha_{2,8}\},\
\{\alpha_{3,3},\alpha_{1,6},\alpha_{3,7}, \alpha_{3,8}\}
$.


\medskip

{\bf n=11,} $H\cong C_2\times C_4$ ($|H|=8$, $i=2$):
$$
H_{11,1}=[
(\alpha_{1,1}\alpha_{3,1})(\alpha_{1,2}\alpha_{3,2})
(\alpha_{1,3}\alpha_{3,3})(\alpha_{1,4}\alpha_{3,4})
(\alpha_{1,5}\alpha_{3,5})(\alpha_{1,6}\alpha_{3,6})
$$
$$
(\alpha_{1,7}\alpha_{3,7})(\alpha_{1,8}\alpha_{3,8}),
$$
$$
(\alpha_{1,2}\alpha_{3,2})(\alpha_{1,3}\alpha_{1,6}\alpha_{1,7}\alpha_{3,8})
(\alpha_{2,3}\alpha_{2,6}\alpha_{2,7}\alpha_{2,8})
(\alpha_{3,3}\alpha_{3,6}\alpha_{3,7}\alpha_{1,8})
$$
$$
(\alpha_{1,4}\alpha_{3,5}\alpha_{3,4}\alpha_{1,5})(\alpha_{2,4}\alpha_{2,5})]
$$
with $Clos(H_{11,1})=H_{22,1}$ above.

\medskip

{\bf n=10,} $H\cong D_8$ ($|H|=8$, $i=3$):
$\rk N_H=15$ and $(N_H)^\ast/N_H\cong
(\bz/4\bz)^5$.
$$
H_{10,1}=[
(\alpha_{1,3}\alpha_{3,3})(\alpha_{1,4}\alpha_{3,5})(\alpha_{2,4}\alpha_{2,5})
(\alpha_{3,4}\alpha_{1,5})(\alpha_{1,6}\alpha_{1,8})
$$
$$
(\alpha_{2,6}\alpha_{2,8})
(\alpha_{3,6}\alpha_{3,8})(\alpha_{1,7}\alpha_{3,7}),
$$
$$
(\alpha_{1,2}\alpha_{3,2})(\alpha_{1,3}\alpha_{3,6})(\alpha_{2,3}\alpha_{2,6})
(\alpha_{3,3}\alpha_{1,6})(\alpha_{1,4}\alpha_{3,4})(\alpha_{1,7}\alpha_{1,8})
$$
$$
(\alpha_{2,7}\alpha_{2,8})(\alpha_{3,7}\alpha_{3,8})]
$$
with orbits
$
\{\alpha_{1,2},\alpha_{3,2}\},\
\{\alpha_{1,3},\alpha_{3,3},\alpha_{3,6},\alpha_{1,6},\alpha_{3,8},
\alpha_{1,8},\alpha_{3,7},\alpha_{1,7}\},\
\{\alpha_{2,3},\alpha_{2,6},\alpha_{2,8},\alpha_{2,7}\},\linebreak
\{\alpha_{1,4},\alpha_{3,5},\alpha_{3,4},\alpha_{1,5}\},\  \{\alpha_{2,4},\alpha_{2,5}\}
$;
$$
H_{10,2}=[
(\alpha_{1,1}\alpha_{3,1})(\alpha_{1,2}\alpha_{3,2})(\alpha_{1,4}\alpha_{1,5})
(\alpha_{2,4}\alpha_{2,5})(\alpha_{3,4}\alpha_{3,5})(\alpha_{1,6}\alpha_{3,8})
$$
$$
(\alpha_{2,6}\alpha_{2,8})(\alpha_{3,6}\alpha_{1,8}),
$$
$$
(\alpha_{1,2}\alpha_{3,2})(\alpha_{1,3}\alpha_{3,6})(\alpha_{2,3}\alpha_{2,6})
(\alpha_{3,3}\alpha_{1,6})(\alpha_{1,4}\alpha_{3,4})(\alpha_{1,7}\alpha_{1,8})
$$
$$
(\alpha_{2,7}\alpha_{2,8})(\alpha_{3,7}\alpha_{3,8})]
$$
with orbits
$
\{\alpha_{1,1},\alpha_{3,1}\},\  \{\alpha_{1,2},\alpha_{3,2}\},\
\{\alpha_{1,3},\alpha_{3,6},\alpha_{1,8},\alpha_{1,7}\},\
\{\alpha_{2,3},\alpha_{2,6},
\alpha_{2,8},\alpha_{2,7}\},\ \linebreak
 \{\alpha_{3,3},\alpha_{1,6},\alpha_{3,8},\alpha_{3,7}\},\
\{\alpha_{1,4},\alpha_{1,5},\alpha_{3,4},\alpha_{3,5}\},\
\{\alpha_{2,4},\alpha_{2,5}\}
$.


\medskip

{\bf n=9,} $H\cong C_2^3$ ($|H|=8$, $i=5$):
$\rk N_H=14$, $(N_H)^\ast/N_H \cong (\bz/4\bz)^2\times (\bz/2\bz)^6$.
$$
H_{9,1}=[
(\alpha_{1,1}\alpha_{3,1})(\alpha_{1,2}\alpha_{3,2})
(\alpha_{1,3}\alpha_{3,3})(\alpha_{1,4}\alpha_{3,4})
(\alpha_{1,5}\alpha_{3,5})
$$
$$
(\alpha_{1,6}\alpha_{3,6})
(\alpha_{1,7}\alpha_{3,7})(\alpha_{1,8}\alpha_{3,8}),
$$
$$
(\alpha_{1,3}\alpha_{1,7})(\alpha_{2,3}\alpha_{2,7})
(\alpha_{3,3}\alpha_{3,7})(\alpha_{1,4}\alpha_{3,4})
(\alpha_{1,5}\alpha_{3,5})
$$
$$
(\alpha_{1,6}\alpha_{3,8})
(\alpha_{2,6}\alpha_{2,8})(\alpha_{3,6}\alpha_{1,8}),
$$
$$
(\alpha_{4 1,2}\alpha_{6 3,2})(\alpha_{7 1,3}\alpha_{18 3,6})
(\alpha_{2,3}\alpha_{2,6})(\alpha_{3,3}\alpha_{1,6})
(\alpha_{1,4}\alpha_{3,4})
$$
$$
(\alpha_{1,7}\alpha_{1,8})
(\alpha_{2,7}\alpha_{2,8})(\alpha_{3,7}\alpha_{3,8})]
$$
with orbits
$
\{\alpha_{1,1},\alpha_{3,1}\},\  \{\alpha_{1,2},\alpha_{3,2}\},\
\{\alpha_{1,3},\alpha_{3,3},
\alpha_{1,7},\alpha_{3,6},\alpha_{3,7},\alpha_{1,6},
\alpha_{1,8},\alpha_{3,8}\},\ \linebreak
\{\alpha_{2,3},\alpha_{2,7},\alpha_{2,6},\alpha_{2,8}\},\
\{\alpha_{1,4},\alpha_{3,4}\},\  \{\alpha_{1,5},\alpha_{3,5}\}
$;
$$
H_{9,2}=[
(\alpha_{1,1}\alpha_{3,1})(\alpha_{1,2}\alpha_{3,2})
(\alpha_{1,3}\alpha_{3,3})(\alpha_{1,4}\alpha_{3,4})
(\alpha_{1,5}\alpha_{3,5})
$$
$$
(\alpha_{1,6}\alpha_{3,6})
(\alpha_{1,7}\alpha_{3,7})(\alpha_{1,8}\alpha_{3,8}),
$$
$$
(\alpha_{1,2}\alpha_{1,4})(\alpha_{2,2}\alpha_{2,4})
(\alpha_{3,2}\alpha_{3,4})(\alpha_{1,3}\alpha_{3,3})
(\alpha_{1,6}\alpha_{3,6})
$$
$$
(\alpha_{1,7}\alpha_{3,8})
(\alpha_{2,7}\alpha_{2,8})(\alpha_{3,7}\alpha_{1,8}),
$$
$$
(\alpha_{1,2}\alpha_{3,2})(\alpha_{1,3}\alpha_{3,6})
(\alpha_{2,3}\alpha_{2,6})(\alpha_{3,3}\alpha_{1,6})
(\alpha_{1,4}\alpha_{3,4})
$$
$$
(\alpha_{1,7}\alpha_{1,8})
(\alpha_{2,7}\alpha_{2,8})(\alpha_{3,7}\alpha_{3,8}) ]
$$
with orbits
$
\{\alpha_{1,1},\alpha_{3,1}\},\  \{\alpha_{1,2},\alpha_{3,2},\alpha_{1,4},
\alpha_{3,4}\},\  \{\alpha_{2,2},\alpha_{2,4}\},\  \{\alpha_{1,3},\alpha_{3,3},
\alpha_{3,6},\alpha_{1,6}\},\  \linebreak
\{\alpha_{2,3},\alpha_{2,6}\},\
\{\alpha_{1,5},\alpha_{3,5}\},\ \{\alpha_{1,7},\alpha_{3,7},\alpha_{3,8},\alpha_{1,8}\},\
\{\alpha_{2,7},\alpha_{2,8}\}
$.

\medskip


{\bf n=8,} $H\cong C_7$ ($|H|=7$, $i=1$):
$$
H_{8,1}=
[(\alpha_{1,2}\alpha_{3,4}\alpha_{1,6}\alpha_{3,5}\alpha_{3,7}\alpha_{3,3}\alpha_{1,8})
(\alpha_{2,2}\alpha_{2,4}\alpha_{2,6}\alpha_{2,5}\alpha_{2,7}\alpha_{2,3}\alpha_{2,8})
$$
$$
(\alpha_{3,2}\alpha_{1,4}\alpha_{3,6}\alpha_{1,5}\alpha_{1,7}\alpha_{1,3}\alpha_{3,8})]
$$
with $Clos(H_{8,1})=H_{33,1}$ above.

\medskip

{\bf n=7,} $H\cong C_6$ ($|H|=6$, $i=2$):
$$
H_{7,1}=[
(\alpha_{1,1}\alpha_{3,1})(\alpha_{1,2}\alpha_{3,2})
(\alpha_{1,3}\alpha_{3,4}\alpha_{1,6}\alpha_{3,3}\alpha_{1,4}\alpha_{3,6})
(\alpha_{2,3}\alpha_{2,4}\alpha_{2,6})
$$
$$
(\alpha_{1,5}\alpha_{3,8}\alpha_{1,7}\alpha_{3,5}\alpha_{1,8}\alpha_{3,7})
(\alpha_{2,5}\alpha_{2,8}\alpha_{2,7})]
$$
with $Clos(H_{7,1})=H_{18,1}$ above.

\medskip

{\bf n=6,} $H\cong D_6$ ($|H|=6$, $i=1$):
$\rk N_H=14$ and $(N_H)^\ast/N_H\cong (\bz/6\bz)^2\times (\bz/3\bz)^3$.
$$
H_{6,1}=[
(\alpha_{1,2}\alpha_{3,2})(\alpha_{1,3}\alpha_{3,3})(\alpha_{1,4}\alpha_{3,6})
(\alpha_{2,4}\alpha_{2,6})
(\alpha_{3,4}\alpha_{1,6})
$$
$$
(\alpha_{1,5}\alpha_{1,8})(\alpha_{2,5}\alpha_{2,8})(\alpha_{3,5}\alpha_{3,8}),
$$
$$
(\alpha_{1,3}\alpha_{1,4}\alpha_{1,6})(\alpha_{2,3}\alpha_{2,4}\alpha_{2,6})
(\alpha_{3,3}\alpha_{3,4}\alpha_{3,6})
(\alpha_{1,5}\alpha_{1,8}\alpha_{1,7})
$$
$$
(\alpha_{2,5}\alpha_{2,8}\alpha_{2,7})(\alpha_{3,5}\alpha_{3,8}\alpha_{3,7})]
$$
with orbits
$
\{\alpha_{1,2},\alpha_{3,2}\},\
\{\alpha_{1,3},\alpha_{3,3},\alpha_{1,4},\alpha_{3,4},\alpha_{3,6},\alpha_{1,6}\},\
\{\alpha_{2,3},\alpha_{2,4},
\alpha_{2,6}\},\  \{\alpha_{1,5},\alpha_{1,8},\alpha_{1,7}\},\ \linebreak
\{\alpha_{2,5},\alpha_{2,8},\alpha_{2,7}\},\
\{\alpha_{3,5},\alpha_{3,8},\alpha_{3,7}\}
$.

\medskip

{\bf n=4,} $H\cong C_4$ ($|H|=4$, $i=1$):
$\rk N_H=14$ and $(N_H)^\ast/N_H \cong
(\bz/4\bz)^4\times (\bz/2\bz)^2$.
$$
H_{4,1}=[
(\alpha_{1,1}\alpha_{3,1})(\alpha_{1,3}\alpha_{3,7},\alpha_{3,3}\alpha_{1,7})
(\alpha_{2,3}\alpha_{2,7})(\alpha_{1,4}\alpha_{3,6}\alpha_{3,5}\alpha_{3,8})
$$
$$
(\alpha_{2,4}\alpha_{2,6}\alpha_{2,5}\alpha_{2,8})
(\alpha_{3,4}\alpha_{1,6}\alpha_{1,5}\alpha_{1,8})]
$$
with orbits
$
\{\alpha_{1,1},\alpha_{3,1}\},\ \{\alpha_{1,3},\alpha_{3,7},
\alpha_{3,3},\alpha_{1,7}\},\
\{\alpha_{2,3},\alpha_{2,7}\},\
\{\alpha_{1,4},\alpha_{3,6},
\alpha_{3,5},\alpha_{3,8}\},\ \linebreak
\{\alpha_{2,4},\alpha_{2,6},\alpha_{2,5},\alpha_{2,8}\},\
\{\alpha_{3,4},\alpha_{1,6},\alpha_{1,5}, \alpha_{1,8}\}
$.

\medskip

{\bf n=3,} $H\cong C_2^2$ ($|H|=4$, $i=2$):
$\rk N_H=12$ and $(N_H)^\ast/N_H\cong (\bz/4\bz)^2\times (\bz/2\bz)^6$.
$$
H_{3,1}=
[(\alpha_{1,1}\alpha_{3,1})(\alpha_{1,2}\alpha_{3,2})
(\alpha_{1,3}\alpha_{3,3})(\alpha_{1,4}\alpha_{3,4})
(\alpha_{1,5}\alpha_{3,5})(\alpha_{1,6}\alpha_{3,6})
$$
$$
(\alpha_{1,7}\alpha_{3,7})(\alpha_{1,8}\alpha_{3,8}),
$$
$$
(\alpha_{1,2}\alpha_{3,2})(\alpha_{1,3}\alpha_{3,6})
(\alpha_{2,3}\alpha_{2,6})(\alpha_{3,3}\alpha_{1,6})
(\alpha_{1,4}\alpha_{3,4})(\alpha_{1,7}\alpha_{1,8})
$$
$$
(\alpha_{2,7}\alpha_{2,8})(\alpha_{3,7}\alpha_{3,8}) ]
$$
with orbits
$
\{\alpha_{1,1},\alpha_{3,1}\},\  \{\alpha_{1,2},\alpha_{3,2}\},\
\{\alpha_{1,3},\alpha_{3,3},\alpha_{3,6},\alpha_{1,6}\},\
\{\alpha_{2,3},
\alpha_{2,6}\},\
\{\alpha_{1,4},\alpha_{3,4}\},\
\linebreak  \{\alpha_{1,5},\alpha_{3,5}\},\
\{\alpha_{1,7},\alpha_{3,7},\alpha_{1,8},\alpha_{3,8}\},\
\{\alpha_{2,7},\alpha_{2,8}\}
$;
$$
H_{3,2}=[
(\alpha_{1,3}\alpha_{1,7})(\alpha_{2,3}\alpha_{2,7})
(\alpha_{3,3}\alpha_{3,7})(\alpha_{1,4}\alpha_{3,4})
(\alpha_{1,5}\alpha_{3,5})(\alpha_{1,6}\alpha_{3,8})
$$
$$
(\alpha_{2,6}\alpha_{2,8})(\alpha_{3,6}\alpha_{1,8}),
$$
$$
(\alpha_{1,2}\alpha_{3,2})(\alpha_{1,3}\alpha_{3,6})
(\alpha_{2,3}\alpha_{2,6})(\alpha_{3,3}\alpha_{1,6})
(\alpha_{1,4}\alpha_{3,4})(\alpha_{1,7}\alpha_{1,8})
$$
$$
(\alpha_{2,7}\alpha_{2,8})(\alpha_{3,7}\alpha_{3,8})]
$$
with orbits
$
\{\alpha_{1,2},\alpha_{3,2}\},\  \{\alpha_{1,3},\alpha_{1,7},\alpha_{3,6},\alpha_{1,8}\},\
\{\alpha_{2,3},\alpha_{2,7},\alpha_{2,6},\alpha_{2,8}\},\
\{\alpha_{3,3},\alpha_{3,7},\alpha_{1,6},\alpha_{3,8}\},\ \linebreak
\{\alpha_{1,4},\alpha_{3,4}\},\ \{\alpha_{1,5},\alpha_{3,5}\}
$;
$$
H_{3,3}=[
(\alpha_{1,2}\alpha_{1,4})(\alpha_{2,2}\alpha_{2,4})(\alpha_{3,2}\alpha_{3,4})
(\alpha_{1,3}\alpha_{3,3})(\alpha_{1,6}\alpha_{3,6})(\alpha_{1,7}\alpha_{3,8})
$$
$$
(\alpha_{2,7}\alpha_{2,8})(\alpha_{3,7}\alpha_{1,8}),
$$
$$
(\alpha_{1,2}\alpha_{3,2})(\alpha_{1,3}\alpha_{3,6})(\alpha_{2,3}\alpha_{2,6})
(\alpha_{3,3}\alpha_{1,6})(\alpha_{1,4}\alpha_{3,4})(\alpha_{1,7}\alpha_{1,8})
$$
$$
(\alpha_{2,7}\alpha_{2,8})(\alpha_{3,7}\alpha_{3,8})]
$$
with orbits
$
\{\alpha_{1,2},\alpha_{1,4},\alpha_{3,2},\alpha_{3,4}\},\
\{\alpha_{2,2},\alpha_{2,4}\},\
\{\alpha_{1,3},\alpha_{3,3},\alpha_{3,6},\alpha_{1,6}\},\
\{\alpha_{2,3},\alpha_{2,6}\},\  \linebreak
\{\alpha_{1,7},\alpha_{3,8},\alpha_{1,8},\alpha_{3,7}\},\
\{\alpha_{2,7},\alpha_{2,8}\}
$;
$$
H_{3,4}=[
(\alpha_{1,1}\alpha_{3,1})(\alpha_{1,2}\alpha_{1,4})
(\alpha_{2,2}\alpha_{2,4})(\alpha_{3,2}\alpha_{3,4})
(\alpha_{1,3}\alpha_{3,6})(\alpha_{2,3}\alpha_{2,6})
$$
$$
(\alpha_{3,3}\alpha_{1,6})(\alpha_{1,5}\alpha_{3,5}),
$$
$$
(\alpha_{1,2}\alpha_{3,2})(\alpha_{1,3}\alpha_{3,6})
(\alpha_{2,3}\alpha_{2,6})(\alpha_{3,3}\alpha_{1,6})
(\alpha_{1,4}\alpha_{3,4})(\alpha_{1,7}\alpha_{1,8})
$$
$$
(\alpha_{2,7}\alpha_{2,8})(\alpha_{3,7}\alpha_{3,8}) ]
$$
with orbits
$
\{\alpha_{1,1},\alpha_{3,1}\},\ \{\alpha_{1,2},\alpha_{1,4},
\alpha_{3,2},\alpha_{3,4}\},\
\{\alpha_{2,2},\alpha_{2,4}\},\ \{\alpha_{1,3},\alpha_{3,6}\},\
\{\alpha_{2,3},\alpha_{2,6}\},\ \linebreak
\{\alpha_{3,3},\alpha_{1,6}\},\
\{\alpha_{1,5},\alpha_{3,5}\},\
\{\alpha_{1,7},\alpha_{1,8}\},\ \{\alpha_{2,7},\alpha_{2,8}\},\
\{\alpha_{3,7},\alpha_{3,8}\}
$.

\medskip


{\bf n=2,} $H\cong C_3$ ($|H|=3$, $i=1$):
$\rk N_H=12$ and $(N_H)^\ast/N_H\cong (\bz/3\bz)^6$.
$$
H_{2,1}=[
(\alpha_{1,3}\alpha_{1,4}\alpha_{1,6})(\alpha_{2,3}\alpha_{2,4}\alpha_{2,6})
(\alpha_{3,3}\alpha_{3,4}\alpha_{3,6})
(\alpha_{1,5}\alpha_{1,8}\alpha_{1,7})
$$
$$
(\alpha_{2,5}\alpha_{2,8}\alpha_{2,7})(\alpha_{3,5}\alpha_{3,8}\alpha_{3,7})]
$$
with orbits
$
\{\alpha_{1,3},\alpha_{1,4},\alpha_{1,6}\},\
\{\alpha_{2,3},\alpha_{2,4},\alpha_{2,6}\},\
\{\alpha_{3,3},\alpha_{3,4},\alpha_{3,6}\},\
\{\alpha_{1,5},\alpha_{1,8},\alpha_{1,7}\},\ \linebreak
\{\alpha_{2,5},\alpha_{2,8},\alpha_{2,7}\},\
\{\alpha_{3,5},\alpha_{3,8},\alpha_{3,7}\}
$.

\medskip

{\bf n=1,} $H\cong C_2$ ($|H|=2$, $i=1$):
$\rk N_H=8$ and $(N_H)^\ast/N_H\cong (\bz/2\bz)^8$.
$$
H_{1,1}=[
(\alpha_{1,1}\alpha_{3,1})(\alpha_{1,2}\alpha_{3,2})
(\alpha_{1,3}\alpha_{3,3})(\alpha_{1,4}\alpha_{3,4})
(\alpha_{1,5}\alpha_{3,5})(\alpha_{1,6}\alpha_{3,6})
$$
$$
(\alpha_{1,7}\alpha_{3,7})(\alpha_{1,8}\alpha_{3,8})]
$$
with orbits
$
\{\alpha_{1,1},\alpha_{3,1}\},\ \{\alpha_{1,2},\alpha_{3,2}\},\
\{\alpha_{1,3},\alpha_{3,3}\},\
\{\alpha_{1,4},\alpha_{3,4}\},\
\{\alpha_{1,5},\alpha_{3,5}\},\
\{\alpha_{1,6},\alpha_{3,6}\},\ \linebreak
\{\alpha_{1,7},\alpha_{3,7}\},\ \{\alpha_{1,8},\alpha_{3,8}\}
$;
$$
H_{1,2}=[
(\alpha_{1,2}\alpha_{3,2})(\alpha_{1,3}\alpha_{3,6})
(\alpha_{2,3}\alpha_{2,6})(\alpha_{3,3}\alpha_{1,6})
(\alpha_{1,4}\alpha_{3,4})(\alpha_{1,7}\alpha_{1,8})
$$
$$
(\alpha_{2,7}\alpha_{2,8})(\alpha_{3,7}\alpha_{3,8})]
$$
with orbits
$
\{\alpha_{1,2},\alpha_{3,2}\},\ \{\alpha_{1,3},\alpha_{3,6}\},\
\{\alpha_{2,3},\alpha_{2,6}\},\  \{\alpha_{3,3},\alpha_{1,6}\},\
\{\alpha_{1,4},\alpha_{3,4}\},\
\{\alpha_{1,7},\alpha_{1,8}\},\ \linebreak
\{\alpha_{2,7},\alpha_{2,8}\},\ \{\alpha_{3,7},\alpha_{3,8}\}
$.

\vskip1cm

In \cite{Nik7} and \cite{Nik8} these conjugacy classes were described by
the direct considerations. Here is the correspondence between these two
descriptions:
$H_{1,1}\cong {\rm (I.10)}$; $H_{1,2}\cong $ (II.11)
and (II,11');
$H_{2,1}\cong {\rm (II.6)}$;
$H_{3,1}\cong {\rm (I.9)}$;
$H_{3,2}\cong {\rm (II.8),\ (II.8')\ and\ (II.8'')}$;
$H_{3,3}\cong {\rm (II.9)}$;
$H_{3,4}\cong {\rm (II.9')}$;
$H_{4,1}\cong {\rm (II.10)\ and\ (II.10')}$;
$H_{6,1}\cong {\rm (II.5)\ and\ (II.5')}$;
$H_{7,1}\cong {\rm (I.4)}$;
$H_{8,1}\cong {\rm (II.2a)}$;
$H_{9,1}\cong {\rm (I.6)}$;
$H_{9,2}\cong {\rm (I.7)}$;
$H_{10,1}\cong {\rm (II.7)\ and\ (II.7''')}$;
$H_{10,2}\cong {\rm (II.7')\ and\ (II.7'')}$;
$H_{11,1}\cong {\rm (I.8)}$;
$H_{17,1}\cong {\rm (II.4)}$;
$H_{17,2}\cong {\rm (IV.2)\ and\ (IV.2')}$;
$H_{18,1}\cong {\rm (I.3)}$;
$H_{22,1}\cong {\rm (I.5)}$;
$H_{33,1}\cong {\rm (II.2)}$;
$H_{34,1}\cong {\rm (II.3)\ and\ (II.3')}$;
$H_{34,2}\cong {\rm (IV.1')}$;
$H_{34,3}\cong {\rm (IV.1)}$;
$H_{35,1}\cong {\rm (I.2)}$;
$H_{35,2}\cong {\rm (III.2)}$;
$H_{51,1}\cong {\rm (I.1)}$;
$H_{51,2}\cong {\rm (III.1)}$;
$H_{74,1}\cong {\rm (II.1)}$.

\vskip1cm

{\bf Case 20.} For the Niemeier lattice
$$
N=N_{20}=N(6A_4)=[6A_4,[1(01441)]]=
$$
$$
[6A_4,\ \varepsilon_{1,1}+\varepsilon_{1,3}-\varepsilon_{1,4}-
\varepsilon_{1,5}+\varepsilon_{1,6},\,
\varepsilon_{1,1}+\varepsilon_{1,2}+\varepsilon_{1,4}-
\varepsilon_{1,5}-\varepsilon_{1,6},\,
\varepsilon_{1,1}-\varepsilon_{1,2}+\varepsilon_{1,3}+
\varepsilon_{1,5}-\varepsilon_{1,6},\,
$$
$$
\varepsilon_{1,1}-\varepsilon_{1,2}-\varepsilon_{1,3}+
\varepsilon_{1,4}+\varepsilon_{1,6},\,
\varepsilon_{1,1}+\varepsilon_{1,2}-\varepsilon_{1,3}-
\varepsilon_{1,4}+\varepsilon_{1,5}]
$$
the group $A=A(N_{20})$ has the order $240$, and it is
generated by
$$
\widetilde{(12)}=
(\alpha_{1,1}\alpha_{1,5}\alpha_{4,1}\alpha_{4,5})
(\alpha_{2,1}\alpha_{2,5}\alpha_{3,1}\alpha_{3,5})
$$
$$
(\alpha_{1,2}\alpha_{4,6}\alpha_{4,2}\alpha_{1,6})
(\alpha_{2,2}\alpha_{3,6}\alpha_{3,2}\alpha_{2,6})
(\alpha_{1,3}\alpha_{4,4}\alpha_{4,3}\alpha_{1,4})
(\alpha_{2,3}\alpha_{3,4}\alpha_{3,3}\alpha_{2,4}),
$$
$$
\widetilde{(23)}=
(\alpha_{1,1}\alpha_{1,6}\alpha_{4,1}\alpha_{4,6})
(\alpha_{2,1}\alpha_{2,6}\alpha_{3,1}\alpha_{3,6})
$$
$$
(\alpha_{1,2}\alpha_{1,3}\alpha_{4,2}\alpha_{4,3})
(\alpha_{2,2}\alpha_{2,3}\alpha_{3,2}\alpha_{3,3})
(\alpha_{1,4}\alpha_{4,5}\alpha_{4,4}\alpha_{1,5})
(\alpha_{2,4}\alpha_{3,5}\alpha_{3,4}\alpha_{2,5}),
$$
$$
\widetilde{(34)}=
(\alpha_{1,1}\alpha_{1,2}\alpha_{4,1}\alpha_{4,2})
(\alpha_{2,1}\alpha_{2,2}\alpha_{3,1}\alpha_{3,2})
$$
$$
(\alpha_{1,3}\alpha_{1,4}\alpha_{4,3}\alpha_{4,4})
(\alpha_{2,3}\alpha_{2,4}\alpha_{3,3}\alpha_{3,4})
(\alpha_{1,5}\alpha_{4,6}\alpha_{4,5}\alpha_{1,6})
(\alpha_{2,5}\alpha_{3,6}\alpha_{3,5}\alpha_{2,6}),
$$
$$
\widetilde{(45)}=
(\alpha_{1,1}\alpha_{1,3}\alpha_{4,1}\alpha_{4,3})
(\alpha_{2,1}\alpha_{2,3}\alpha_{3,1}\alpha_{3,3})
$$
$$
(\alpha_{1,2}\alpha_{1,6}\alpha_{4,2}\alpha_{4,6})
(\alpha_{2,2}\alpha_{2,6}\alpha_{3,2}\alpha_{3,6})
(\alpha_{1,4}\alpha_{1,5}\alpha_{4,4}\alpha_{4,5})
(\alpha_{2,4}\alpha_{2,5}\alpha_{3,4}\alpha_{3,5}).
$$
(see \cite[Chs. 16, 18]{CS} and \cite{Nik7}, \cite{Nik8}).

\vskip1cm

\centerline {\bf Classification of KahK3 conjugacy classes for $A(N_{20})$:}

\vskip1cm

{\bf n=32,} $H\cong Hol(C_5)$ ($|H|=20$, $i=3$):
$\rk N_H=18$  and
$(N_H)^\ast/N_H \cong (\bz/10\bz)^2 \times \bz/5\bz$.
$$
H_{32,1}=[
(\alpha_{1,2}\alpha_{4,2})(\alpha_{2,2}\alpha_{3,2})
(\alpha_{1,3}\alpha_{4,4}\alpha_{1,6}\alpha_{4,5})
(\alpha_{2,3}\alpha_{3,4}\alpha_{2,6}\alpha_{3,5})
$$
$$
(\alpha_{3,3}\alpha_{2,4}\alpha_{3,6}\alpha_{2,5})
(\alpha_{4,3}\alpha_{1,4}\alpha_{4,6}\alpha_{1,5}),
$$
$$
(\alpha_{1,2}\alpha_{1,3}\alpha_{1,4}\alpha_{1,5}\alpha_{1,6})
(\alpha_{2,2}\alpha_{2,3}\alpha_{2,4}\alpha_{2,5}\alpha_{2,6})
(\alpha_{3,2}\alpha_{3,3}\alpha_{3,4}\alpha_{3,5}\alpha_{3,6})
$$
$$
(\alpha_{4,2}\alpha_{4,3}\alpha_{4,4}\alpha_{4,5}\alpha_{4,6})]
$$
with orbits
$
\{\alpha_{1,2},\alpha_{4,2},\alpha_{1,3},\alpha_{4,3},\alpha_{4,4},\alpha_{1,4},
\alpha_{1,6},\alpha_{4,5},\alpha_{4,6},\alpha_{1,5}\},\
\{\alpha_{2,2}, \alpha_{3,2},\alpha_{2,3},\alpha_{3,3},
\linebreak
\alpha_{3,4},\alpha_{2,4},
\alpha_{2,6},\alpha_{3,5},\alpha_{3,6},\alpha_{2,5}\}
$;
$$
H_{32,2}=[
(\alpha_{1,1}\alpha_{4,1})(\alpha_{2,1}\alpha_{3,1})
(\alpha_{1,3}\alpha_{1,4}\alpha_{1,6}\alpha_{1,5})
(\alpha_{2,3}\alpha_{2,4}\alpha_{2,6}\alpha_{2,5})
$$
$$
(\alpha_{3,3}\alpha_{3,4}\alpha_{3,6}\alpha_{3,5})
(\alpha_{4,3}\alpha_{4,4}\alpha_{4,6}\alpha_{4,5}),
$$
$$
(\alpha_{1,2}\alpha_{1,3}\alpha_{1,4}\alpha_{1,5}\alpha_{1,6})
(\alpha_{2,2}\alpha_{2,3}\alpha_{2,4}\alpha_{2,5}\alpha_{2,6})
(\alpha_{3,2}\alpha_{3,3}\alpha_{3,4}\alpha_{3,5}\alpha_{3,6})
$$
$$
(\alpha_{4,2}\alpha_{4,3}\alpha_{4,4}\alpha_{4,5}\alpha_{4,6})]
$$
with orbits
$
\{\alpha_{1,1},\alpha_{4,1}\},\  \{\alpha_{2,1},\alpha_{3,1}\},\
\{\alpha_{1,2},\alpha_{1,3},\alpha_{1,4},\alpha_{1,6},\alpha_{1,5}\},\
\{\alpha_{2,2},\alpha_{2,3},\alpha_{2,4},\alpha_{2,6},\alpha_{2,5}\},\ \linebreak
\{\alpha_{3,2},\alpha_{3,3},\alpha_{3,4},\alpha_{3,6},\alpha_{3,5}\},\
\{\alpha_{4,2},\alpha_{4,3},\alpha_{4,4},
\alpha_{4,6},\alpha_{4,5}\}
$.

\medskip

{\bf n=16,} $H\cong D_{10}$ ($|H|=10$, $i=1$):
$\rk N_H=16$ and $(N_H)^\ast/N_H\cong (\bz/5\bz)^4$.
$$
H_{16,1}=[
(\alpha_{1,3}\alpha_{1,6})(\alpha_{2,3}\alpha_{2,6})
(\alpha_{3,3}\alpha_{3,6})(\alpha_{4,3}\alpha_{4,6})
(\alpha_{1,4}\alpha_{1,5})(\alpha_{2,4}\alpha_{2,5})
$$
$$
(\alpha_{3,4}\alpha_{3,5})(\alpha_{4,4}\alpha_{4,5}),
$$
$$
(\alpha_{1,2}\alpha_{1,3}\alpha_{1,4}\alpha_{1,5}\alpha_{1,6})
(\alpha_{2,2}\alpha_{2,3}\alpha_{2,4}\alpha_{2,5}\alpha_{2,6})
(\alpha_{3,2}\alpha_{3,3}\alpha_{3,4}\alpha_{3,5}\alpha_{3,6})
$$
$$
(\alpha_{4,2}\alpha_{4,3}\alpha_{4,4}\alpha_{4,5}\alpha_{4,6})]
$$
with orbits
$
\{\alpha_{1,2},\alpha_{1,3},\alpha_{1,6},\alpha_{1,4},\alpha_{1,5}\},\
\{\alpha_{2,2},\alpha_{2,3},\alpha_{2,6},\alpha_{2,4},\alpha_{2,5}\},\
\{\alpha_{3,2},\alpha_{3,3},\alpha_{3,6},\alpha_{3,4},\alpha_{3,5}\},\
\linebreak
\{\alpha_{4,2},\alpha_{4,3},\alpha_{4,6},\alpha_{4,4},\alpha_{4,5}\}
$.

\medskip

{\bf n=5,} $H\cong C_5$ ($|H|=5$, $i=1$):
$$
H_{5,1}=[
(\alpha_{1,2}\alpha_{1,3}\alpha_{1,4}\alpha_{1,5}\alpha_{1,6})
(\alpha_{2,2}\alpha_{2,3}\alpha_{2,4}\alpha_{2,5}\alpha_{2,6})
(\alpha_{3,2}\alpha_{3,3}\alpha_{3,4}\alpha_{3,5}\alpha_{3,6})
$$
$$
(\alpha_{4,2}\alpha_{4,3}\alpha_{4,4}\alpha_{4,5}\alpha_{4,6})]
$$
with $Clos(H_{5,1})=H_{16,1}$ above.

\medskip

{\bf n=4,} $H\cong C_4$ ($|H|=4$, $i=1$):
$\rk N_H=14$ and $(N_H)^\ast/N_H \cong
(\bz/4\bz)^4\times (\bz/2\bz)^2$.
$$
H_{4,1}=[
(\alpha_{1,2}\alpha_{4,2})(\alpha_{2,2}\alpha_{3,2})
(\alpha_{1,3}\alpha_{4,4}\alpha_{1,6}\alpha_{4,5})
(\alpha_{2,3}\alpha_{3,4}\alpha_{2,6}\alpha_{3,5})
$$
$$
(\alpha_{3,3}\alpha_{2,4}\alpha_{3,6}\alpha_{2,5})
(\alpha_{4,3}\alpha_{1,4}\alpha_{4,6}\alpha_{1,5})]
$$
with orbits
$
\{\alpha_{1,2},\alpha_{4,2}\},\  \{\alpha_{2,2},\alpha_{3,2}\},\
\{\alpha_{1,3},\alpha_{4,4},\alpha_{1,6},\alpha_{4,5}\},\
\{\alpha_{2,3},\alpha_{3,4},\alpha_{2,6},\alpha_{3,5}\},\ \linebreak
\{\alpha_{3,3},\alpha_{2,4},\alpha_{3,6},\alpha_{2,5}\},\
\{\alpha_{4,3},\alpha_{1,4},\alpha_{4,6},\alpha_{1,5}\}
$.

\medskip

{\bf n=1,} $H\cong C_2$ ($|H|=2$, $i=1$):
$\rk N_H=8$ and $(N_H)^\ast/N_H\cong (\bz/2\bz)^8$.
$$
H_{1,1}=
[(\alpha_{1,3}\alpha_{1,6})(\alpha_{2,3}\alpha_{2,6})
(\alpha_{3,3}\alpha_{3,6})(\alpha_{4,3}\alpha_{4,6})
(\alpha_{1,4}\alpha_{1,5})(\alpha_{2,4}\alpha_{2,5})
$$
$$
(\alpha_{3,4}\alpha_{3,5})(\alpha_{4,4}\alpha_{4,5})]
$$
with orbits
$
\{\alpha_{1,3},\alpha_{1,6}\},\  \{\alpha_{2,3},\alpha_{2,6}\},\
\{\alpha_{3,3},\alpha_{3,6}\},\  \{\alpha_{4,3},\alpha_{4,6}\},\
\{\alpha_{1,4},\alpha_{1,5}\},\  \{\alpha_{2,4},\alpha_{2,5}\},\
\linebreak
\{\alpha_{3,4},\alpha_{3,5}\},\  \{\alpha_{4,4},\alpha_{4,5}\}
$.

\medskip

In \cite{Nik7} and \cite{Nik8}, these conjugacy classes were described by
the direct considerations. Here is the correspondence between these two
descriptions:
$H_{32,1}=[\widetilde{(12)} \widetilde{(23)} \widetilde{(34)} \widetilde{(45)},
\linebreak
\widetilde{(45)}\widetilde{(23)}\widetilde{(34)}\varphi]$ in 1);
$H_{32,2}=[\widetilde{(12)} \widetilde{(23)} \widetilde{(34)} \widetilde{(45)},
\widetilde{(45)}\widetilde{(23)}\widetilde{(34)}]$ in 1);
$H_{16,1}=[\widetilde{(12)}\widetilde{(23)}\widetilde{(34)}\widetilde{(45)},
\linebreak
\widetilde{(45)}\widetilde{(23)}\widetilde{(34)}\widetilde{(23)}
\widetilde{(45)}\widetilde{(34)}]$ in 2);
$H_{5,1}=[\widetilde{(12)}\widetilde{(23)}\widetilde{(34)}\widetilde{(45)}]$ in 3);
$H_{4,1}\cong [\widetilde{(12)}\widetilde{(23)}
\widetilde{(34)}]$ in 4);
 $H_{1,1}\cong [\widetilde{(12)}\widetilde{(34)}]$ in 5).



\vskip0.5cm

{\bf Case 19.} For the Niemeier lattice
$$
N=N_{19}=N(6D_4)=
$$
$$
[6D_4,\ \varepsilon_{1,3}+\varepsilon_{1,4}+\varepsilon_{1,5}+\varepsilon_{1,6},\,
\varepsilon_{1,2}+\varepsilon_{1,4}+\varepsilon_{2,5}+\varepsilon_{3,6},\,
\varepsilon_{1,1}+\varepsilon_{1,4}+\varepsilon_{3,5}+\varepsilon_{2,6},\,
$$
$$
\varepsilon_{2,3}+\varepsilon_{2,4}+\varepsilon_{2,5}+\varepsilon_{2,6},\,
\varepsilon_{2,2}+\varepsilon_{2,4}+\varepsilon_{3,5}+\varepsilon_{1,6},\,
\varepsilon_{2,1}+\varepsilon_{2,4}+\varepsilon_{1,5}+\varepsilon_{3,6}]
$$
the group $A=A(N)$ has order $2160$ and it is generated by
$$
\varphi =
(\alpha_{1,1}\alpha_{3,1}\alpha_{4,1})
(\alpha_{1,2}\alpha_{3,2}\alpha_{4,2})
(\alpha_{1,3}\alpha_{3,3}\alpha_{4,3})
(\alpha_{1,4}\alpha_{3,4}\alpha_{4,4})
$$
$$
(\alpha_{1,5}\alpha_{3,5}\alpha_{4,5})
(\alpha_{1,6}\alpha_{3,6}\alpha_{4,6}),
$$
$$
\widetilde{(12)}=
(\alpha_{4,1}\alpha_{4,2})(\alpha_{1,1}\alpha_{3,2})
(\alpha_{3,1}\alpha_{1,2})(\alpha_{1,3}\alpha_{3,3})
(\alpha_{1,4}\alpha_{3,4})(\alpha_{1,5}\alpha_{3,5})
$$
$$
(\alpha_{1,6}\alpha_{3,6})(\alpha_{2,1}\alpha_{2,2}),
$$
$$
\widetilde{(23)}=
(\alpha_{1,1}\alpha_{3,1})(\alpha_{1,2}\alpha_{3,3})
(\alpha_{3,2}\alpha_{1,3})(\alpha_{4,2}\alpha_{4,3})
(\alpha_{1,4}\alpha_{3,4})(\alpha_{1,5}\alpha_{4,5})
$$
$$
(\alpha_{3,6}\alpha_{4,6})(\alpha_{2,2}\alpha_{2,3}),
$$
$$
\widetilde{(34)}=
(\alpha_{1,1}\alpha_{3,1})(\alpha_{1,2}\alpha_{3,2})
(\alpha_{4,3}\alpha_{4,4})(\alpha_{1,3}\alpha_{3,4})
(\alpha_{3,3}\alpha_{1,4})(\alpha_{1,5}\alpha_{3,5})
$$
$$
(\alpha_{1,6}\alpha_{3,6})(\alpha_{2,3}\alpha_{2,4}),
$$
$$
\widetilde{(45)}=
(\alpha_{1,1}\alpha_{4,1})(\alpha_{3,2}\alpha_{4,2})
(\alpha_{1,3}\alpha_{3,3})(\alpha_{4,4}\alpha_{4,5})
(\alpha_{1,4}\alpha_{3,5})(\alpha_{3,4}\alpha_{1,5})
$$
$$
(\alpha_{1,6}\alpha_{3,6})(\alpha_{2,4}\alpha_{2,5}),
$$
$$
\widetilde{(56)}=
(\alpha_{1,1}\alpha_{3,1})(\alpha_{1,2}\alpha_{3,2})
(\alpha_{1,3}\alpha_{3,3})(\alpha_{1,4}\alpha_{3,4})
(\alpha_{1,5}\alpha_{3,6})(\alpha_{3,5}\alpha_{1,6})
$$
$$
(\alpha_{4,5}\alpha_{4,6})(\alpha_{2,5}\alpha_{2,6})
$$
(see \cite[Chs. 16, 18]{CS}) and \cite{Nik7}, \cite{Nik8}).

\medskip

\centerline {\bf Classification of KahK3 conjugacy classes for $A(N_{19})$:}

\medskip

{\bf n=70,} $H\cong {\mathfrak S}_5$ ($|H|=120$, $i=34$):
$\rk N_H=19$ and
$(N_H)^\ast/N_H\cong \bz/60\bz \times
\bz/5\bz$.
$$
H_{70,1}=[
(\alpha_{1,1}\alpha_{1,2}\alpha_{3,3}\alpha_{4,5}\alpha_{3,4})
(\alpha_{2,1}\alpha_{2,2}\alpha_{2,3}\alpha_{2,5}\alpha_{2,4})
$$
$$
(\alpha_{3,1}\alpha_{3,2}\alpha_{4,3}\alpha_{1,5}\alpha_{4,4})
(\alpha_{4,1}\alpha_{4,2}\alpha_{1,3}\alpha_{3,5}\alpha_{1,4}),
$$
$$
(\alpha_{1,1}\alpha_{4,2})(\alpha_{2,1}\alpha_{2,2})
(\alpha_{3,1}\alpha_{3,2})
(\alpha_{4,1}\alpha_{1,2})(\alpha_{1,3}\alpha_{4,3})
(\alpha_{1,4}\alpha_{4,4})
$$
$$
(\alpha_{1,5}\alpha_{4,5})
(\alpha_{1,6}\alpha_{4,6})]
$$
with orbits
$
\{\alpha_{1,1},\alpha_{1,2},\alpha_{4,2},\alpha_{3,3},\alpha_{4,1},
\alpha_{1,3},\alpha_{4,5},
\alpha_{3,5},\alpha_{4,3},\alpha_{3,4},\alpha_{1,5},
\alpha_{1,4},\alpha_{4,4},\alpha_{3,1},\alpha_{3,2}\},\ \linebreak
\{\alpha_{2,1},\alpha_{2,2},\alpha_{2,3},\alpha_{2,5},\alpha_{2,4}\},\
\{\alpha_{1,6}, \alpha_{4,6}\}
$.

\medskip

{\bf n=61,} $H\cong {\mathfrak A}_{4,3}$ ($|H|=72$, $i=43$):
$\rk N_H=18$ and
$(N_H)^\ast/N_H\cong (\bz/12\bz)^2 \times \bz/3\bz$.
$$
H_{61,1}=[
(\alpha_{1,1}\alpha_{1,3})(\alpha_{2,1}\alpha_{2,3})(\alpha_{3,1}\alpha_{4,3})
(\alpha_{4,1}\alpha_{3,3})(\alpha_{3,2}\alpha_{4,2})(\alpha_{3,4}\alpha_{4,4})
$$
$$
(\alpha_{1,5}\alpha_{4,5})(\alpha_{1,6}\alpha_{3,6}),
$$
$$
(\alpha_{1,2}\alpha_{1,3}\alpha_{1,4})(\alpha_{2,2}\alpha_{2,3}\alpha_{2,4})
(\alpha_{3,2}\alpha_{3,3}\alpha_{3,4})(\alpha_{4,2}\alpha_{4,3}\alpha_{4,4})
$$
$$
(\alpha_{1,5}\alpha_{3,5}\alpha_{4,5})(\alpha_{1,6}\alpha_{4,6}\alpha_{3,6}),
$$
$$
(\alpha_{1,1}\alpha_{4,2})(\alpha_{2,1}\alpha_{2,2})(\alpha_{3,1}\alpha_{3,2})
(\alpha_{4,1}\alpha_{1,2})(\alpha_{1,3}\alpha_{4,3})(\alpha_{1,4}\alpha_{4,4})
$$
$$
(\alpha_{1,5}\alpha_{4,5})(\alpha_{1,6}\alpha_{4,6})]
$$
with orbits
$
\{\alpha_{1,1},\alpha_{1,3},\alpha_{4,2},\alpha_{1,4},
\alpha_{4,3},\alpha_{3,2},\alpha_{1,2},
\alpha_{4,4},\alpha_{3,1},\alpha_{3,3},\alpha_{4,1},
\alpha_{3,4}\},\
\{\alpha_{2,1},\alpha_{2,3},\alpha_{2,2},\alpha_{2,4}\},\ \linebreak
\{\alpha_{1,5},\alpha_{4,5},\alpha_{3,5}\},\
\{\alpha_{1,6},\alpha_{3,6},\alpha_{4,6}\}
$.

\medskip

{\bf n=55,} $H\cong {\mathfrak A}_5$ ($|H|=60$, $i=5$):
$\rk N_H=18$ and $(N_H)^\ast/N_H\cong \bz/30\bz\times \bz/10\bz$.
$$
H_{55,1}=[
(\alpha_{1,1}\alpha_{1,3}\alpha_{1,6}\alpha_{4,2}\alpha_{4,5})
(\alpha_{2,1}\alpha_{2,3}\alpha_{2,6}\alpha_{2,2}\alpha_{2,5})
$$
$$
(\alpha_{3,1}\alpha_{3,3}\alpha_{3,6}\alpha_{1,2}\alpha_{1,5})
(\alpha_{4,1}\alpha_{4,3}\alpha_{4,6}\alpha_{3,2}\alpha_{3,5}),
$$
$$
(\alpha_{1,1}\alpha_{1,2})(\alpha_{2,1}\alpha_{2,2})
(\alpha_{3,1}\alpha_{3,2})(\alpha_{4,1}\alpha_{4,2})
(\alpha_{1,5}\alpha_{1,6})(\alpha_{2,5}\alpha_{2,6})
$$
$$
(\alpha_{3,5}\alpha_{3,6})(\alpha_{4,5}\alpha_{4,6})]
$$
with orbits
$
\{\alpha_{1,1},\alpha_{1,3},\alpha_{1,2},\alpha_{1,6},\alpha_{1,5},\alpha_{4,2},
\alpha_{3,1},\alpha_{4,5},\alpha_{4,1},\alpha_{3,3},\alpha_{3,2},
\alpha_{4,6},\alpha_{4,3},\alpha_{3,6},\alpha_{3,5}\},\ \linebreak
\{\alpha_{2,1},\alpha_{2,3},\alpha_{2,2},\alpha_{2,6},\alpha_{2,5}\}
$.

\medskip

{\bf n=48,} $H\cong {\mathfrak S}_{3,3}$ ($|H|=36$, $i=10$):
$\rk N_H=18$, $(N_H)^\ast/N_H \cong
\bz/18\bz\times
\bz/6\bz\times (\bz/3\bz)^2$.
$$
H_{48,1}=[
(\alpha_{3,1}\alpha_{4,1})(\alpha_{1,2}\alpha_{3,2})
(\alpha_{1,3}\alpha_{4,3})(\alpha_{1,4}\alpha_{4,5})
(\alpha_{2,4}\alpha_{2,5})(\alpha_{3,4}\alpha_{3,5})
$$
$$
(\alpha_{4,4}\alpha_{1,5})(\alpha_{1,6}\alpha_{4,6}),
$$
$$
(\alpha_{1,1}\alpha_{4,2}\alpha_{3,1}\alpha_{1,2}\alpha_{4,1}\alpha_{3,2})
(\alpha_{2,1}\alpha_{2,2})
(\alpha_{1,3}\alpha_{4,3}\alpha_{3,3})(\alpha_{1,4}\alpha_{4,4}\alpha_{3,4})
$$
$$
(\alpha_{1,5}\alpha_{4,6}\alpha_{3,5}\alpha_{1,6}\alpha_{4,5}\alpha_{3,6})
(\alpha_{2,5}\alpha_{2,6})]
$$
with orbits
$
\{\alpha_{1,1},\alpha_{4,2},\alpha_{3,1},\alpha_{4,1},\alpha_{1,2},\alpha_{3,2}\},\
\{\alpha_{2,1},\alpha_{2,2}\},\  \{\alpha_{1,3},\alpha_{4,3},
\alpha_{3,3}\},\ \linebreak
\{\alpha_{1,4},\alpha_{4,5},\alpha_{4,4},\alpha_{3,6},\alpha_{1,5},
\alpha_{3,4},\alpha_{4,6},\alpha_{3,5},\alpha_{1,6}\},\
\{\alpha_{2,4},\alpha_{2,5},\alpha_{2,6}\}
$.

\medskip

{\bf n=47,} $H\cong C_3\times {\mathfrak A}_4$ ($|H|=36$, $i=11$):
$$
H_{47,1}=[
(\alpha_{1,2}\alpha_{1,3}\alpha_{1,4})(\alpha_{2,2}\alpha_{2,3}\alpha_{2,4})
(\alpha_{3,2}\alpha_{3,3}\alpha_{3,4})(\alpha_{4,2}\alpha_{4,3}\alpha_{4,4})
$$
$$
(\alpha_{1,5}\alpha_{3,5}\alpha_{4,5})(\alpha_{1,6}\alpha_{4,6}\alpha_{3,6}),
$$
$$
(\alpha_{1,1}\alpha_{3,2}\alpha_{4,1}\alpha_{1,2}\alpha_{3,1}\alpha_{4,2})
(\alpha_{2,1}\alpha_{2,2})
(\alpha_{1,3}\alpha_{3,4}\alpha_{4,3}\alpha_{1,4}\alpha_{3,3}\alpha_{4,4})
(\alpha_{2,3}\alpha_{2,4})
$$
$$
(\alpha_{1,5}\alpha_{3,5}\alpha_{4,5})
(\alpha_{1,6}\alpha_{3,6}\alpha_{4,6})]
$$
with $Clos(H_{47,1})=H_{61,1}$ above.

\medskip

{\bf n=36,} $H\cong C_3\rtimes D_8$ ($|H|=24$, $i=8$):
$$
H_{36,1}=[
(\alpha_{3,1}\alpha_{4,1})(\alpha_{1,2}\alpha_{1,3})(\alpha_{2,2}\alpha_{2,3})
(\alpha_{3,2}\alpha_{4,3})(\alpha_{4,2}\alpha_{3,3})(\alpha_{3,4}\alpha_{4,4})
$$
$$
(\alpha_{1,5}\alpha_{3,5})(\alpha_{1,6}\alpha_{4,6}),
$$
$$
(\alpha_{1,1}\alpha_{3,2}\alpha_{4,1}\alpha_{1,2}\alpha_{3,1}\alpha_{4,2})
(\alpha_{2,1}\alpha_{2,2})
(\alpha_{1,3}\alpha_{3,4}\alpha_{4,3}\alpha_{1,4}\alpha_{3,3}\alpha_{4,4})
(\alpha_{2,3}\alpha_{2,4})
$$
$$
(\alpha_{1,5}\alpha_{3,5}\alpha_{4,5})(\alpha_{1,6}\alpha_{3,6}\alpha_{4,6})]
$$
with $Clos(H_{36,1})=H_{61,1}$ above.

\medskip

{\bf n=34}, $H\cong {\mathfrak S}_4$ ($|H|=24$, $i=12$):
$\rk N_H=17$ and $(N_H)^\ast/N_H\cong
(\bz/12\bz)^2\times \bz/4\bz$.
$$
H_{34,1}=[
(\alpha_{1,2}\alpha_{1,5}\alpha_{1,6})(\alpha_{2,2}\alpha_{2,5}\alpha_{2,6})
(\alpha_{3,2}\alpha_{3,5}\alpha_{3,6})(\alpha_{4,2}\alpha_{4,5}\alpha_{4,6})
$$
$$
(\alpha_{1,3}\alpha_{3,3}\alpha_{4,3})(\alpha_{1,4}\alpha_{4,4}\alpha_{3,4}),
$$
$$
(\alpha_{1,1}\alpha_{3,2}\alpha_{1,6}\alpha_{3,5})
(\alpha_{2,1}\alpha_{2,2}\alpha_{2,6}\alpha_{2,5})
(\alpha_{3,1}\alpha_{1,2}\alpha_{3,6}\alpha_{1,5})
$$
$$
(\alpha_{4,1}\alpha_{4,2}\alpha_{4,6}\alpha_{4,5})
(\alpha_{1,3}\alpha_{4,3})(\alpha_{3,4}\alpha_{4,4}) ]
$$
with orbits
$
\{\alpha_{1,1},\alpha_{3,2},\alpha_{3,5},\alpha_{1,6},
\alpha_{3,6},\alpha_{1,2},\alpha_{1,5},\alpha_{3,1}\},\
\{\alpha_{2,1},\alpha_{2,2},\alpha_{2,5},
\alpha_{2,6}\},\
\{\alpha_{4,1},\alpha_{4,2}, \linebreak \alpha_{4,5},\alpha_{4,6}\},\
\{\alpha_{1,3},\alpha_{3,3},\alpha_{4,3}\},\
\{\alpha_{1,4},\alpha_{4,4},\alpha_{3,4}\}
$;
$$
H_{34,2}=[
(\alpha_{1,1}\alpha_{3,1}\alpha_{4,1})(\alpha_{1,2}\alpha_{3,5}\alpha_{4,6})
(\alpha_{2,2}\alpha_{2,5}\alpha_{2,6})(\alpha_{3,2}\alpha_{4,5}\alpha_{1,6})
$$
$$
(\alpha_{4,2}\alpha_{1,5}\alpha_{3,6})(\alpha_{1,3}\alpha_{4,3}\alpha_{3,3}),
$$
$$
(\alpha_{1,1}\alpha_{3,2}\alpha_{1,6}\alpha_{3,5})
(\alpha_{2,1}\alpha_{2,2}\alpha_{2,6}\alpha_{2,5})
(\alpha_{3,1}\alpha_{1,2}\alpha_{3,6}\alpha_{1,5})
$$
$$
(\alpha_{4,1}\alpha_{4,2}\alpha_{4,6}\alpha_{4,5})
(\alpha_{1,3}\alpha_{4,3})(\alpha_{3,4}\alpha_{4,4})]
$$
with orbits
$
\{\alpha_{1,1},\alpha_{3,1},\alpha_{3,2},\alpha_{4,1},\alpha_{1,2},\alpha_{4,5},
\alpha_{1,6},\alpha_{4,2},\alpha_{3,5},\alpha_{3,6},\alpha_{1,5},
\alpha_{4,6}\},\
\{\alpha_{2,1},\alpha_{2,2},\alpha_{2,5},\alpha_{2,6}\},\
\linebreak
\{\alpha_{1,3},\alpha_{4,3},\alpha_{3,3}\},\ \{\alpha_{3,4},\alpha_{4,4}\}
$.

\medskip

{\bf n=32,} $H\cong Hol(C_5)$ ($|H|=20$, $i=3$):
$\rk N_H=18$  and
$(N_H)^\ast/N_H \linebreak \cong (\bz/10\bz)^2 \times \bz/5\bz$.
$$
H_{32,1}=[
(\alpha_{3,1}\alpha_{4,1})(\alpha_{1,2}\alpha_{4,6}\alpha_{1,4}\alpha_{3,5})
(\alpha_{2,2}\alpha_{2,6}\alpha_{2,4}\alpha_{2,5})
(\alpha_{3,2}\alpha_{3,6}\alpha_{3,4}\alpha_{1,5})
$$
$$
(\alpha_{4,2}\alpha_{1,6}\alpha_{4,4}\alpha_{4,5})
(\alpha_{3,3}\alpha_{4,3}),
$$
$$
(\alpha_{1,1}\alpha_{3,2}\alpha_{3,6}\alpha_{1,5}\alpha_{3,4})
(\alpha_{2,1}\alpha_{2,2}\alpha_{2,6}\alpha_{2,5}\alpha_{2,4})
(\alpha_{3,1}\alpha_{4,2}\alpha_{4,6}\alpha_{3,5}\alpha_{4,4})
$$
$$
(\alpha_{4,1}\alpha_{1,2}\alpha_{1,6}\alpha_{4,5}\alpha_{1,4})]
$$
with orbits
$
\{\alpha_{1,1},\alpha_{3,2},\alpha_{3,6},\alpha_{3,4},\alpha_{1,5}\},\
\{\alpha_{2,1},\alpha_{2,2},\alpha_{2,6},\alpha_{2,4},\alpha_{2,5}\},\
\{\alpha_{3,1},\alpha_{4,1},\alpha_{4,2},\alpha_{1,2},\linebreak
\alpha_{1,6},\alpha_{4,6},
\alpha_{4,4},\alpha_{4,5},\alpha_{1,4},\alpha_{3,5}\},\
\{\alpha_{3,3},\alpha_{4,3}\}
$.

\medskip

{\bf n=31,} $H\cong C_3\times D_6$ ($|H|=18$, $i=3$):
$$
H_{31,1}=[
(\alpha_{1,1}\alpha_{3,1}\alpha_{4,1})(\alpha_{1,2}\alpha_{3,2}\alpha_{4,2})
(\alpha_{1,3}\alpha_{3,3}\alpha_{4,3})(\alpha_{1,4}\alpha_{3,4}\alpha_{4,4})
$$
$$
(\alpha_{1,5}\alpha_{3,5}\alpha_{4,5})(\alpha_{1,6}\alpha_{3,6}\alpha_{4,6}),
$$
$$
(\alpha_{1,1}\alpha_{1,2}\alpha_{3,1}\alpha_{4,2}\alpha_{4,1}\alpha_{3,2})
(\alpha_{2,1}\alpha_{2,2})(\alpha_{1,3}\alpha_{4,3})
(\alpha_{1,4}\alpha_{4,6}\alpha_{1,5}\alpha_{4,4}\alpha_{1,6}\alpha_{4,5})
$$
$$
(\alpha_{2,4}\alpha_{2,6}\alpha_{2,5})(\alpha_{3,4}\alpha_{3,6}\alpha_{3,5})]
$$
with $Clos(H_{31,1})=H_{48,1}$ above;
$$
H_{31,2}=[
(\alpha_{1,1}\alpha_{3,1}\alpha_{4,1})(\alpha_{1,2}\alpha_{4,2}\alpha_{3,2})
(\alpha_{1,4}\alpha_{1,5}\alpha_{1,6})(\alpha_{2,4}\alpha_{2,5}\alpha_{2,6})
$$
$$
(\alpha_{3,4}\alpha_{3,5}\alpha_{3,6})(\alpha_{4,4}\alpha_{4,5}\alpha_{4,6}),
$$
$$
(\alpha_{1,1}\alpha_{1,2}\alpha_{3,1}\alpha_{3,2}\alpha_{4,1}\alpha_{4,2})
(\alpha_{2,1}\alpha_{2,2})(\alpha_{1,3}\alpha_{4,3}\alpha_{3,3})
(\alpha_{1,4}\alpha_{4,6}\alpha_{3,4}\alpha_{1,6}\alpha_{4,4}\alpha_{3,6})
$$
$$
(\alpha_{2,4}\alpha_{2,6})(\alpha_{1,5}\alpha_{4,5}\alpha_{3,5})]
$$
with $Clos(H_{31,2})=H_{48,1}$ above.

\medskip

{\bf n=30,} $H\cong {\mathfrak A}_{3,3}$ ($|H|=18$, $i=4$):
$\rk N_H=16$ and
$(N_H)^\ast/N_H\cong \bz/9\bz\times (\bz/3\bz)^4$.
$$
H_{30,1}=[
(\alpha_{3,1}\alpha_{4,1})(\alpha_{1,2}\alpha_{1,3})
(\alpha_{2,2}\alpha_{2,3})(\alpha_{3,2}\alpha_{4,3})
(\alpha_{4,2}\alpha_{3,3})(\alpha_{3,4}\alpha_{4,4})
$$
$$
(\alpha_{1,5}\alpha_{3,5})(\alpha_{1,6}\alpha_{4,6}),
$$
$$
(\alpha_{1,1}\alpha_{3,1}\alpha_{4,1})(\alpha_{1,2}\alpha_{3,2}\alpha_{4,2})
(\alpha_{1,3}\alpha_{3,3}\alpha_{4,3})(\alpha_{1,4}\alpha_{3,4}\alpha_{4,4})
$$
$$
(\alpha_{1,5}\alpha_{3,5}\alpha_{4,5})(\alpha_{1,6}\alpha_{3,6}\alpha_{4,6}),
$$
$$
(\alpha_{1,1}\alpha_{3,2})(\alpha_{2,1}\alpha_{2,2})
(\alpha_{3,1}\alpha_{1,2})(\alpha_{4,1}\alpha_{4,2})
(\alpha_{1,3}\alpha_{3,3})(\alpha_{1,4}\alpha_{3,4})
$$
$$
(\alpha_{1,5}\alpha_{3,5})(\alpha_{1,6}\alpha_{3,6})]
$$
with orbits
$
\{\alpha_{1,1},\alpha_{3,1},\alpha_{3,2},\alpha_{4,1},\alpha_{1,2},
\alpha_{4,3},\alpha_{4,2},\alpha_{1,3},\alpha_{3,3}\},\
\{\alpha_{2,1},\alpha_{2,2},
\alpha_{2,3}\},\  \linebreak
\{\alpha_{1,4},\alpha_{3,4},\alpha_{4,4}\},\
\{\alpha_{1,5},\alpha_{3,5},\alpha_{4,5}\},\
\{\alpha_{1,6},\alpha_{4,6},\alpha_{3,6}\}
$.

\medskip

{\bf n=20,} $H\cong Q_{12}$ ($|H|=12$, $i=1$):
$$
H_{20,1}=[
(\alpha_{1,1}\alpha_{1,3}\alpha_{1,2}\alpha_{1,4})
(\alpha_{2,1}\alpha_{2,3}\alpha_{2,2}\alpha_{2,4})
(\alpha_{3,1}\alpha_{4,3}\alpha_{3,2}\alpha_{4,4})
$$
$$
(\alpha_{4,1}\alpha_{3,3}\alpha_{4,2}\alpha_{3,4})
(\alpha_{3,5}\alpha_{4,5})(\alpha_{3,6}\alpha_{4,6}),
$$
$$
(\alpha_{1,1}\alpha_{3,2}\alpha_{4,1}\alpha_{1,2}\alpha_{3,1}\alpha_{4,2})
(\alpha_{2,1}\alpha_{2,2})
(\alpha_{1,3}\alpha_{3,4}\alpha_{4,3}\alpha_{1,4}\alpha_{3,3}\alpha_{4,4})
(\alpha_{2,3}\alpha_{2,4})
$$
$$
(\alpha_{1,5}\alpha_{3,5}\alpha_{4,5})(\alpha_{1,6}\alpha_{3,6}\alpha_{4,6})]
$$
with $Clos(H_{20,1})=H_{61,1}$ above.

\medskip

{\bf n=19,} $H\cong C_2\times C_6$ ($|H|=12$, $i=5$):
$$
H_{19,1}=[
(\alpha_{1,1}\alpha_{1,3})(\alpha_{2,1}\alpha_{2,3})
(\alpha_{3,1}\alpha_{3,3})(\alpha_{4,1}\alpha_{4,3})
(\alpha_{1,2}\alpha_{1,4})(\alpha_{2,2}\alpha_{2,4})
$$
$$
(\alpha_{3,2}\alpha_{3,4})(\alpha_{4,2}\alpha_{4,4}),
$$
$$
(\alpha_{1,1}\alpha_{3,2}\alpha_{4,1}\alpha_{1,2}
\alpha_{3,1}\alpha_{4,2})(\alpha_{2,1}\alpha_{2,2})
(\alpha_{1,3}\alpha_{3,4}\alpha_{4,3}\alpha_{1,4}\alpha_{3,3}\alpha_{4,4})
(\alpha_{2,3}\alpha_{2,4})
$$
$$
(\alpha_{1,5}\alpha_{3,5}\alpha_{4,5})(\alpha_{1,6}\alpha_{3,6}\alpha_{4,6})]
$$
with $Clos(H_{19,1})=H_{61,1}$ above.

\medskip

{\bf n=18,} $H\cong D_{12}$ ($|H|=12$, $i=4$):
$\rk N_H=16$ and $(N_H)^\ast/N_H\cong (\bz/6\bz)^4$.

$$
H_{18,1}=[
(\alpha_{3,1}\alpha_{4,1})(\alpha_{3,2}\alpha_{4,2})
(\alpha_{3,3}\alpha_{4,3})(\alpha_{3,4}\alpha_{4,4})
(\alpha_{1,5}\alpha_{1,6})(\alpha_{2,5}\alpha_{2,6})
$$
$$
(\alpha_{3,5}\alpha_{4,6})(\alpha_{4,5}\alpha_{3,6}),
$$
$$
(\alpha_{1,1}\alpha_{3,2}\alpha_{4,1}\alpha_{1,2}\alpha_{3,1}\alpha_{4,2})
(\alpha_{2,1}\alpha_{2,2})
(\alpha_{1,3}\alpha_{3,3}\alpha_{4,3})(\alpha_{1,4}\alpha_{3,4}\alpha_{4,4})
$$
$$
(\alpha_{1,5}\alpha_{3,6}\alpha_{4,5}\alpha_{1,6}\alpha_{3,5}\alpha_{4,6})
(\alpha_{2,5}\alpha_{2,6}) ]
$$
with orbits
$
\{\alpha_{1,1},\alpha_{3,2},\alpha_{4,2},\alpha_{4,1},\alpha_{3,1},\alpha_{1,2}\},\
\{\alpha_{2,1},\alpha_{2,2}\},\  \{\alpha_{1,3},\alpha_{3,3},
\alpha_{4,3}\},\ \linebreak
\{\alpha_{1,4},\alpha_{3,4},\alpha_{4,4}\},\
\{\alpha_{1,5},\alpha_{1,6},\alpha_{3,6},\alpha_{3,5},\alpha_{4,5},\alpha_{4,6}\},\
\{\alpha_{2,5},\alpha_{2,6}\}
$;
$$
H_{18,2}=[
(\alpha_{1,1}\alpha_{4,1})(\alpha_{1,2}\alpha_{4,6})
(\alpha_{2,2}\alpha_{2,6})(\alpha_{3,2}\alpha_{3,6})
(\alpha_{4,2}\alpha_{1,6})(\alpha_{1,3}\alpha_{3,3})
$$
$$
(\alpha_{3,4}\alpha_{4,4})(\alpha_{1,5}\alpha_{4,5}),
$$
$$
(\alpha_{1,1}\alpha_{3,2}\alpha_{3,6}\alpha_{4,1}\alpha_{4,2}\alpha_{1,6})
(\alpha_{2,1}\alpha_{2,2}\alpha_{2,6})(\alpha_{3,1}\alpha_{1,2}\alpha_{4,6})
(\alpha_{1,3}\alpha_{3,3})
$$
$$
(\alpha_{1,4}\alpha_{1,5}\alpha_{4,4}\alpha_{3,5}\alpha_{3,4}\alpha_{4,5})
(\alpha_{2,4}\alpha_{2,5})]
$$
with orbits
$
\{\alpha_{1,1},\alpha_{4,1},\alpha_{3,2},\alpha_{4,2},
\alpha_{3,6},\alpha_{1,6}\},\
\{\alpha_{2,1},\alpha_{2,2},\alpha_{2,6}\},\  \{\alpha_{3,1},
\alpha_{1,2},\alpha_{4,6}\},\ \linebreak
\{\alpha_{1,3},\alpha_{3,3}\},\
\{\alpha_{1,4},\alpha_{1,5},\alpha_{4,5},\alpha_{4,4},\alpha_{3,4},\alpha_{3,5}\},\
\{\alpha_{2,4},\alpha_{2,5}\}
$.

\medskip

{\bf n=17,} $H\cong {\mathfrak A}_4$ ($|H|=12$, $i=3$):
$\rk N_H=16$ and $(N_H)^\ast/N_H\cong (\bz/12\bz)^2\times (\bz/2\bz)^2$.
$$
H_{17,1}=[
(\alpha_{1,2}\alpha_{1,5}\alpha_{1,6})(\alpha_{2,2}\alpha_{2,5}\alpha_{2,6})
(\alpha_{3,2}\alpha_{3,5}\alpha_{3,6})(\alpha_{4,2}\alpha_{4,5}\alpha_{4,6})
$$
$$
(\alpha_{9}\alpha_{11}\alpha_{12})(\alpha_{13}\alpha_{16}\alpha_{15}),
$$
$$
(\alpha_{1,1}\alpha_{1,2})(\alpha_{2,1}\alpha_{2,2})
(\alpha_{3,1}\alpha_{3,2})(\alpha_{4,1}\alpha_{4,2})
(\alpha_{1,5}\alpha_{1,6})(\alpha_{2,5}\alpha_{2,6})
$$
$$
(\alpha_{3,5}\alpha_{3,6})(\alpha_{4,5}\alpha_{4,6})]
$$
with orbits
$
\{\alpha_{1,1},\alpha_{1,2},\alpha_{1,5},\alpha_{1,6}\},\
\{\alpha_{2,1},\alpha_{2,2},\alpha_{2,5},\alpha_{2,6}\},\
\{\alpha_{3,1},\alpha_{3,2},
\alpha_{3,5},\alpha_{3,6}\},\ \linebreak
\{\alpha_{4,1}, \alpha_{4,2}, \alpha_{4,5},\alpha_{4,6}\},\
\{\alpha_{1,3},\alpha_{3,3},\alpha_{4,3}\},\
\{\alpha_{1,4},\alpha_{4,4},\alpha_{3,4}\}
$;
$$
H_{17,2}=[
(\alpha_{1,1}\alpha_{3,1}\alpha_{4,1})(\alpha_{1,2}\alpha_{3,5}\alpha_{4,6})
(\alpha_{2,2}\alpha_{2,5}\alpha_{2,6})(\alpha_{3,2}\alpha_{4,5}\alpha_{1,6})
$$
$$
(\alpha_{4,2}\alpha_{1,5}\alpha_{3,6})(\alpha_{1,3}\alpha_{4,3}\alpha_{3,3}),
$$
$$
(\alpha_{1,1}\alpha_{1,2})(\alpha_{2,1}\alpha_{2,2})
(\alpha_{3,1}\alpha_{3,2})(\alpha_{4,1}\alpha_{4,2})
(\alpha_{1,5}\alpha_{1,6})(\alpha_{2,5}\alpha_{2,6})
$$
$$
(\alpha_{3,5}\alpha_{3,6})(\alpha_{4,5}\alpha_{4,6})]
$$
with orbits
$
\{\alpha_{1,1},\alpha_{3,1},\alpha_{1,2},\alpha_{4,1},\alpha_{3,2},\alpha_{3,5},
\alpha_{4,2},\alpha_{4,5},\alpha_{4,6},\alpha_{3,6},\alpha_{1,5},
\alpha_{1,6}\},\ \linebreak
\{\alpha_{2,1},\alpha_{2,2},\alpha_{2,5},\alpha_{2,6}\},\
\{\alpha_{1,3},\alpha_{4,3},\alpha_{3,3}\}
$.

\medskip

{\bf n=16,} $H\cong D_{10}$ ($|H|=10$, $i=1$):
$\rk N_H=16$ and $(N_H)^\ast/N_H\cong (\bz/5\bz)^4$.
$$
H_{16,1}=[
(\alpha_{1,2}\alpha_{1,4})(\alpha_{2,2}\alpha_{2,4})
(\alpha_{3,2}\alpha_{3,4})(\alpha_{4,2}\alpha_{4,4})
(\alpha_{1,5}\alpha_{3,6})(\alpha_{2,5}\alpha_{2,6})
$$
$$
(\alpha_{3,5}\alpha_{4,6})(\alpha_{4,5}\alpha_{1,6}),
$$
$$
(\alpha_{1,1}\alpha_{3,2}\alpha_{3,6}\alpha_{1,5}\alpha_{3,4})
(\alpha_{2,1}\alpha_{2,2}\alpha_{2,6}\alpha_{2,5}\alpha_{2,4})
(\alpha_{3,1}\alpha_{4,2}\alpha_{4,6}\alpha_{3,5}\alpha_{4,4})
$$
$$
(\alpha_{4,1}\alpha_{1,2}\alpha_{1,6}\alpha_{4,5}\alpha_{1,4})]
$$
with orbits
$
\{\alpha_{1,1},\alpha_{3,2},\alpha_{3,4},\alpha_{3,6},\alpha_{1,5}\},\
\{\alpha_{2,1},\alpha_{2,2},\alpha_{2,4},\alpha_{2,6},\alpha_{2,5}\},\
\{\alpha_{3,1},\alpha_{4,2},\alpha_{4,4},\alpha_{4,6},\alpha_{3,5}\},\
\linebreak
\{\alpha_{4,1},\alpha_{1,2},\alpha_{1,4},\alpha_{1,6},\alpha_{4,5}\}
$.

\medskip

{\bf n=15,} $H\cong C_3^2$ ($|H|=9$, $i=2$):
$$
H_{15,1}=[
(\alpha_{1,1}\alpha_{3,1}\alpha_{4,1})(\alpha_{1,2}\alpha_{3,2}\alpha_{4,2})
(\alpha_{1,3}\alpha_{3,3}\alpha_{4,3})(\alpha_{1,4}\alpha_{3,4}\alpha_{4,4})
$$
$$
(\alpha_{1,5}\alpha_{3,5}\alpha_{4,5})(\alpha_{1,6}\alpha_{3,6}\alpha_{4,6}),
$$
$$
(\alpha_{1,1}\alpha_{1,2}\alpha_{1,3})(\alpha_{2,1}\alpha_{2,2}\alpha_{2,3})
(\alpha_{3,1}\alpha_{3,2}\alpha_{3,3})(\alpha_{4,1}\alpha_{4,2}\alpha_{4,3})
$$
$$
(\alpha_{1,5}\alpha_{4,5}\alpha_{3,5})(\alpha_{1,6}\alpha_{3,6}\alpha_{4,6})]
$$
with $Clos(H_{15,1})=H_{30,1}$ above.

\medskip

{\bf n=10,} $H\cong D_8$ ($|H|=8$, $i=3$):
$\rk N_H=15$ and $(N_H)^\ast/N_H\cong
(\bz/4\bz)^5$.
$$
H_{10,1}=[
(\alpha_{1,1}\alpha_{3,5})(\alpha_{2,1}\alpha_{2,5})
(\alpha_{3,1}\alpha_{1,5})(\alpha_{4,1}\alpha_{4,5})
(\alpha_{1,2}\alpha_{3,2})(\alpha_{3,3}\alpha_{4,3})
$$
$$
(\alpha_{1,4}\alpha_{4,4})(\alpha_{1,6}\alpha_{3,6}),
$$
$$
(\alpha_{1,1}\alpha_{1,2})(\alpha_{2,1}\alpha_{2,2})
(\alpha_{3,1}\alpha_{3,2})(\alpha_{4,1}\alpha_{4,2})
(\alpha_{1,5}\alpha_{1,6})(\alpha_{2,5}\alpha_{2,6})
$$
$$
(\alpha_{3,5}\alpha_{3,6})(\alpha_{4,5}\alpha_{4,6})]
$$
with orbits
$
\{\alpha_{1,1},\alpha_{3,5},\alpha_{1,2},\alpha_{3,6},\alpha_{3,2},
\alpha_{1,6},\alpha_{3,1},\alpha_{1,5}\},\
\{\alpha_{2,1},\alpha_{2,5},\alpha_{2,2},
\alpha_{2,6} \},\ \{\alpha_{4,1},\alpha_{4,5}, \linebreak
\alpha_{4,2},\alpha_{4,6}\},\
\{\alpha_{3,3},\alpha_{4,3}\},\  \{\alpha_{1,4},\alpha_{4,4}\}
$.

\medskip

{\bf n=7,} $H\cong C_6$ ($|H|=6$, $i=2$):
$$
H_{7,1}=[
(\alpha_{1,1}\alpha_{3,2}\alpha_{4,1}\alpha_{1,2}\alpha_{3,1}\alpha_{4,2})
(\alpha_{2,1}\alpha_{2,2})
(\alpha_{1,3}\alpha_{3,3}\alpha_{4,3})(\alpha_{1,4}\alpha_{3,4}\alpha_{4,4})
$$
$$
(\alpha_{1,5}\alpha_{3,6}\alpha_{4,5}\alpha_{1,6}\alpha_{3,5}\alpha_{4,6})
(\alpha_{2,5}\alpha_{2,6}) ]
$$
with $Clos(H_{7,1})=H_{18,1}$ above;
$$
H_{7,2}=[
(\alpha_{1,1}\alpha_{3,2}\alpha_{3,6}\alpha_{4,1}\alpha_{4,2}\alpha_{1,6})
(\alpha_{2,1}\alpha_{2,2}\alpha_{2,6})(\alpha_{3,1}\alpha_{1,2}\alpha_{4,6})
(\alpha_{1,3}\alpha_{3,3})
$$
$$
(\alpha_{1,4}\alpha_{1,5}\alpha_{4,4}\alpha_{3,5}\alpha_{3,4}\alpha_{4,5})
(\alpha_{2,4}\alpha_{2,5})]
$$
with $Clos(H_{7,2})=H_{18,2}$ above.

\medskip

{\bf n=6,} $H\cong D_6$ ($|H|=6$, $i=1$):
$\rk N_H=14$ and $(N_H)^\ast/N_H\cong (\bz/6\bz)^2\times (\bz/3\bz)^3$.
$$
H_{6,1}=[
(\alpha_{1,1}\alpha_{4,2})(\alpha_{2,1}\alpha_{2,2})
(\alpha_{3,1}\alpha_{3,2})(\alpha_{4,1}\alpha_{1,2})
(\alpha_{1,3}\alpha_{4,3})(\alpha_{1,4}\alpha_{4,4})
$$
$$
(\alpha_{1,5}\alpha_{4,5})(\alpha_{1,6}\alpha_{4,6}),
$$
$$
(\alpha_{1,1}\alpha_{3,2})(\alpha_{2,1}\alpha_{2,2})
(\alpha_{3,1}\alpha_{1,2})(\alpha_{4,1}\alpha_{4,2})
(\alpha_{9}\alpha_{11})(\alpha_{13}\alpha_{15})
$$
$$
(\alpha_{1,5}\alpha_{3,5})(\alpha_{1,6}\alpha_{3,6})]
$$
with orbits
$
\{\alpha_{1,1},\alpha_{4,2},\alpha_{3,2},\alpha_{4,1},\alpha_{3,1},\alpha_{1,2}\},\
\{\alpha_{2,1},\alpha_{2,2}\},\  \{\alpha_{1,3},\alpha_{4,3},
\alpha_{3,3}\},\  \{\alpha_{1,4},\alpha_{4,4},\alpha_{3,4}\},\ \linebreak
\{\alpha_{1,5},\alpha_{4,5},\alpha_{3,5}\},\  \{\alpha_{1,6},\alpha_{4,6},\alpha_{3,6}\}
$;
$$
H_{6,2}=[
(\alpha_{1,1}\alpha_{3,3})(\alpha_{2,1}\alpha_{2,3})
(\alpha_{3,1}\alpha_{1,3})(\alpha_{4,1}\alpha_{4,3})
(\alpha_{1,2}\alpha_{3,2})(\alpha_{1,4}\alpha_{3,4})
$$
$$
(\alpha_{3,5}\alpha_{4,5})(\alpha_{1,6}\alpha_{4,6}),
$$
$$
(\alpha_{1,1},\alpha_{3,2})(\alpha_{2,1},\alpha_{2,2})
(\alpha_{3,1},\alpha_{1,2})(\alpha_{4,1},\alpha_{4,2})
(\alpha_{1,3}\alpha_{3,3})(\alpha_{1,4}\alpha_{3,4})
$$
$$
(\alpha_{1,5}\alpha_{3,5})(\alpha_{1,6}\alpha_{3,6}) ])
$$
with orbits
$
\{\alpha_{1,1},\alpha_{3,3},\alpha_{3,2},\alpha_{1,3},\alpha_{1,2},\alpha_{3,1}\},\
\{\alpha_{2,1},\alpha_{2,3},\alpha_{2,2}\},\  \{\alpha_{4,1},
\alpha_{4,3},\alpha_{4,2}\},\ \linebreak
\{\alpha_{1,4},\alpha_{3,4}\},\  \{\alpha_{1,5},\alpha_{3,5},\alpha_{4,5}\},\
\{\alpha_{1,6},\alpha_{4,6},\alpha_{3,6}\}
$;
$$
H_{6,3}=[
(\alpha_{1,1}\alpha_{3,2})(\alpha_{2,1}\alpha_{2,2})
(\alpha_{3,1}\alpha_{4,2})(\alpha_{4,1}\alpha_{1,2})
(\alpha_{1,4}\alpha_{1,6})(\alpha_{2,4}\alpha_{2,6})
$$
$$
(\alpha_{3,4}\alpha_{3,6})(\alpha_{4,4}\alpha_{4,6}),
$$
$$
(\alpha_{1,1}\alpha_{1,2})(\alpha_{2,1}\alpha_{2,2})
(\alpha_{3,1}\alpha_{3,2})(\alpha_{4,1}\alpha_{4,2})
(\alpha_{1,5}\alpha_{1,6})(\alpha_{2,5}\alpha_{2,6})
$$
$$
(\alpha_{3,5}\alpha_{3,6})(\alpha_{4,5}\alpha_{4,6})]
$$
with orbits
$
\{\alpha_{1,1},\alpha_{3,2},\alpha_{1,2},\alpha_{3,1},\alpha_{4,1},\alpha_{4,2}\},\
\{\alpha_{2,1},\alpha_{2,2}\},\  \{\alpha_{1,4},\alpha_{1,6},
\alpha_{1,5}\},\  \{\alpha_{2,4},\alpha_{2,6},\alpha_{2,5}\},\ \linebreak
\{\alpha_{3,4},\alpha_{3,6},\alpha_{3,5}\},\ \{\alpha_{4,4},\alpha_{4,6},\alpha_{4,5}\}
$.

\medskip

{\bf n=5,} $H\cong C_5$ ($|H|=5$, $i=1$):
$$
H_{5,1}=[
(\alpha_{1,1}\alpha_{3,2}\alpha_{3,6}\alpha_{1,5}\alpha_{3,4})
(\alpha_{2,1}\alpha_{2,2}\alpha_{2,6}\alpha_{2,5}\alpha_{2,4})
(\alpha_{3,1}\alpha_{4,2}\alpha_{4,6}\alpha_{3,5}\alpha_{4,4})
$$
$$
(\alpha_{4,1}\alpha_{1,2}\alpha_{1,6}\alpha_{4,5}\alpha_{1,4})]
$$
with $Clos(H_{5,1})=H_{16,1}$ above.

\medskip

{\bf n=4,} $H\cong C_4$ ($|H|=4$, $i=1$):
$\rk N_H=14$ and $(N_H)^\ast/N_H \cong
(\bz/4\bz)^4\times (\bz/2\bz)^2$.
$$
H_{4,1}=[
(\alpha_{1,1}\alpha_{3,2}\alpha_{1,6}\alpha_{3,5})
(\alpha_{2,1}\alpha_{2,2}\alpha_{2,6}\alpha_{2,5})
(\alpha_{3,1}\alpha_{1,2}\alpha_{3,6}\alpha_{1,5})
$$
$$
(\alpha_{4,1}\alpha_{4,2}\alpha_{4,6}\alpha_{4,5})
(\alpha_{1,3}\alpha_{4,3})(\alpha_{3,4}\alpha_{4,4})]
$$
with orbits
$
\{\alpha_{1,1},\alpha_{3,2},\alpha_{1,6},\alpha_{3,5}\},\
\{\alpha_{2,1},\alpha_{2,2},\alpha_{2,6},\alpha_{2,5}\},\
\{\alpha_{3,1},\alpha_{1,2},
\alpha_{3,6},\alpha_{1,5}\},\ \linebreak
\{\alpha_{4,1},\alpha_{4,2},\alpha_{4,6},\alpha_{4,5}\},\
\{\alpha_{1,3},\alpha_{4,3}\},\  \{\alpha_{3,4},\alpha_{4,4}\}
$.

\medskip

{\bf n=3,} $H\cong C_2^2$ ($|H|=4$, $i=2$):
$\rk N_H=12$ and $(N_H)^\ast/N_H\cong (\bz/4\bz)^2\times (\bz/2\bz)^6$.
$$
H_{3,1}=[
(\alpha_{1,1}\alpha_{3,1})(\alpha_{1,2}\alpha_{3,2})
(\alpha_{1,3}\alpha_{3,4})(\alpha_{2,3}\alpha_{2,4})
(\alpha_{3,3}\alpha_{1,4})(\alpha_{4,3}\alpha_{4,4})
$$
$$
(\alpha_{1,5}\alpha_{3,5})(\alpha_{1,6}\alpha_{3,6}),
$$
$$
(\alpha_{1,1}\alpha_{3,2})(\alpha_{2,1}\alpha_{2,2})
(\alpha_{3,1}\alpha_{1,2})(\alpha_{4,1}\alpha_{4,2})
(\alpha_{1,3}\alpha_{3,3})(\alpha_{1,4}\alpha_{3,4})
$$
$$
(\alpha_{1,5}\alpha_{3,5})(\alpha_{1,6}\alpha_{3,6})]
$$
with orbits
$
\{\alpha_{1,1},\alpha_{3,1},\alpha_{3,2},\alpha_{1,2}\},\
\{\alpha_{2,1},\alpha_{2,2}\},\  \{\alpha_{4,1},\alpha_{4,2}\},\
\{\alpha_{1,3},\alpha_{3,4},
\alpha_{3,3},\alpha_{1,4}\},\ \linebreak
\{\alpha_{2,3},\alpha_{2,4}\},\  \{\alpha_{4,3},\alpha_{4,4}\},\
\{\alpha_{1,5},\alpha_{3,5}\},\  \{\alpha_{1,6},\alpha_{3,6}\}
$;
$$
H_{3,2}=[
(\alpha_{1,1}\alpha_{1,5})(\alpha_{2,1}\alpha_{2,5})
(\alpha_{3,1}\alpha_{3,5})(\alpha_{4,1}\alpha_{4,5})
(\alpha_{1,2}\alpha_{1,6})(\alpha_{2,2}\alpha_{2,6})
$$
$$
(\alpha_{3,2}\alpha_{3,6})(\alpha_{4,2}\alpha_{4,6}),
$$
$$
(\alpha_{1,1}\alpha_{1,2})(\alpha_{2,1}\alpha_{2,2})
(\alpha_{3,1}\alpha_{3,2})(\alpha_{4,1}\alpha_{4,2})
(\alpha_{1,5}\alpha_{1,6})(\alpha_{2,5}\alpha_{2,6})
$$
$$
(\alpha_{3,5}\alpha_{3,6})(\alpha_{4,5}\alpha_{4,6})]
$$
with orbits
$
\{\alpha_{1,1},\alpha_{1,5},\alpha_{1,2},\alpha_{1,6}\},\
\{\alpha_{2,1},\alpha_{2,5},\alpha_{2,2},\alpha_{2,6}\},\
\{\alpha_{3,1},\alpha_{3,5},
\alpha_{3,2},\alpha_{3,6}\},\ \linebreak
\{\alpha_{4,1},\alpha_{4,5},\alpha_{4,2},\alpha_{4,6}\}
$;
$$
H_{3,3}=[
(\alpha_{1,1}\alpha_{1,2})(\alpha_{2,1}\alpha_{2,2})
(\alpha_{3,1}\alpha_{3,2})(\alpha_{4,1}\alpha_{4,2})
(\alpha_{1,3}\alpha_{1,4})(\alpha_{2,3}\alpha_{2,4})
$$
$$
(\alpha_{3,3}\alpha_{3,4})(\alpha_{4,3}\alpha_{4,4}),
$$
$$
(\alpha_{1,1}\alpha_{1,2})(\alpha_{2,1}\alpha_{2,2})
(\alpha_{3,1}\alpha_{3,2})(\alpha_{4,1}\alpha_{4,2})
(\alpha_{1,5}\alpha_{1,6})(\alpha_{2,5}\alpha_{2,6})
$$
$$
(\alpha_{3,5}\alpha_{3,6})(\alpha_{4,5}\alpha_{4,6})]
$$
with orbits
$
\{\alpha_{1,1},\alpha_{1,2}\},\  \{\alpha_{2,1},\alpha_{2,2}\},\
\{\alpha_{3,1},\alpha_{3,2}\},\ \{\alpha_{4,1},\alpha_{4,2}\},\
\{\alpha_{1,3},
\alpha_{1,4}\},\  \{\alpha_{2,3},\alpha_{2,4}\},\ \linebreak
\{\alpha_{3,3},\alpha_{3,4}\},\  \{\alpha_{4,3},\alpha_{4,4}\},\
\{\alpha_{1,5},\alpha_{1,6}\},\  \{\alpha_{2,5},\alpha_{2,6}\},\
\{\alpha_{3,5},\alpha_{3,6}\},\ \{\alpha_{4,5},\alpha_{4,6}\}
$.

\medskip

{\bf n=2,} $H\cong C_3$ ($|H|=3$, $i=1$):
$\rk N_H=12$ and $(N_H)^\ast/N_H\cong (\bz/3\bz)^6$.
$$
H_{2,1}=[
(\alpha_{1,1}\alpha_{3,1}\alpha_{4,1})(\alpha_{1,2}\alpha_{3,2}\alpha_{4,2})
(\alpha_{1,3}\alpha_{3,3}\alpha_{4,3})
(\alpha_{1,4}\alpha_{3,4}\alpha_{4,4})
$$
$$
(\alpha_{1,5}\alpha_{3,5}\alpha_{4,5})(\alpha_{1,6}\alpha_{3,6}\alpha_{4,6})]
$$
with orbits
$
\{\alpha_{1,1},\alpha_{3,1},\alpha_{4,1}\},\
\{\alpha_{1,2},\alpha_{3,2},\alpha_{4,2}\},\
\{\alpha_{1,3},\alpha_{3,3},\alpha_{4,3}\},\  \{\alpha_{1,4},
\alpha_{3,4},\alpha_{4,4}\},\ \linebreak
\{\alpha_{1,5},\alpha_{3,5},\alpha_{4,5}\},\
\{\alpha_{1,6},\alpha_{3,6},\alpha_{4,6}\}
$;
$$
H_{2,2}=[
(\alpha_{1,1}\alpha_{3,2}\alpha_{4,6})(\alpha_{2,1}\alpha_{2,2}\alpha_{2,6})
(\alpha_{3,1}\alpha_{4,2}\alpha_{1,6})(\alpha_{4,1}\alpha_{1,2}\alpha_{3,6})
$$
$$
(\alpha_{1,3}\alpha_{4,3}\alpha_{3,3})(\alpha_{1,5}\alpha_{3,5}\alpha_{4,5})]
$$
with orbits
$
\{\alpha_{1,1},\alpha_{3,2},\alpha_{4,6}\},\
\{\alpha_{2,1},\alpha_{2,2},\alpha_{2,6}\},\
\{\alpha_{3,1},\alpha_{4,2},\alpha_{1,6}\},\
\{\alpha_{4,1},\alpha_{1,2},\alpha_{3,6}\},\ \linebreak
\{\alpha_{1,3},\alpha_{4,3},\alpha_{3,3}\},\
\{\alpha_{1,5},\alpha_{3,5},\alpha_{4,5}\}
$.

\medskip

{\bf n=1,} $H\cong C_2$ ($|H|=2$, $i=1$):
$\rk N_H=8$ and $(N_H)^\ast/N_H\cong (\bz/2\bz)^8$.
$$
H_{1,1}=[
(\alpha_{1,1}\alpha_{3,2})(\alpha_{2,1}\alpha_{2,2})
(\alpha_{3,1}\alpha_{1,2})(\alpha_{4,1}\alpha_{4,2})
(\alpha_{1,3}\alpha_{3,3})(\alpha_{1,4}\alpha_{3,4})
$$
$$
(\alpha_{1,5}\alpha_{3,5})(\alpha_{1,6}\alpha_{3,6})]
$$
with orbits
$
\{\alpha_{1,1},\alpha_{3,2}\},\ \{\alpha_{2,1},\alpha_{2,2}\},\
\{\alpha_{3,1},\alpha_{1,2}\},\ \{\alpha_{4,1},\alpha_{4,2}\},\
\{\alpha_{1,3},
\alpha_{3,3}\},\ \{\alpha_{1,4},\alpha_{3,4}\},\ \linebreak
\{\alpha_{1,5},\alpha_{3,5}\},\ \{\alpha_{1,6},\alpha_{3,6}\}
$;
$$
H_{1,2}=[
(\alpha_{1,1}\alpha_{1,2})(\alpha_{2,1}\alpha_{2,2})
(\alpha_{3,1}\alpha_{3,2})(\alpha_{4,1}\alpha_{4,2})
(\alpha_{1,5}\alpha_{1,6})(\alpha_{2,5}\alpha_{2,6})
$$
$$
(\alpha_{3,5}\alpha_{3,6})(\alpha_{4,5}\alpha_{4,6})]
$$
with orbits
$
\{\alpha_{1,1},\alpha_{1,2}\},\ \{\alpha_{2,1},\alpha_{2,2}\},\
\{\alpha_{3,1},\alpha_{3,2}\},\  \{\alpha_{4,1},\alpha_{4,2}\},\
\{\alpha_{1,5},
\alpha_{1,6}\},\  \{\alpha_{2,5},\alpha_{2,6}\},\ \linebreak
\{\alpha_{3,5},\alpha_{3,6}\},\  \{\alpha_{4,5},\alpha_{4,6}\}
$;

\medskip

In \cite{Nik7} and \cite{Nik8}, these conjugacy classes were described by
direct considerations. Here is the correspondence between these two
descriptions:
$H_{70,1}\cong {\rm (II.1)}$;
$H_{61,1}\cong {\rm (I.1)}$;
$H_{55,1}\cong {\rm (II.2)}$;
$H_{48,1}\cong {\rm (I'.1)}$;
$H_{47,1}\cong {\rm (I.2)}$;
$H_{36,1}\cong {\rm (I.3)}$;
$H_{34,1}\cong {\rm (II.6')}$;
$H_{34,2}\cong {\rm (II.6)}$;
$H_{32,1}\cong {\rm (II.3)}$;
$H_{31,1}\cong {\rm (I'.2)}$;
$H_{31,2}\cong {\rm (I'.3)}=
[\varphi, \widetilde{(45)}\widetilde{(23)},
\widetilde{(45)}\widetilde{(56)}]$;
$H_{30,1}\cong {\rm (I.6)}$;
$H_{20,1}\cong {\rm (I.5)}$;
$H_{19,1}\cong {\rm (I.4)}$;
$H_{18,1}\cong {\rm (I.8)}$;
$H_{18,2}\cong {\rm (II.8)}$;
$H_{17,1}\cong {\rm (II.7')}$;
$H_{17,2}\cong {\rm (II.7)}$;
$H_{16,1}\cong {\rm (II.4)}$;
$H_{15,1}\cong {\rm (I.7)}$;
$H_{10,1}\cong {\rm (II.12)}$;
$H_{7,1}\cong {\rm (I.9)}$;
$H_{7,2}\cong {\rm (II.9)}$;
$H_{6,1}\cong {\rm (I.10)}$;
$H_{6,2}\cong {\rm (II.11)}$;
$H_{6,3}\cong {\rm (II.10)}$;
$H_{5,1}\cong {\rm (II.5)}$;
$H_{4,1}\cong {\rm (II.13)}$;
$H_{3,1}\cong {\rm (II.16)}$;
$H_{3,2}\cong {\rm (II.15)}$;
$H_{3,3}\cong {\rm (II'.1)}$;
$H_{2,1}\cong {\rm (I.11)}$;
$H_{2,2}\cong {\rm (II.14)}$;
$H_{1,1}\cong {\rm (II.18)}$;
$H_{1,2}\cong {\rm (II.17)}$.


\vskip1cm


{\bf Case 18.} For the Niemeier lattice
$$
N=N_{18}=N(4A_5\oplus D_4)=
[4A_5\oplus D_4,\
2\varepsilon_{1,1}+2\varepsilon_{1,3}+4\varepsilon_{1,4},\
2\varepsilon_{1,1}+4\varepsilon_{1,2}+2\varepsilon_{1,4},\,
$$
$$
2\varepsilon_{1,1}+2\varepsilon_{1,2}+4\varepsilon_{1,3},\,
3\varepsilon_{1,1}+3\varepsilon_{1,2}+\varepsilon_{1,5},\
3\varepsilon_{1,1}+3\varepsilon_{1,3}+\varepsilon_{2,5},\,
3\varepsilon_{1,1}+3\varepsilon_{1,4}+\varepsilon_{3,5}]
$$
the group $A=A(N)$ has the order $48$, and it is
generated by involutions
$$
\varphi_0=(\alpha_{1,1}\alpha_{5,1})(\alpha_{1,2}\alpha_{5,2})
(\alpha_{1,3}\alpha_{5,3})(\alpha_{1,4}\alpha_{5,4})
(\alpha_{2,1}\alpha_{4,1})(\alpha_{2,2}\alpha_{4,2})
$$
$$
(\alpha_{2,3}\alpha_{4,3})(\alpha_{2,4}\alpha_{4,4}),
$$
$$
\widetilde{(12)}=(\alpha_{1,1}\alpha_{1,2})(\alpha_{5,1}\alpha_{5,2})
(\alpha_{1,4}\alpha_{5,4})(\alpha_{1,5}\alpha_{3,5})
(\alpha_{2,1}\alpha_{2,2})(\alpha_{4,1}\alpha_{4,2})
$$
$$
(\alpha_{2,4}\alpha_{4,4})(\alpha_{3,1}\alpha_{3,2}),
$$
$$
\widetilde{(23)}=
(\alpha_{1,1}\alpha_{5,1})(\alpha_{1,2}\alpha_{1,3})
(\alpha_{5,2}\alpha_{5,3})(\alpha_{1,5}\alpha_{4,5})
(\alpha_{2,1}\alpha_{4,1})(\alpha_{2,2}\alpha_{2,3})
$$
$$
(\alpha_{4,2}\alpha_{4,3})(\alpha_{3,2}\alpha_{3,3}),
$$
$$
\widetilde{(34)}=
(\alpha_{1,1}\alpha_{5,1})(\alpha_{1,3}\alpha_{1,4})
(\alpha_{5,3}\alpha_{5,4})(\alpha_{1,5}\alpha_{3,5})
(\alpha_{2,1}\alpha_{4,1})(\alpha_{2,3}\alpha_{2,4})
$$
$$
(\alpha_{4,3}\alpha_{4,4})(\alpha_{3,3}\alpha_{3,4}).
$$
(see \cite[Ch. 16]{CS} and \cite{Nik7}, \cite{Nik8}).

\medskip

\centerline {\bf Classification of KahK3 conjugacy classes for $A(N_{18})$:}

\medskip

{\bf n=54,} $H\cong T_{48}$ ($|H|=48$, $i=29$):
$\rk N_H=19$ and $(N_H)^\ast/N_H\cong
\bz/24\bz\times \bz/8\bz\times \bz/2\bz$.
$$
H_{54,1}=[
(\alpha_{1,2}\alpha_{1,4}\alpha_{1,3})(\alpha_{2,2}\alpha_{2,4}\alpha_{2,3})
(\alpha_{3,2}\alpha_{3,4}\alpha_{3,3})(\alpha_{4,2}\alpha_{4,4}\alpha_{4,3})
$$
$$
(\alpha_{5,2}\alpha_{5,4}\alpha_{5,3})(\alpha_{1,5}\alpha_{4,5}\alpha_{3,5}),
$$
$$
(\alpha_{1,1}\alpha_{1,2}\alpha_{1,3}\alpha_{5,4}
\alpha_{5,1}\alpha_{5,2}\alpha_{5,3}\alpha_{1,4})
(\alpha_{2,1}\alpha_{2,2}\alpha_{2,3}\alpha_{4,4}
\alpha_{4,1}\alpha_{4,2}\alpha_{4,3}\alpha_{2,4})
$$
$$
(\alpha_{3,1}\alpha_{3,2}\alpha_{3,3}\alpha_{3,4})
(\alpha_{3,5}\alpha_{4,5})]
$$
with orbits
$
\{\alpha_{1,1},\alpha_{1,2},\alpha_{1,4},\alpha_{1,3},
\alpha_{5,4},\alpha_{5,3},\alpha_{5,1},\alpha_{5,2}\},\
\{\alpha_{2,1},\alpha_{2,2},\alpha_{2,4},
\alpha_{2,3},\alpha_{4,4},\alpha_{4,3},
\linebreak
\alpha_{4,1},\alpha_{4,2}\},\
\{\alpha_{3,1},\alpha_{3,2},\alpha_{3,4},\alpha_{3,3}\},\
\{\alpha_{1,5},\alpha_{4,5},\alpha_{3,5}\}
$.

\medskip

{\bf n=38,} $H\cong T_{24}$ ($|H|=24$, $i=3$):
$$
H_{38,1}=[
(\alpha_{1,2}\alpha_{1,4}\alpha_{1,3})(\alpha_{2,2}\alpha_{2,4}\alpha_{2,3})
(\alpha_{3,2}\alpha_{3,4}\alpha_{3,3})(\alpha_{4,2}\alpha_{4,4}\alpha_{4,3})
$$
$$
(\alpha_{5,2}\alpha_{5,4}\alpha_{5,3})(\alpha_{1,5}\alpha_{4,5}\alpha_{3,5}),
$$
$$
(\alpha_{1,1}\alpha_{1,2}\alpha_{5,1}\alpha_{5,2})
(\alpha_{2,1}\alpha_{2,2}\alpha_{4,1}\alpha_{4,2})
(\alpha_{3,1}\alpha_{3,2})(\alpha_{1,3}\alpha_{1,4}
\alpha_{5,3}\alpha_{5,4})
$$
$$
(\alpha_{2,3}\alpha_{2,4}\alpha_{4,3}\alpha_{4,4})
(\alpha_{3,3}\alpha_{3,4})]
$$
with $Clos(H_{38,1})=H_{54,1}$ above.

\medskip


{\bf n=26,} $H\cong SD_{16}$ ($|H|=16$, $i=8$):
$\rk N_H=18$ and
$(N_H)^\ast/N_H\cong (\bz/8\bz)^2\times
\bz/4\bz\times \bz/2\bz$.
$$
H_{26,1}=[
(\alpha_{1,1}\alpha_{1,3}\alpha_{5,1}\alpha_{5,3})
(\alpha_{2,1}\alpha_{2,3}\alpha_{4,1}\alpha_{4,3})
(\alpha_{3,1}\alpha_{3,3})
(\alpha_{1,2}\alpha_{5,4}\alpha_{5,2}\alpha_{1,4})
$$
$$
(\alpha_{2,2}\alpha_{4,4}\alpha_{4,2}\alpha_{2,4})
(\alpha_{3,2}\alpha_{3,4}),
$$
$$
(\alpha_{1,2}\alpha_{5,2})(\alpha_{2,2}\alpha_{4,2})
(\alpha_{1,3}\alpha_{5,4})(\alpha_{2,3}\alpha_{4,4})
(\alpha_{3,3}\alpha_{3,4})(\alpha_{4,3}\alpha_{2,4})
$$
$$
(\alpha_{5,3}\alpha_{1,4})(\alpha_{1,5}\alpha_{3,5})]
$$
with orbits
$
\{\alpha_{1,1},\alpha_{1,3},\alpha_{5,1},\alpha_{5,4},
\alpha_{5,3},\alpha_{5,2},\alpha_{1,4},\alpha_{1,2}\},\
\{\alpha_{2,1},\alpha_{2,3},\alpha_{4,1},
\alpha_{4,4},\alpha_{4,3},\alpha_{4,2},
\linebreak
\alpha_{2,4},\alpha_{2,2}\},\
\{\alpha_{3,1},\alpha_{3,3},\alpha_{3,4},\alpha_{3,2}\},\
\{\alpha_{1,5},\alpha_{3,5}\}
$.

\medskip

{\bf n=18,} $H\cong D_{12}$ ($|H|=12$, $i=4$):
$\rk N_H=16$ and $(N_H)^\ast/N_H\cong (\bz/6\bz)^4$.

$$
H_{18,1}=[
(\alpha_{1,2}\alpha_{5,2})(\alpha_{2,2}\alpha_{4,2})
(\alpha_{1,3}\alpha_{5,4})(\alpha_{2,3}\alpha_{4,4})
(\alpha_{3,3}\alpha_{3,4})(\alpha_{4,3}\alpha_{2,4})
$$
$$
(\alpha_{5,3}\alpha_{1,4})(\alpha_{1,5}\alpha_{3,5}),
$$
$$
(\alpha_{1,1}\alpha_{5,1})(\alpha_{2,1}\alpha_{4,1})
(\alpha_{1,2}\alpha_{5,3}\alpha_{1,4}\alpha_{5,2}\alpha_{1,3}\alpha_{5,4})
(\alpha_{2,2}\alpha_{4,3}\alpha_{2,4}\alpha_{4,2}\alpha_{2,3}\alpha_{4,4})
$$
$$
(\alpha_{3,2}\alpha_{3,3}\alpha_{3,4})(\alpha_{1,5}\alpha_{3,5}\alpha_{4,5})]
$$
with orbits
$
\{\alpha_{1,1},\alpha_{5,1}\},\  \{\alpha_{2,1},\alpha_{4,1}\},\
\{\alpha_{1,2},\alpha_{5,2},\alpha_{5,3},\alpha_{1,3},
\alpha_{1,4},\alpha_{5,4}\},\
\{\alpha_{2,2},\alpha_{4,2},\alpha_{4,3},
\linebreak
\alpha_{2,3},\alpha_{2,4},\alpha_{4,4}\},\
\{\alpha_{3,2},\alpha_{3,3},\alpha_{3,4}\},\
\{\alpha_{1,5},\alpha_{3,5},\alpha_{4,5}\}
$.

\medskip

{\bf n=14,} $H\cong C_8$ ($|H|=8$, $i=1$):
$$
H_{14,1}=[
(\alpha_{1,1}\alpha_{1,3}\alpha_{5,2}\alpha_{5,4}
\alpha_{5,1}\alpha_{5,3}\alpha_{1,2}\alpha_{1,4})
(\alpha_{2,1}\alpha_{2,3}\alpha_{4,2}\alpha_{4,4}
\alpha_{4,1}\alpha_{4,3}\alpha_{2,2}\alpha_{2,4})
$$
$$
(\alpha_{3,1}\alpha_{3,3}\alpha_{3,2}\alpha_{3,4})
(\alpha_{1,5}\alpha_{3,5})]
$$
with $Clos(H_{14,1})=H_{26,1}$ above.

\medskip

{\bf n=12,} $H\cong Q_8$ ($|H|=8$, $i=4$):
$\rk N_H=17$, $(N_H)^\ast/N_H \cong
(\bz/8\bz)^2\times (\bz/2\bz)^3$.
$$
H_{12,1}=[
(\alpha_{1,1}\alpha_{5,4}\alpha_{5,1}\alpha_{1,4})
(\alpha_{2,1}\alpha_{4,4}\alpha_{4,1}\alpha_{2,4})
(\alpha_{3,1}\alpha_{3,4})(\alpha_{1,2}\alpha_{5,3}\alpha_{5,2}\alpha_{1,3})
$$
$$
(\alpha_{2,2}\alpha_{4,3}\alpha_{4,2}\alpha_{2,3})(\alpha_{3,2}\alpha_{3,3}),
$$
$$
(\alpha_{1,1}\alpha_{1,2}\alpha_{5,1}\alpha_{5,2})
(\alpha_{2,1}\alpha_{2,2}\alpha_{4,1}\alpha_{4,2})
(\alpha_{3,1}\alpha_{3,2})(\alpha_{1,3}\alpha_{1,4}\alpha_{5,3}\alpha_{5,4})
$$
$$
(\alpha_{2,3}\alpha_{2,4}\alpha_{4,3}\alpha_{4,4})(\alpha_{3,3}\alpha_{3,4})]
$$
with orbits
$
\{\alpha_{1,1},\alpha_{5,4},\alpha_{1,2},\alpha_{5,1},
\alpha_{1,3},\alpha_{5,3},\alpha_{1,4},\alpha_{5,2}\},\
\{\alpha_{2,1},\alpha_{4,4},\alpha_{2,2},
\alpha_{4,1},\alpha_{2,3},\alpha_{4,3},\alpha_{2,4},\alpha_{4,2}\},\
\linebreak
\{\alpha_{3,1},\alpha_{3,4},\alpha_{3,2},\alpha_{3,3}\}
$.

\medskip

{\bf n=10,} $H\cong D_8$ ($|H|=8$, $i=3$):
$\rk N_H=15$ and $(N_H)^\ast/N_H\cong
(\bz/4\bz)^5$.
$$
H_{10,1}=[
(\alpha_{1,1}\alpha_{1,2}\alpha_{5,1}\alpha_{5,2})
(\alpha_{2,1}\alpha_{2,2}\alpha_{4,1}\alpha_{4,2})
(\alpha_{3,1}\alpha_{3,2})
(\alpha_{1,3}\alpha_{1,4}\alpha_{5,3}\alpha_{5,4})
$$
$$
(\alpha_{2,3}\alpha_{2,4}\alpha_{4,3}\alpha_{4,4})
(\alpha_{3,3}\alpha_{3,4}),
$$
$$
(\alpha_{1,2}\alpha_{5,2})(\alpha_{2,2}\alpha_{4,2})
(\alpha_{1,3}\alpha_{5,4})(\alpha_{2,3}\alpha_{4,4})
(\alpha_{3,3}\alpha_{3,4})(\alpha_{4,3}\alpha_{2,4})
$$
$$
(\alpha_{5,3}\alpha_{1,4})(\alpha_{1,5}\alpha_{3,5})]
$$
with orbits
$
\{\alpha_{1,1},\alpha_{1,2},\alpha_{5,1},\alpha_{5,2}\},\
\{\alpha_{2,1},\alpha_{2,2},\alpha_{4,1},\alpha_{4,2}\},\
\{\alpha_{3,1},\alpha_{3,2}\},\
\{\alpha_{1,3},\alpha_{1,4},\alpha_{5,4},\alpha_{5,3}\},\
\linebreak
\{\alpha_{2,3},\alpha_{2,4},\alpha_{4,4},\alpha_{4,3}\},\
\{\alpha_{3,3},\alpha_{3,4}\},\  \{\alpha_{1,5},\alpha_{3,5}\}
$.

\medskip

{\bf n=7,} $H\cong C_6$ ($|H|=6$, $i=2$):
$$
H_{7,1}=
[(\alpha_{1,1}\alpha_{5,1})(\alpha_{2,1}\alpha_{4,1})
(\alpha_{1,2}\alpha_{5,3}\alpha_{1,4}\alpha_{5,2}\alpha_{1,3}\alpha_{5,4})
(\alpha_{2,2}\alpha_{4,3}\alpha_{2,4}\alpha_{4,2}\alpha_{2,3}\alpha_{4,4})
$$
$$
(\alpha_{3,2}\alpha_{3,3}\alpha_{3,4})(\alpha_{1,5}\alpha_{3,5}\alpha_{4,5})]
$$
with $Clos(H_{7,1})=H_{18,1}$ above.

\medskip

{\bf n=6,} $H\cong D_6$ ($|H|=6$, $i=1$):
$\rk N_H=14$ and $(N_H)^\ast/N_H\cong (\bz/6\bz)^2\times (\bz/3\bz)^3$.
$$
H_{6,1}=[
(\alpha_{1,2}\alpha_{5,2})(\alpha_{2,2}\alpha_{4,2})
(\alpha_{1,3}\alpha_{5,4})(\alpha_{2,3}\alpha_{4,4})
(\alpha_{3,3}\alpha_{3,4})(\alpha_{4,3}\alpha_{2,4})
$$
$$
(\alpha_{5,3}\alpha_{1,4})(\alpha_{1,5}\alpha_{3,5}),
$$
$$
(\alpha_{1,2}\alpha_{1,3}\alpha_{1,4})(\alpha_{2,2}\alpha_{2,3}\alpha_{2,4})
(\alpha_{3,2}\alpha_{3,3}\alpha_{3,4})(\alpha_{4,2}\alpha_{4,3}\alpha_{4,4})
$$
$$
(\alpha_{5,2}\alpha_{5,3}\alpha_{5,4})(\alpha_{1,5}\alpha_{3,5}\alpha_{4,5})]
$$
with orbits
$
\{\alpha_{1,2},\alpha_{5,2},\alpha_{1,3},\alpha_{5,3},\alpha_{5,4},\alpha_{1,4}\},\
\{\alpha_{2,2},\alpha_{4,2},\alpha_{2,3},\alpha_{4,3},\alpha_{4,4},
\alpha_{2,4}\},\  \{\alpha_{3,2},\alpha_{3,3},\alpha_{3,4}\},\ \linebreak
\{\alpha_{1,5},\alpha_{3,5},\alpha_{4,5}\}
$;
$$
H_{6,2}=[
(\alpha_{1,1}\alpha_{5,1})(\alpha_{2,1}\alpha_{4,1})(\alpha_{1,3}\alpha_{1,4})
(\alpha_{2,3}\alpha_{2,4})(\alpha_{3,3}\alpha_{3,4})(\alpha_{4,3}\alpha_{4,4})
$$
$$
(\alpha_{5,3}\alpha_{5,4})(\alpha_{1,5}\alpha_{3,5}),
$$
$$
(\alpha_{1,2}\alpha_{1,3}\alpha_{1,4})(\alpha_{2,2}\alpha_{2,3}\alpha_{2,4})
(\alpha_{3,2}\alpha_{3,3}\alpha_{3,4})(\alpha_{4,2}\alpha_{4,3}\alpha_{4,4})
$$
$$
(\alpha_{5,2}\alpha_{5,3}\alpha_{5,4})(\alpha_{1,5}\alpha_{3,5}\alpha_{4,5})]
$$
with orbits
$
\{\alpha_{1,1},\alpha_{5,1}\},\  \{\alpha_{2,1},\alpha_{4,1}\},\
\{\alpha_{1,2},\alpha_{1,3},\alpha_{1,4}\},\
\{\alpha_{2,2},\alpha_{2,3},\alpha_{2,4}\},\
\linebreak
\{\alpha_{3,2},\alpha_{3,3},\alpha_{3,4}\},\
\{\alpha_{4,2},\alpha_{4,3},\alpha_{4,4}\},\
\{\alpha_{5,2},\alpha_{5,3},\alpha_{5,4}\},\
\{\alpha_{1,5},\alpha_{3,5},\alpha_{4,5}\}
$.

\medskip

{\bf n=4,} $H\cong C_4$ ($|H|=4$, $i=1$):
$\rk N_H=14$ and $(N_H)^\ast/N_H \cong
(\bz/4\bz)^4\times (\bz/2\bz)^2$.
$$
H_{4,1}=[
(\alpha_{1,1}\alpha_{1,2}\alpha_{5,1}\alpha_{5,2})
(\alpha_{2,1}\alpha_{2,2}\alpha_{4,1}\alpha_{4,2})
(\alpha_{3,1}\alpha_{3,2})
(\alpha_{1,3}\alpha_{1,4}\alpha_{5,3}\alpha_{5,4})
$$
$$
(\alpha_{2,3}\alpha_{2,4}\alpha_{4,3}\alpha_{4,4})
(\alpha_{3,3}\alpha_{3,4})]
$$
with orbits
$
\{\alpha_{1,1},\alpha_{1,2},\alpha_{5,1},\alpha_{5,2}\},\
\{\alpha_{2,1},\alpha_{2,2},\alpha_{4,1},\alpha_{4,2}\},\
\{\alpha_{3,1},\alpha_{3,2}\},\ \linebreak
\{\alpha_{1,3},\alpha_{1,4},\alpha_{5,3},\alpha_{5,4}\},\
\{\alpha_{2,3},\alpha_{2,4},\alpha_{4,3},\alpha_{4,4}\},\
\{\alpha_{3,3},\alpha_{3,4}\}
$.

\medskip

{\bf n=3,} $H\cong C_2^2$ ($|H|=4$, $i=2$):
$\rk N_H=12$ and $(N_H)^\ast/N_H\cong (\bz/4\bz)^2\times (\bz/2\bz)^6$.
$$
H_{3,1}=[
(\alpha_{1,1}\alpha_{5,1})(\alpha_{2,1}\alpha_{4,1})
(\alpha_{1,3}\alpha_{1,4})
(\alpha_{2,3}\alpha_{2,4})(\alpha_{3,3}\alpha_{3,4})
(\alpha_{4,3}\alpha_{4,4})
$$
$$
(\alpha_{5,3}\alpha_{5,4})(\alpha_{1,5}\alpha_{3,5}),
$$
$$
(\alpha_{1,2}\alpha_{5,2})(\alpha_{2,2}\alpha_{4,2})
(\alpha_{1,3}\alpha_{5,4})(\alpha_{2,3}\alpha_{4,4})
(\alpha_{3,3}\alpha_{3,4})(\alpha_{4,3}\alpha_{2,4})
$$
$$
(\alpha_{5,3}\alpha_{1,4})(\alpha_{1,5}\alpha_{3,5})]
$$
with orbits
$
\{\alpha_{1,1},\alpha_{5,1}\},\  \{\alpha_{2,1},\alpha_{4,1}\},\
\{\alpha_{1,2},\alpha_{5,2}\},\
\{\alpha_{2,2},\alpha_{4,2}\},\
\{\alpha_{1,3},
\alpha_{1,4},\alpha_{5,4},\alpha_{5,3}\},\ \linebreak
\{\alpha_{2,3},\alpha_{2,4},\alpha_{4,4},\alpha_{4,3}\},\
\{\alpha_{3,3},\alpha_{3,4}\},\  \{\alpha_{1,5},\alpha_{3,5}\}
$.

\medskip

{\bf n=2,} $H\cong C_3$ ($|H|=3$, $i=1$):
$\rk N_H=12$ and $(N_H)^\ast/N_H\cong (\bz/3\bz)^6$.
$$
H_{2,1}=[
(\alpha_{1,2}\alpha_{1,3}\alpha_{1,4})
(\alpha_{2,2}\alpha_{2,3}\alpha_{2,4})
(\alpha_{3,2}\alpha_{3,3}\alpha_{3,4})
(\alpha_{4,2}\alpha_{4,3}\alpha_{4,4})
$$
$$
(\alpha_{5,2}\alpha_{5,3}\alpha_{5,4})
(\alpha_{1,5}\alpha_{3,5}\alpha_{4,5})]
$$
with orbits
$
\{\alpha_{1,2},\alpha_{1,3},\alpha_{1,4}\},\
\{\alpha_{2,2},\alpha_{2,3},\alpha_{2,4}\},\
\{\alpha_{3,2},\alpha_{3,3},\alpha_{3,4}\},\
\{\alpha_{4,2},
\alpha_{4,3},\alpha_{4,4}\},\ \linebreak
\{\alpha_{5,2},\alpha_{5,3},\alpha_{5,4}\},\
\{\alpha_{1,5},\alpha_{3,5},\alpha_{4,5}\}
$.

\medskip

{\bf n=1,} $H\cong C_2$ ($|H|=2$, $i=1$):
$\rk N_H=8$ and $(N_H)^\ast/N_H\cong (\bz/2\bz)^8$.
$$
H_{1,1}=[
(\alpha_{1,2}\alpha_{5,2})(\alpha_{2,2}\alpha_{4,2})(\alpha_{1,3}\alpha_{5,4})
(\alpha_{2,3}\alpha_{4,4})(\alpha_{3,3}\alpha_{3,4})(\alpha_{4,3}\alpha_{2,4})
$$
$$
(\alpha_{5,3}\alpha_{1,4})(\alpha_{1,5},\alpha_{3,5})]
$$
with orbits
$
\{\alpha_{1,2},\alpha_{5,2}\},\  \{\alpha_{2,2},\alpha_{4,2}\},\
\{\alpha_{1,3},\alpha_{5,4}\},\  \{\alpha_{2,3},\alpha_{4,4}\},\
\{\alpha_{3,3},
\alpha_{3,4}\},\ \{\alpha_{4,3},\alpha_{2,4}\},\ \linebreak
\{\alpha_{5,3},\alpha_{1,4}\},\  \{\alpha_{1,5},\alpha_{3,5}\}
$;
$$
H_{1,2}=[
(\alpha_{1,1}\alpha_{5,1})(\alpha_{2,1}\alpha_{4,1})
(\alpha_{1,2}\alpha_{5,2})(\alpha_{2,2}\alpha_{4,2})
(\alpha_{1,3}\alpha_{5,3})(\alpha_{2,3}\alpha_{4,3})
$$
$$
(\alpha_{1,4}\alpha_{5,4})(\alpha_{2,4}\alpha_{4,4})]
$$
with orbits
$
\{\alpha_{1,1},\alpha_{5,1}\},\  \{\alpha_{2,1},\alpha_{4,1}\},\
\{\alpha_{1,2},\alpha_{5,2}\},\  \{\alpha_{2,2},\alpha_{4,2}\},\
\{\alpha_{1,3},
\alpha_{5,3}\},\  \{\alpha_{2,3},\alpha_{4,3}\},\ \linebreak
\{\alpha_{1,4},\alpha_{5,4}\},\  \{\alpha_{2,4},\alpha_{4,4}\}
$.

\medskip

In \cite{Nik7} and \cite{Nik8}, these conjugacy classes were described by
direct considerations. Here is the correspondence between these two
descriptions: $H_{54,1}=A(N_{18})$; $H_{38,1}\cong 1$
(subgroups 1 in \cite{Nik7}
and \cite{Nik8}); $H_{26,1}\cong 2$; $H_{18,1}\cong 5$;
$H_{14,1}\cong 3$; $H_{12,1}\cong 4$;
$H_{10,1}\cong 7$; $H_{7,1}\cong 6$; $H_{6,1}\cong
[\varphi_0\widetilde{(23)}\widetilde{(12)}\widetilde{(23)},
\widetilde{(23)}]$ in 8;
$H_{6,2}\cong [\varphi_0\widetilde{(23)},\widetilde{(12)}]$ in 8;
$H_{4,1}\cong 9$;
$H_{3,1}\cong$ 10; $H_{2,1}\cong$ 11; $H_{1,1}\cong [\widetilde{(ij)}]$ and
$[\widetilde{(ij)}\varphi_0]$ in 12; $H_{1,2}\cong [\varphi_0]$ in 12.


\medskip

{\bf Case 17.} For the Niemeier lattice
$$
N=N_{17}=N(4A_6)=
$$
$$
=[4A_6,\
\varepsilon_{1,1}+2\varepsilon_{1,2}+\varepsilon_{1,3}+6\varepsilon_{1,4},\
\varepsilon_{1,1}+6\varepsilon_{1,2}+2\varepsilon_{1,3}+\varepsilon_{1,4},\,
\varepsilon_{1,1}+\varepsilon_{1,2}+6\varepsilon_{1,3}+2\varepsilon_{1,4}],
$$
the group $A(N)$ has order $24$.
It is generated by
$$
\varphi_0=(\alpha_{1,1}\alpha_{6,1})(\alpha_{1,2}\alpha_{6,2})
(\alpha_{1,3}\alpha_{6,3})(\alpha_{1,4}\alpha_{6,4})
(\alpha_{2,1}\alpha_{5,1})(\alpha_{3,1}\alpha_{4,1})
$$
$$
(\alpha_{2,2}\alpha_{5,2})(\alpha_{3,2}\alpha_{4,2})
(\alpha_{2,3}\alpha_{5,3})(\alpha_{3,3}\alpha_{4,3})
(\alpha_{2,4}\alpha_{5,4})(\alpha_{3,4}\alpha_{4,4}),
$$
$$
\varphi_1=(\alpha_{1,2}\alpha_{1,3}\alpha_{1,4})(\alpha_{6,2}\alpha_{6,3}\alpha_{6,4})
(\alpha_{2,2}\alpha_{2,3}\alpha_{2,4})(\alpha_{3,2}\alpha_{3,3}\alpha_{3,4})
$$
$$
(\alpha_{4,2}\alpha_{4,3}\alpha_{4,4})(\alpha_{5,2}\alpha_{5,3}\alpha_{5,4}),
$$
$$
\varphi_2=(\alpha_{1,1}\alpha_{1,2}\alpha_{6,3})
(\alpha_{6,1}\alpha_{6,2}\alpha_{1,3})
(\alpha_{2,1}\alpha_{2,2}\alpha_{5,3})(\alpha_{3,1}\alpha_{3,2}\alpha_{4,3})
$$
$$
(\alpha_{4,1}\alpha_{4,2}\alpha_{3,3})(\alpha_{5,1}\alpha_{5,2}\alpha_{2,3}).
$$

\medskip

\centerline {\bf Classification of KahK3 conjugacy classes for $A(N_{17})$:}

\vskip1cm

{\bf n=2,} $H\cong C_3$ ($|H|=3$, $i=1$):
$\rk N_H=12$ and $(N_H)^\ast/N_H\cong (\bz/3\bz)^6$.
$$
H_{2,1}=[
\varphi_1=(\alpha_{1,2}\alpha_{1,3}\alpha_{1,4})
(\alpha_{2,2}\alpha_{2,3}\alpha_{2,4})
(\alpha_{3,2}\alpha_{3,3}\alpha_{3,4})
(\alpha_{4,2}\alpha_{4,3}\alpha_{4,4})
$$
$$
(\alpha_{5,2}\alpha_{5,3}\alpha_{5,4})
(\alpha_{6,2}\alpha_{6,3}\alpha_{6,4})]
$$
with orbits
$
\{\alpha_{1,2},\alpha_{1,3},\alpha_{1,4}\},\
\{\alpha_{2,2},\alpha_{2,3},\alpha_{2,4}\},\
\{\alpha_{3,2},\alpha_{3,3},\alpha_{3,4}\},\
\{\alpha_{4,2},
\alpha_{4,3},\alpha_{4,4}\},\ \linebreak
\{\alpha_{5,2},\alpha_{5,3},\alpha_{5,4}\},\
\{\alpha_{6,2},\alpha_{6,3},\alpha_{6,4}\}
$.

\medskip

In \cite{Nik7} and \cite{Nik8}, the same result was obtained by
direct considerations: $H_{2,1}=[\varphi_1]\cong C_3$ is the only
non-trivial KahK3 conjugacy class in $A(N_{17})$.

\medskip

{\bf Case 16.} For the Niemeier lattice
$$
N=N_{16}=N(2A_7\oplus 2D_5)=
$$
$$
=[2A_7\oplus 2D_5,\
\varepsilon_{1,1}+\varepsilon_{1,2}+\varepsilon_{1,3}+\varepsilon_{2,4},\
\varepsilon_{1,1}+7\varepsilon_{1,2}+\varepsilon_{2,3}+\varepsilon_{1,4}],
$$
the group $A(N)$ has order $8$, and it is generated by
$$
\varphi_0=
(\alpha_{1,1}\alpha_{7,1})(\alpha_{1,2}\alpha_{7,2})
(\alpha_{4,3}\alpha_{5,3})(\alpha_{4,4}\alpha_{5,4})
$$
$$
(\alpha_{2,1}\alpha_{6,1})(\alpha_{3,1}\alpha_{5,1})
(\alpha_{2,2}\alpha_{6,2})(\alpha_{3,2}\alpha_{5,2}),
$$
$$
\varphi_1=(\alpha_{1,1}\alpha_{1,2})(\alpha_{7,1}\alpha_{7,2})
(\alpha_{4,4}\alpha_{5,4})
$$
$$
(\alpha_{2,1}\alpha_{2,2})(\alpha_{3,1}\alpha_{3,2})
(\alpha_{4,1}\alpha_{4,2})
(\alpha_{5,1}\alpha_{5,2})(\alpha_{6,1}\alpha_{6,2}),
$$
$$
\varphi_2=(\alpha_{1,2}\alpha_{7,2})(\alpha_{4,3}\alpha_{4,4})
(\alpha_{5,3}\alpha_{5,4})
$$
$$
(\alpha_{2,2}\alpha_{6,2})(\alpha_{3,2}\alpha_{5,2})
(\alpha_{1,3}\alpha_{1,4})(\alpha_{2,3}\alpha_{2,4})
(\alpha_{3,3}\alpha_{3,4}).
$$
(see \cite[Ch. 16]{CS} and \cite{Nik7}, \cite{Nik8}).

\medskip

\centerline {\bf Classification of KahK3 conjugacy classes for $A(N_{16})$:}

\medskip

{\bf n=3,} $H\cong C_2^2$ ($|H|=4$, $i=2$):
$\rk N_H=12$ and $(N_H)^\ast/N_H\cong (\bz/4\bz)^2\times (\bz/2\bz)^6$.
$$
H_{3,1}=[
(\alpha_{1,1}\alpha_{7,1})(\alpha_{2,1}\alpha_{6,1})
(\alpha_{3,1}\alpha_{5,1})
(\alpha_{1,3}\alpha_{1,4})
(\alpha_{2,3}\alpha_{2,4})(\alpha_{3,3}\alpha_{3,4})
$$
$$
(\alpha_{4,3}\alpha_{5,4})(\alpha_{5,3}\alpha_{4,4}),
$$
$$
(\alpha_{1,2}\alpha_{7,2})(\alpha_{2,2}\alpha_{6,2})
(\alpha_{3,2}\alpha_{5,2})(\alpha_{1,3}\alpha_{1,4})
(\alpha_{2,3}\alpha_{2,4})(\alpha_{3,3}\alpha_{3,4})
$$
$$
(\alpha_{4,3}\alpha_{4,4})(\alpha_{5,3}\alpha_{5,4})]
$$
with orbits
$
\{\alpha_{1,1},\alpha_{7,1}\},\  \{\alpha_{2,1},\alpha_{6,1}\},\
\{\alpha_{3,1},\alpha_{5,1}\},\
\{\alpha_{1,2},\alpha_{7,2}\},\
\{\alpha_{2,2}, \alpha_{6,2}\},\  \{\alpha_{3,2},\alpha_{5,2}\},\
\linebreak
\{\alpha_{1,3},\alpha_{1,4}\},\
\{\alpha_{2,3},\alpha_{2,4}\},\  \{\alpha_{3,3},\alpha_{3,4}\},\
\{\alpha_{4,3},\alpha_{5,4},
\alpha_{4,4},\alpha_{5,3}\};
$
$$
H_{3,2}=[
(\alpha_{1,1}\alpha_{7,1})(\alpha_{2,1}\alpha_{6,1})
(\alpha_{3,1}\alpha_{5,1})
(\alpha_{1,2}\alpha_{7,2})(\alpha_{2,2}\alpha_{6,2})
(\alpha_{3,2}\alpha_{5,2})
$$
$$
(\alpha_{4,3}\alpha_{5,3})(\alpha_{4,4}\alpha_{5,4}),
$$
$$
(\alpha_{1,1}\alpha_{1,2})(\alpha_{2,1}\alpha_{2,2})
(\alpha_{3,1}\alpha_{3,2})(\alpha_{4,1}\alpha_{4,2})
(\alpha_{5,1}\alpha_{5,2})(\alpha_{6,1}\alpha_{6,2})
$$
$$
(\alpha_{7,1}\alpha_{7,2})(\alpha_{4,4}\alpha_{5,4})]
$$
with orbits
$
\{\alpha_{1,1},\alpha_{7,1},\alpha_{1,2},\alpha_{7,2}\},\
\{\alpha_{2,1},\alpha_{6,1},\alpha_{2,2},\alpha_{6,2}\},\
\{\alpha_{3,1},\alpha_{5,1},
\alpha_{3,2},\alpha_{5,2}\},\ \linebreak
\{\alpha_{4,1},\alpha_{4,2}\},\
\{\alpha_{4,3},\alpha_{5,3}\},\  \{\alpha_{4,4},\alpha_{5,4}\}
$.

\medskip

{\bf n=1,} $H\cong C_2$ ($|H|=2$, $i=1$):
$\rk N_H=8$ and $(N_H)^\ast/N_H\cong (\bz/2\bz)^8$.
$$
H_{1,1}=[
(\alpha_{1,2}\alpha_{7,2})(\alpha_{2,2}\alpha_{6,2})
(\alpha_{3,2}\alpha_{5,2})(\alpha_{1,3}\alpha_{1,4})
(\alpha_{2,3}\alpha_{2,4})(\alpha_{3,3}\alpha_{3,4})
$$
$$
(\alpha_{4,3}\alpha_{4,4})(\alpha_{5,3}\alpha_{5,4})]
$$
with orbits
$
\{\alpha_{1,2},\alpha_{7,2}\},\ \{\alpha_{2,2},\alpha_{6,2}\},\
\{\alpha_{3,2},\alpha_{5,2}\},\ \{\alpha_{1,3},\alpha_{1,4}\},\
\{\alpha_{2,3},\alpha_{2,4}\},\ \{\alpha_{3,3},\alpha_{3,4}\},\
\linebreak
\{\alpha_{4,3},\alpha_{4,4}\},\ \{\alpha_{5,3},\alpha_{5,4}\}
$;

\medskip

$$
H_{1,2}=[
(\alpha_{1,1}\alpha_{7,1})(\alpha_{2,1}\alpha_{6,1})
(\alpha_{3,1}\alpha_{5,1})
(\alpha_{1,2}\alpha_{7,2})(\alpha_{2,2}\alpha_{6,2})
(\alpha_{3,2}\alpha_{5,2})
$$
$$
(\alpha_{4,3}\alpha_{5,3})(\alpha_{4,4}\alpha_{5,4})]
$$
with orbits
$
\{\alpha_{1,1},\alpha_{7,1}\},\ \{\alpha_{2,1},\alpha_{6,1}\},\
\{\alpha_{3,1},\alpha_{5,1}\},\ \{\alpha_{1,2},\alpha_{7,2}\},\
\{\alpha_{2,2},\alpha_{6,2}\},\ \{\alpha_{3,2},\alpha_{5,2}\},\ \linebreak
\{\alpha_{4,3},\alpha_{5,3}\},\ \{\alpha_{4,4},\alpha_{5,4}\}
$;

\medskip

$$
H_{1,3}=[
(\alpha_{1,1}\alpha_{1,2})(\alpha_{2,1}\alpha_{2,2})
(\alpha_{3,1}\alpha_{3,2})(\alpha_{4,1}\alpha_{4,2})
(\alpha_{5,1}\alpha_{5,2})(\alpha_{6,1}\alpha_{6,2})
$$
$$
(\alpha_{7,1}\alpha_{7,2})(\alpha_{4,4}\alpha_{5,4})]
$$
with orbits
$
\{\alpha_{1,1},\alpha_{1,2}\},\ \{\alpha_{2,1},\alpha_{2,2}\},\
\{\alpha_{3,1},\alpha_{3,2}\},\ \{\alpha_{4,1},\alpha_{4,2}\},\
\{\alpha_{5,1},\alpha_{5,2}\},\ \{\alpha_{6,1},\alpha_{6,2}\},\ \linebreak
\{\alpha_{7,1},\alpha_{7,2}\},\ \{\alpha_{4,4},\alpha_{5,4})\}
$.

\medskip

In \cite{Nik7} and \cite{Nik8}, these conjugacy classes were described by
direct considerations. Here is the correspondence between these two
descriptions: $H_{3,1}=[\varphi_0,\varphi_2]$; $H_{3,2}=[\varphi_0,\varphi_1]$;
$H_{1,1}=[\varphi_2]$; $H_{1,2}=[\varphi_0]$; $H_{1,3}=[\varphi_1]$.

\medskip


{\bf Case 15.} For the Niemeier lattice
$$
N=N_{15}=N(3A_8)=
[3A_8,\
4\varepsilon_{1,1}+\varepsilon_{1,2}+\varepsilon_{1,3},\
\varepsilon_{1,1}+4\varepsilon_{1,2}+\varepsilon_{1,3},\
\varepsilon_{1,1}+\varepsilon_{1,2}+4\varepsilon_{1,3}]
$$
the group $A(N)$ has order 12, and it is generated by
$$
\varphi_0=(\alpha_{1,1}\alpha_{8,1})(\alpha_{1,2}\alpha_{8,2})(\alpha_{1,3}\alpha_{8,3})
$$
$$
(\alpha_{2,1}\alpha_{7,1})(\alpha_{3,1}\alpha_{6,1})
(\alpha_{4,1}\alpha_{5,1})(\alpha_{2,2}\alpha_{7,2})
(\alpha_{3,2}\alpha_{6,2})(\alpha_{4,2}\alpha_{5,2})
$$
$$
(\alpha_{2,3}\alpha_{7,3})(\alpha_{3,3}\alpha_{6,3})(\alpha_{4,3}\alpha_{5,3}),
$$
$$
\widetilde{(12)}=(\alpha_{1,1}\alpha_{1,2})(\alpha_{8,1}\alpha_{8,2})
$$
$$
(\alpha_{2,1}\alpha_{2,2})(\alpha_{3,1}\alpha_{3,2})
(\alpha_{4,1}\alpha_{4,2})(\alpha_{5,1}\alpha_{5,2})
$$
$$
(\alpha_{6,1}\alpha_{6,2})(\alpha_{7,1}\alpha_{7,2}),
$$
$$
\widetilde{(23)}=(\alpha_{1,2}\alpha_{1,3})(\alpha_{8,2}\alpha_{8,3})
$$
$$
(\alpha_{2,2}\alpha_{2,3})(\alpha_{3,2}\alpha_{3,3})
(\alpha_{4,2}\alpha_{4,3})(\alpha_{5,2}\alpha_{5,3})
$$
$$
(\alpha_{6,2}\alpha_{6,3})(\alpha_{7,2}\alpha_{7,3}).
$$
(see \cite[Ch. 16]{CS} and \cite{Nik7}, \cite{Nik8}).

\medskip

\centerline {\bf Classification of KahK3 conjugacy classes for $A(N_{15})$:}

\medskip

{\bf n=1,} $H\cong C_2$ ($|H|=2$, $i=1$):
$\rk N_H=8$ and $(N_H)^\ast/N_H\cong (\bz/2\bz)^8$.
$$
H_{1,1}=[
\widetilde{(23)}=(\alpha_{1,2}\alpha_{1,3})(\alpha_{2,2}\alpha_{2,3})
(\alpha_{3,2}\alpha_{3,3})(\alpha_{4,2}\alpha_{4,3})(\alpha_{5,2}\alpha_{5,3})
(\alpha_{6,2}\alpha_{6,3})
$$
$$
(\alpha_{7,2}\alpha_{7,3})(\alpha_{8,2}\alpha_{8,3})]
$$
with orbits
$
\{\alpha_{1,2},\alpha_{1,3} \},\  \{\alpha_{2,2},\alpha_{2,3}\},\
\{\alpha_{3,2},\alpha_{3,3}\},\ \{\alpha_{4,2},\alpha_{4,3}\},\
\{\alpha_{5,2},\alpha_{5,3}\},\
\{\alpha_{6,2},\alpha_{6,3}\},\ \linebreak
\{\alpha_{7,2},\alpha_{7,3}\},\ \{\alpha_{8,2},\alpha_{8,3}\}
$.

In \cite{Nik7} and \cite{Nik8}, the  conjugacy class of
$H_{1,1}=[\widetilde{(23)}]$ was described by direct
considerations.

\medskip


\vskip1cm

{\bf Case 14.} For the Niemeier lattice
$$
N=N_{14}=N(4D_6)=
[4D_6,\ even\ permutations\ of\
0_1+\varepsilon_{1,2}+\varepsilon_{2,3}+\varepsilon_{3,4}]
$$
$$
=[4D_6,\ \varepsilon_{1,2}+\varepsilon_{2,3}+\varepsilon_{3,4},\
\varepsilon_{3,2}+\varepsilon_{1,3}+\varepsilon_{2,4},\
\varepsilon_{2,2}+\varepsilon_{3,3}+\varepsilon_{1,4};\
$$
$$
\varepsilon_{1,1}+\varepsilon_{3,3}+\varepsilon_{2,4},\
\varepsilon_{2,1}+\varepsilon_{1,3}+\varepsilon_{3,4},\
\varepsilon_{3,1}+\varepsilon_{2,3}+\varepsilon_{1,4};\
$$
$$
\varepsilon_{2,1}+\varepsilon_{3,2}+\varepsilon_{1,4},\
\varepsilon_{1,1}+\varepsilon_{2,2}+\varepsilon_{3,4},\
\varepsilon_{3,1}+\varepsilon_{1,2}+\varepsilon_{2,4};\
$$
$$
\varepsilon_{3,1}+\varepsilon_{2,2}+\varepsilon_{1,3},\
\varepsilon_{1,1}+\varepsilon_{3,2}+\varepsilon_{2,3},\
\varepsilon_{2,1}+\varepsilon_{1,2}+\varepsilon_{3,3}]
$$
the group $A(N)$ has order $24$, and it is generated by
$$
(12)=(\alpha_{1,1}\alpha_{1,2})
(\alpha_{5,1}\alpha_{6,2})(\alpha_{6,1}\alpha_{5,2})(\alpha_{5,3}\alpha_{6,3})(\alpha_{5,4}\alpha_{6,4})
$$
$$
(\alpha_{2,1}\alpha_{2,2})(\alpha_{3,1}\alpha_{3,2})(\alpha_{4,1}\alpha_{4,2}),
$$
$$
(23)=(\alpha_{1,2}\alpha_{1,3})
(\alpha_{5,2}\alpha_{6,3})(\alpha_{6,2}\alpha_{5,3})(\alpha_{5,1}\alpha_{6,1})(\alpha_{5,4}\alpha_{6,4})
$$
$$
(\alpha_{2,2}\alpha_{2,3})(\alpha_{3,2}\alpha_{3,3})(\alpha_{4,2}\alpha_{4,3}),
$$
$$
(34)=(\alpha_{1,3}\alpha_{1,4})
(\alpha_{5,3}\alpha_{6,4})(\alpha_{6,3}\alpha_{5,4})(\alpha_{5,1}\alpha_{6,1})(\alpha_{5,2}\alpha_{6,2})
$$
$$
(\alpha_{2,3}\alpha_{2,4})(\alpha_{3,3}\alpha_{3,4})(\alpha_{4,3}\alpha_{4,4})
$$
(see \cite[Ch. 16]{CS} and \cite{Nik7}, \cite{Nik8}).

\medskip

\centerline {\bf Classification of KahK3 conjugacy classes for $A(N_{14})$:}

\medskip

{\bf n=6,} $H\cong D_6$ ($|H|=6$, $i=1$):
$\rk N_H=14$ and $(N_H)^\ast/N_H\cong (\bz/6\bz)^2\times (\bz/3\bz)^3$.
$$
H_{6,1}=[
(\alpha_{5,1}\alpha_{6,1})(\alpha_{5,2}\alpha_{6,2})(\alpha_{1,3}\alpha_{1,4})
(\alpha_{2,3}\alpha_{2,4})(\alpha_{3,3}\alpha_{3,4})(\alpha_{4,3}\alpha_{4,4})
$$
$$
(\alpha_{5,3}\alpha_{6,4})(\alpha_{6,3}\alpha_{5,4}),
$$
$$
(\alpha_{1,2}\alpha_{1,3}\alpha_{1,4})(\alpha_{2,2}\alpha_{2,3}\alpha_{2,4})
(\alpha_{3,2}\alpha_{3,3}\alpha_{3,4})(\alpha_{4,2}\alpha_{4,3}\alpha_{4,4})
$$
$$
(\alpha_{5,2}\alpha_{5,3}\alpha_{5,4})(\alpha_{6,2}\alpha_{6,3}\alpha_{6,4})]
$$
with orbits
$
\{\alpha_{5,1},\alpha_{6,1}\},\  \{\alpha_{1,2},\alpha_{1,3},\alpha_{1,4}\},\
\{\alpha_{2,2},\alpha_{2,3},\alpha_{2,4}\},\
 \{\alpha_{3,2},\alpha_{3,3},\alpha_{3,4}\},\ \linebreak
\{\alpha_{4,2},\alpha_{4,3},\alpha_{4,4}\},\
\{\alpha_{5,2},\alpha_{6,2},\alpha_{5,3},\alpha_{6,3},\alpha_{6,4}, \alpha_{5,4}\}
$.

\medskip

{\bf n=2,} $H\cong C_3$ ($|H|=3$, $i=1$):
$\rk N_H=12$ and $(N_H)^\ast/N_H\cong (\bz/3\bz)^6$.
$$
H_{2,1}=[
(\alpha_{1,2}\alpha_{1,3}\alpha_{1,4})(\alpha_{2,2}\alpha_{2,3}\alpha_{2,4})
(\alpha_{3,2}\alpha_{3,3}\alpha_{3,4})(\alpha_{4,2}\alpha_{4,3}\alpha_{4,4})
$$
$$
(\alpha_{5,2}\alpha_{5,3}\alpha_{5,4})(\alpha_{6,2}\alpha_{6,3}\alpha_{6,4})]
$$
with orbits
$
\{\alpha_{1,2},\alpha_{1,3},\alpha_{1,4}\},\
\{\alpha_{2,2},\alpha_{2,3},\alpha_{2,4}\},\
\{\alpha_{3,2},\alpha_{3,3},\alpha_{3,4}\},\
\{\alpha_{4,2},\alpha_{4,3},\alpha_{4,4}\},\
\linebreak
\{\alpha_{5,2},\alpha_{5,3},\alpha_{5,4}\},\
\{\alpha_{6,2},\alpha_{6,3},\alpha_{6,4}\}
$.

\medskip

{\bf n=1,} $H\cong C_2$ ($|H|=2$, $i=1$):
$\rk N_H=8$ and $(N_H)^\ast/N_H\cong (\bz/2\bz)^8$.
$$
H_{1,1}=[
(\alpha_{5,1}\alpha_{6,1})(\alpha_{5,2}\alpha_{6,2})(\alpha_{1,3}\alpha_{1,4})
(\alpha_{2,3}\alpha_{2,4})(\alpha_{3,3}\alpha_{3,4})(\alpha_{4,3}\alpha_{4,4})
$$
$$
(\alpha_{5,3}\alpha_{6,4})(\alpha_{6,3}\alpha_{5,4})]
$$
with orbits
$
\{\alpha_{5,1},\alpha_{6,1}\},\ \{\alpha_{5,2},\alpha_{6,2}\},\
\{\alpha_{1,3},\alpha_{1,4}\},\ \{\alpha_{2,3},\alpha_{2,4}\},\
\{\alpha_{3,3},\alpha_{3,4}\},\
\{\alpha_{4,3},\alpha_{4,4}\},\ \linebreak
\{\alpha_{5,3},\alpha_{6,4}\},\ \{\alpha_{6,3},\alpha_{5,4}\}
$.

\medskip

In \cite{Nik7} and \cite{Nik8}, these conjugacy classes were described by
direct considerations. Here is the correspondence between these two
descriptions: $H_{6,1}=({\mathfrak D}_6)_1$; $H_{2,1}=(C_3)_1$;
$H_{1,1}=[(34)]$.

\vskip1cm

{\bf Case 13.} For the Niemeier lattice
$$
N=N_{13}=N(2A_9\oplus D_6)
=[2A_9\oplus D_6,\
2\varepsilon_{1,1}+4\varepsilon_{1,2},\,
5\varepsilon_{1,1}+\varepsilon_{1,3},\,
5\varepsilon_{1,2}+\varepsilon_{3,3}],
$$
the group $A(N)$ has order $4$, and
it is generated by
$$
\varphi=(\alpha_{1,1}\alpha_{1,2}\alpha_{9,1}\alpha_{9,2})(\alpha_{5,3}\alpha_{6,3})
$$
$$
(\alpha_{2,1}\alpha_{2,2}\alpha_{8,1}\alpha_{8,2})
(\alpha_{3,1}\alpha_{3,2}\alpha_{7,1}\alpha_{7,2})
(\alpha_{4,1}\alpha_{4,2}\alpha_{6,1}\alpha_{6,2})
(\alpha_{5,1}\alpha_{5,2})
$$
(see \cite[Ch. 16]{CS} and \cite{Nik7}, \cite{Nik8}).

\medskip

\centerline {\bf Classification of KahK3 conjugacy classes for $A(N_{13})$:}

\medskip

{\bf n=4,} $H\cong C_4$ ($|H|=4$, $i=1$):
$\rk N_H=14$ and $(N_H)^\ast/N_H \cong
(\bz/4\bz)^4\times (\bz/2\bz)^2$.
$$
H_{4,1}=[
\varphi=(\alpha_{1,1}\alpha_{1,2}\alpha_{9,1}\alpha_{9,2})(\alpha_{5,3}\alpha_{6,3})
$$
$$
(\alpha_{2,1}\alpha_{2,2}\alpha_{8,1}\alpha_{8,2})
(\alpha_{3,1}\alpha_{3,2}\alpha_{7,1}\alpha_{7,2})
(\alpha_{4,1}\alpha_{4,2}\alpha_{6,1}\alpha_{6,2})
(\alpha_{5,1}\alpha_{5,2})
$$
with orbits
$
\{\alpha_{1,1},\alpha_{1,2},\alpha_{9,1},\alpha_{9,2}\},\
\{\alpha_{2,1},\alpha_{2,2},\alpha_{8,1},\alpha_{8,2}\},\
\{\alpha_{3,1},\alpha_{3,2},\alpha_{7,1},\alpha_{7,2}\},\
\linebreak
\{\alpha_{4,1},\alpha_{4,2},\alpha_{6,1},\alpha_{6,2}\},\
\{\alpha_{5,1},\alpha_{5,2}\},\
\{\alpha_{5,3},\alpha_{6,3}\}
$.

\medskip

{\bf n=1,} $H\cong C_2$ ($|H|=2$, $i=1$):
$\rk N_H=8$ and $(N_H)^\ast/N_H\cong (\bz/2\bz)^8$.
$$
H_{1,1}=[
\varphi^2=(\alpha_{1,1}\alpha_{9,1})(\alpha_{1,2}\alpha_{9,2})
(\alpha_{2,1}\alpha_{8,1})(\alpha_{3,1}\alpha_{7,1})(\alpha_{4,1}\alpha_{6,1})
$$
$$
(\alpha_{2,2}\alpha_{8,2})(\alpha_{3,2}\alpha_{7,2})(\alpha_{4,2}\alpha_{6,2})]
$$
with orbits
$
\{\alpha_{1,1},\alpha_{9,1}\},\ \{\alpha_{1,2},\alpha_{9,2}\},\
\{\alpha_{2,1},\alpha_{8,1},\},\ \{\alpha_{3,1},\alpha_{7,1}\},\
\{\alpha_{4,1},\alpha_{6,1}\},\
\{\alpha_{2,2},\alpha_{8,2}\},\
\linebreak
\{\alpha_{3,2},\alpha_{7,2}\},\
\{\alpha_{4,2},\alpha_{6,2}\}
$.


\vskip1cm

{\bf Case 12.} For the Niemeier lattice
$$
N_{12}=N(4E_6)=
[4E_6,\
\varepsilon_{1,1}+\varepsilon_{1,3}+\varepsilon_{2,4},\
\varepsilon_{1,1}+\varepsilon_{2,2}+\varepsilon_{1,4},\,
$$
$$
\varepsilon_{1,1}+\varepsilon_{1,2}+\varepsilon_{2,3}],
$$
the group $A(N)$ has order $48$, and it is generated by
$$
\varphi_0=(\alpha_{1,1}\alpha_{6,1})(\alpha_{1,2}\alpha_{6,2})
(\alpha_{1,3}\alpha_{6,3})(\alpha_{1,4}\alpha_{6,4})
$$
$$
(\alpha_{3,1}\alpha_{5,1})(\alpha_{3,2}\alpha_{5,2})
(\alpha_{3,3}\alpha_{5,3})(\alpha_{3,4}\alpha_{5,4}),
$$
$$
\widetilde{(12)}=(\alpha_{1,1}\alpha_{1,2})(\alpha_{6,1}\alpha_{6,2})
(\alpha_{1,4}\alpha_{6,4})
$$
$$
(\alpha_{2,1}\alpha_{2,2})(\alpha_{3,1}\alpha_{3,2})(\alpha_{4,1}\alpha_{4,2})
(\alpha_{5,1}\alpha_{5,2})(\alpha_{3,4},\alpha_{5,4}),
$$
$$
\widetilde{(23)}=(\alpha_{1,1}\alpha_{6,1})(\alpha_{1,2}\alpha_{1,3})(\alpha_{6,2}\alpha_{6,3})
$$
$$
(\alpha_{3,1}\alpha_{5,1})
(\alpha_{2,2}\alpha_{2,3})(\alpha_{3,2}\alpha_{3,3})(\alpha_{4,2}\alpha_{4,3})(\alpha_{5,2}\alpha_{5,3}),
$$
$$
\widetilde{(34)}=(\alpha_{1,1}\alpha_{6,1})(\alpha_{1,3}\alpha_{1,4})(\alpha_{6,3}\alpha_{6,4})
$$
$$
(\alpha_{3,1}\alpha_{5,1})
(\alpha_{2,3}\alpha_{2,4})(\alpha_{3,3}\alpha_{3,4})(\alpha_{4,3}\alpha_{4,4})(\alpha_{5,3}\alpha_{5,4}).
$$

\medskip

\centerline {\bf Classification of KahK3 conjugacy classes for $A(N_{12})$:}

\medskip

{\bf n=18,} $H\cong D_{12}$ ($|H|=12$, $i=4$):
$\rk N_H=16$ and $(N_H)^\ast/N_H\cong (\bz/6\bz)^4$.
$$
H_{18,1}=[
(\alpha_{1,2}\alpha_{6,2})(\alpha_{3,2}\alpha_{5,2})(\alpha_{1,3}\alpha_{6,4})
(\alpha_{2,3}\alpha_{2,4})
(\alpha_{3,3}\alpha_{5,4})(\alpha_{4,3}\alpha_{4,4})
$$
$$
(\alpha_{5,3}\alpha_{3,4})(\alpha_{6,3}\alpha_{1,4}),
$$
$$
(\alpha_{1,1}\alpha_{6,1})
(\alpha_{3,1}\alpha_{5,1})(\alpha_{1,2}\alpha_{6,3}\alpha_{1,4}\alpha_{6,2}\alpha_{1,3}\alpha_{6,4})
(\alpha_{2,2}\alpha_{2,3}\alpha_{2,4})
$$
$$
(\alpha_{3,2}\alpha_{5,3}\alpha_{3,4}\alpha_{5,2}\alpha_{3,3}\alpha_{5,4})
(\alpha_{4,2}\alpha_{4,3}\alpha_{4,4})]
$$
with orbits
$
\{\alpha_{1,1},\alpha_{6,1}\},\  \{\alpha_{3,1},\alpha_{5,1}\},\
\{\alpha_{1,2},\alpha_{6,2},\alpha_{6,3},\alpha_{1,3},\alpha_{1,4},\alpha_{6,4}\},\
\{\alpha_{2,2},\alpha_{2,3},\alpha_{2,4}\},\
\linebreak
\{\alpha_{3,2},\alpha_{5,2},\alpha_{5,3},\alpha_{3,3},\alpha_{3,4},\alpha_{5,4}\},\
\{\alpha_{4,2},\alpha_{4,3},\alpha_{4,4}\}
$.

\medskip

{\bf n=7,} $H\cong C_6$ ($|H|=6$, $i=2$):
$$
H_{7,1}=[
(\alpha_{1,1}\alpha_{6,1})
(\alpha_{3,1}\alpha_{5,1})(\alpha_{1,2}\alpha_{6,3}\alpha_{1,4}\alpha_{6,2}\alpha_{1,3}\alpha_{6,4})
(\alpha_{2,2}\alpha_{2,3}\alpha_{2,4})
$$
$$
(\alpha_{3,2}\alpha_{5,3}\alpha_{3,4}\alpha_{5,2}\alpha_{3,3}\alpha_{5,4})
(\alpha_{4,2}\alpha_{4,3}\alpha_{4,4})]
$$
has $Clos(H_{7,1})=H_{18,1}$ above.

\medskip

{\bf n=6,} $H\cong D_6$ ($|H|=6$, $i=1$):
$\rk N_H=14$ and $(N_H)^\ast/N_H\cong (\bz/6\bz)^2\times (\bz/3\bz)^3$.
$$
H_{6,1}=[
(\alpha_{1,2}\alpha_{6,2})(\alpha_{3,2}\alpha_{5,2})
(\alpha_{1,3}\alpha_{6,4})(\alpha_{2,3}\alpha_{2,4})
(\alpha_{3,3}\alpha_{5,4})
(\alpha_{4,3}\alpha_{4,4})
$$
$$
(\alpha_{5,3}\alpha_{3,4})(\alpha_{6,3}\alpha_{1,4}),
$$
$$
(\alpha_{1,2}\alpha_{1,3}\alpha_{1,4})(\alpha_{2,2}\alpha_{2,3}\alpha_{2,4})
(\alpha_{3,2}\alpha_{3,3}\alpha_{3,4})(\alpha_{4,2}\alpha_{4,3}\alpha_{4,4})
$$
$$
(\alpha_{5,2}\alpha_{5,3}\alpha_{5,4})(\alpha_{6,2}\alpha_{6,3}\alpha_{6,4})]
$$
 with orbits
$
\{\alpha_{1,2},\alpha_{6,2},\alpha_{1,3},\alpha_{6,3},\alpha_{6,4},\alpha_{1,4}\},\
\{\alpha_{2,2},\alpha_{2,3},\alpha_{2,4}\},\
\{\alpha_{3,2},\alpha_{5,2},\alpha_{3,3},\alpha_{5,3},
\alpha_{5,4},\alpha_{3,4}\},\ \linebreak
\{\alpha_{4,2},\alpha_{4,3},\alpha_{4,4}\}
$;
$$
H_{6,2}=[
(\alpha_{1,1}\alpha_{6,1})(\alpha_{3,1}\alpha_{5,1})(\alpha_{1,3}\alpha_{1,4})
(\alpha_{2,3}\alpha_{2,4})(\alpha_{3,3}\alpha_{3,4})(\alpha_{4,3}\alpha_{4,4})
$$
$$
(\alpha_{5,3}\alpha_{5,4})(\alpha_{6,3}\alpha_{6,4}),
$$
$$
(\alpha_{1,2}\alpha_{1,3}\alpha_{1,4})(\alpha_{2,2}\alpha_{2,3}\alpha_{2,4})
(\alpha_{3,2}\alpha_{3,3}\alpha_{3,4})(\alpha_{4,2}\alpha_{4,3}\alpha_{4,4})
$$
$$
(\alpha_{5,2}\alpha_{5,3}\alpha_{5,4})
(\alpha_{6,2}\alpha_{6,3}\alpha_{6,4})]
$$
with orbits
$
\{\alpha_{1,1},\alpha_{6,1}\},\  \{\alpha_{3,1},\alpha_{5,1}\},\
\{\alpha_{1,2},\alpha_{1,3},\alpha_{1,4}\},\  \{\alpha_{2,2},\alpha_{2,3},\alpha_{2,4}\},\
\{\alpha_{3,2},\alpha_{3,3},\alpha_{3,4}\},\
\linebreak
\{\alpha_{4,2},\alpha_{4,3},\alpha_{4,4}\},\
\{\alpha_{5,2},\alpha_{5,3},\alpha_{5,4}\},\  \{\alpha_{6,2},\alpha_{6,3},\alpha_{6,4}\}
$.

\medskip

{\bf n=3,} $H\cong C_2^2$ ($|H|=4$, $i=2$):
$\rk N_H=12$ and $(N_H)^\ast/N_H\cong (\bz/4\bz)^2\times (\bz/2\bz)^6$.
$$
H_{3,1}=[
\varphi_0=(\alpha_{1,1}\alpha_{6,1})(\alpha_{1,2}\alpha_{6,2})
(\alpha_{1,3}\alpha_{6,3})(\alpha_{1,4}\alpha_{6,4})
$$
$$
(\alpha_{3,1}\alpha_{5,1})(\alpha_{3,2}\alpha_{5,2})
(\alpha_{3,3}\alpha_{5,3})(\alpha_{3,4}\alpha_{5,4}),
$$
$$
\widetilde{(34)}=(\alpha_{1,1}\alpha_{6,1})(\alpha_{1,3}\alpha_{1,4})(\alpha_{6,3}\alpha_{6,4})
$$
$$
(\alpha_{3,1}\alpha_{5,1})
(\alpha_{2,3}\alpha_{2,4})(\alpha_{3,3}\alpha_{3,4})(\alpha_{4,3}\alpha_{4,4})(\alpha_{5,3}\alpha_{5,4})]
$$
with orbits
$
\{\alpha_{1,1},\alpha_{6,1}\},\  \{\alpha_{3,1},\alpha_{5,1}\},\  \{\alpha_{1,2},\alpha_{6,2}\},\
\{\alpha_{3,2},\alpha_{5,2}\},\
\{\alpha_{1,3},\alpha_{6,3},\alpha_{6,4},\alpha_{1,4}\},\
\linebreak
\{\alpha_{2,3},\alpha_{2,4}\},\
\{\alpha_{3,3},\alpha_{5,3},\alpha_{5,4},\alpha_{3,4}\},\  \{\alpha_{4,3},\alpha_{4,4}\}
$.

\medskip

{\bf n=2,} $H\cong C_3$ ($|H|=3$, $i=1$):
$\rk N_H=12$ and $(N_H)^\ast/N_H\cong (\bz/3\bz)^6$.
$$
H_{2,1}=[
(\alpha_{1,2}\alpha_{1,3}\alpha_{1,4})(\alpha_{2,2}\alpha_{2,3}\alpha_{2,4})
(\alpha_{3,2}\alpha_{3,3}\alpha_{3,4})(\alpha_{4,2}\alpha_{4,3}\alpha_{4,4})
$$
$$
(\alpha_{5,2}\alpha_{5,3}\alpha_{5,4})(\alpha_{6,2}\alpha_{6,3}\alpha_{6,4})]
$$
with orbits
$
\{\alpha_{1,2},\alpha_{1,3},\alpha_{1,4}\},\
\{\alpha_{2,2},\alpha_{2,3},\alpha_{2,4}\},\
\{\alpha_{3,2},\alpha_{3,3},\alpha_{3,4}\},\
\{\alpha_{4,2},\alpha_{4,3},\alpha_{4,4}\},\
\linebreak
\{\alpha_{5,2},\alpha_{5,3},\alpha_{5,4}\},\
\{\alpha_{6,2},\alpha_{6,3},\alpha_{6,4}\}
$.

\medskip

{\bf n=1,} $H\cong C_2$ ($|H|=2$, $i=1$):
$\rk N_H=8$ and $(N_H)^\ast/N_H\cong (\bz/2\bz)^8$.
$$
H_{1,1}=[
\varphi_0=(\alpha_{1,1}\alpha_{6,1})(\alpha_{1,2}\alpha_{6,2})
(\alpha_{1,3}\alpha_{6,3})(\alpha_{1,4}\alpha_{6,4})
$$
$$
(\alpha_{3,1}\alpha_{5,1})(\alpha_{3,2}\alpha_{5,2})
(\alpha_{3,3}\alpha_{5,3})(\alpha_{3,4}\alpha_{5,4})],
$$
with orbits
$
\{\alpha_{1,1},\alpha_{6,1}\},\ \{\alpha_{1,2},\alpha_{6,2}\},\
\{\alpha_{1,3},\alpha_{6,3}\},\ \{\alpha_{1,4},\alpha_{6,4}\},\
\{\alpha_{3,1},\alpha_{5,1}\},\ \{\alpha_{3,2},\alpha_{5,2}\},\
\linebreak
\{\alpha_{3,3},\alpha_{5,3}\},\ \{\alpha_{3,4},\alpha_{5,4}\}
$;

$$
H_{1,2}=[
\widetilde{(34)}=(\alpha_{1,1}\alpha_{6,1})
(\alpha_{1,3}\alpha_{1,4})(\alpha_{6,3}\alpha_{6,4})
$$
$$
(\alpha_{3,1}\alpha_{5,1})
(\alpha_{2,3}\alpha_{2,4})(\alpha_{3,3}\alpha_{3,4})
(\alpha_{4,3}\alpha_{4,4})(\alpha_{5,3}\alpha_{5,4})]
$$
with orbits
$
\{\alpha_{1,1},\alpha_{6,1}\},\ \{\alpha_{1,3},\alpha_{1,4}\},\
\{\alpha_{6,3},\alpha_{6,4}\},\ \{\alpha_{3,1},\alpha_{5,1}\},\
\{\alpha_{2,3},\alpha_{2,4}\},\ \{\alpha_{3,3},\alpha_{3,4}\},\
\linebreak
\{\alpha_{4,3},\alpha_{4,4}\},\ \{\alpha_{5,3},\alpha_{5,4}\}
$.

In \cite{Nik7} and \cite{Nik8}, these conjugacy classes were described by
direct considerations. Here is the correspondence between these two
descriptions: $H_{12,1}=({\mathfrak D}_{12})_1$;
$H_{7,1}=(C_6)_1$; \newline
$\{H_{6,1},H_{6,2}\}=\{({\mathfrak D}_{6})_{11},({\mathfrak D}_{6})_{12}\}$;
$H_{3,1}=[\varphi_0,\widetilde{(34)}]$;
$H_{2,1}=(C_3)_1$; $H_{1,1}=[\varphi_0]$; $H_{1,2}=[\widetilde{(34)}]$.

\vskip1cm

{\bf Case 11.} For the Niemeier lattice
$N_{11}=N(A_{11}\oplus D_7\oplus E_6)=[A_{11}\oplus D_{7}\oplus E_6,\
\varepsilon_{1,1}+\varepsilon_{1,2}+\varepsilon_{1,3}]$,
the group $A(N_{11})$ has order $2$, and it is generated by
$$
\varphi=(\alpha_{1,1}\alpha_{11,1})(\alpha_{2,1}\alpha_{10,1})
(\alpha_{3,1}\alpha_{9,1})(\alpha_{4,1}\alpha_{8,1})
(\alpha_{5,1}\alpha_{7,1})(\alpha_{6,2}\alpha_{7,2})
$$
$$
(\alpha_{1,3}\alpha_{6,3})(\alpha_{3,3}\alpha_{5,3}).
$$

\medskip

\centerline {\bf Classification of KahK3 conjugacy classes for $A(N_{11})$:}

\medskip

{\bf n=1,} $H\cong C_2$ ($|H|=2$, $i=1$):
$\rk N_H=8$ and $(N_H)^\ast/N_H\cong (\bz/2\bz)^8$.
$$
H_{1,1}=[\varphi]
$$
with orbits
$
\{\alpha_{1,1},\alpha_{11,1}\},\ \{\alpha_{2,1},\alpha_{10,1}\},\
\{\alpha_{3,1},\alpha_{9,1}\},\ \{\alpha_{4,1},\alpha_{8,1}\},\
\{\alpha_{5,1},\alpha_{7,1}\},\
\{\alpha_{6,2},\alpha_{7,2}\},\ \linebreak
\{\alpha_{1,3},\alpha_{6,3}\},\ \{\alpha_{3,3},\alpha_{5,3}\}
$.
See \cite[Ch. 16]{CS} and \cite{Nik7}, \cite{Nik8} for details.

\vskip1cm

{\bf Case 9.} For the Niemeier lattice
$$
N_9=N(3D_8)=
[3D_8,\
\varepsilon_{1,1}+\varepsilon_{2,2}+\varepsilon_{2,3},\
\varepsilon_{2,1}+\varepsilon_{1,2}+\varepsilon_{2,3},\
\varepsilon_{2,1}+\varepsilon_{2,2}+\varepsilon_{1,3}],
$$
the group $A(N_9)$ has order $6$, and it is generated by
$$
(12)=(\alpha_{1,1}\alpha_{1,2})(\alpha_{2,1}\alpha_{2,2})
(\alpha_{3,1}\alpha_{3,2})(\alpha_{4,1}\alpha_{4,2})
(\alpha_{5,1}\alpha_{5,2})(\alpha_{6,1}\alpha_{6,2})
$$
$$
(\alpha_{7,1}\alpha_{7,2})(\alpha_{8,1}\alpha_{8,2}),
$$
$$
(23)=(\alpha_{1,2}\alpha_{1,3})(\alpha_{2,2}\alpha_{2,3})
(\alpha_{3,2}\alpha_{3,3})(\alpha_{4,2}\alpha_{4,3})
(\alpha_{5,2}\alpha_{5,3})(\alpha_{6,2}\alpha_{6,3})
$$
$$
(\alpha_{7,2}\alpha_{7,3})(\alpha_{8,2}\alpha_{8,3}).
$$

\medskip

\centerline {\bf Classification of KahK3 conjugacy classes for $A(N_{9})$:}

\medskip

{\bf n=1,} $H\cong C_2$ ($|H|=2$, $i=1$):
$\rk N_H=8$ and $(N_H)^\ast/N_H\cong (\bz/2\bz)^8$.
$$
H_{1,1}=[(12)]
$$
with orbits
$
\{\alpha_{1,1},\alpha_{1,2}\},\ \{\alpha_{2,1},\alpha_{2,2}\},\
\{\alpha_{3,1},\alpha_{3,2}\},\ \{\alpha_{4,1},\alpha_{4,2}\},\
\{\alpha_{5,1},\alpha_{5,2}\},\ \{\alpha_{6,1},\alpha_{6,2}\},\
\linebreak
\{\alpha_{7,1}\alpha_{7,2}\},\ \{\alpha_{8,1}\alpha_{8,2}\}
$.
See \cite[Ch. 16]{CS} and \cite{Nik7}, \cite{Nik8} for details.

\vskip1cm

{\bf Case 8.} For the Niemeier lattice
$N_8=N(A_{15}\oplus D_9)=[A_{15}\oplus D_9,\
2\varepsilon_{1,1}+\varepsilon_{1,2}]$,
the group $A(N_8)$ has order $2$, and it is
generated by
$$
\varphi=(\alpha_{1,1}\alpha_{15,1})(\alpha_{2,1}\alpha_{14,1})
(\alpha_{3,1}\alpha_{13,1})(\alpha_{4,1}\alpha_{12,1})
(\alpha_{5,1}\alpha_{11,1})(\alpha_{6,1}\alpha_{10,1})
$$
$$
(\alpha_{7,1}\alpha_{9,1})(\alpha_{8,2}\alpha_{9,2}).
$$

\medskip

\centerline {\bf Classification of KahK3 conjugacy classes for $A(N_{8})$:}

\medskip

{\bf n=1,} $H\cong C_2$ ($|H|=2$, $i=1$):
$\rk N_H=8$ and $(N_H)^\ast/N_H\cong (\bz/2\bz)^8$.
$$
H_{1,1}=[\varphi]
$$
with orbits
$
\{\alpha_{1,1},\alpha_{15,1}\},\ \{\alpha_{2,1},\alpha_{14,1}\},\
\{\alpha_{3,1},\alpha_{13,1}\},\ \{\alpha_{4,1},\alpha_{12,1}\},\
\{\alpha_{5,1},\alpha_{11,1}\},\ \{\alpha_{6,1},\alpha_{10,1}\},\
\linebreak
\{\alpha_{7,1},\alpha_{9,1}\},\ \{\alpha_{8,2},\alpha_{9,2}\}
$.
See \cite[Ch. 16]{CS} and \cite{Nik7}, \cite{Nik8} for details.

\vskip1cm

{\bf Case 7.} For the Niemeier lattice
$N_7=N(D_{10}\oplus 2E_7)=[D_{10}\oplus 2E_7,\
\varepsilon_{1,1}+\varepsilon_{1,2},\ \varepsilon_{3,1}+\varepsilon_{1,3}]$,
the group $A(N_7)$ has order $2$, and it is generated by
$$
\varphi=(\alpha_{9,1}\alpha_{10,1})(\alpha_{1,2}\alpha_{1,3})
(\alpha_{2,2}\alpha_{2,3})(\alpha_{3,2}\alpha_{3,3})
(\alpha_{4,2}\alpha_{4,3})(\alpha_{5,2}\alpha_{5,3})
$$
$$
(\alpha_{6,2}\alpha_{6,3})(\alpha_{7,2}\alpha_{7,3}).
$$

\medskip

\centerline {\bf Classification of KahK3 conjugacy classes for $A(N_{7})$:}

\medskip

{\bf n=1,} $H\cong C_2$ ($|H|=2$, $i=1$):
$\rk N_H=8$ and $(N_H)^\ast/N_H\cong (\bz/2\bz)^8$.
$$
H_{1,1}=[\varphi]
$$
with orbits
$
\{\alpha_{9,1},\alpha_{10,1}\},\ \{\alpha_{1,2},\alpha_{1,3}\},\
\{\alpha_{2,2},\alpha_{2,3}\},\ \{\alpha_{3,2},\alpha_{3,3}\},\
\{\alpha_{4,2},\alpha_{4,3}\},\ \{\alpha_{5,2},\alpha_{5,3}\},\
\linebreak
\{\alpha_{6,2},\alpha_{6,3}\},\ \{\alpha_{7,2},\alpha_{7,3}\}
$.
See \cite[Ch. 16]{CS} and \cite{Nik7}, \cite{Nik8} for details.

\vskip1cm

{\bf Case 6.} For the Niemeier lattice
$N_6=N(A_{17}\oplus E_7)=[A_{17}\oplus E_7,
3\varepsilon_{1,1}+\varepsilon_{1,2}]$,
the group $A(N_6)$ has order $2$, and it is generated by
$$
\varphi=(\alpha_{1,1}\alpha_{17,1})(\alpha_{2,1}\alpha_{16,1})
(\alpha_{3,1}\alpha_{15,1})(\alpha_{4,1}\alpha_{14,1})
(\alpha_{5,1}\alpha_{13,1})(\alpha_{6,1}\alpha_{12,1})
$$
$$
(\alpha_{7,1}\alpha_{11,1})(\alpha_{8,1}\alpha_{10,1}).
$$

\medskip

\centerline {\bf Classification of KahK3 conjugacy classes for $A(N_{6})$:}

\medskip

{\bf n=1,} $H\cong C_2$ ($|H|=2$, $i=1$):
$\rk N_H=8$ and $(N_H)^\ast/N_H\cong (\bz/2\bz)^8$.
$$
H_{1,1}=[\varphi]
$$
with orbits
$
\{\alpha_{1,1},\alpha_{17,1}\},\ \{\alpha_{2,1},\alpha_{16,1}\},\
\{\alpha_{3,1},\alpha_{15,1}\},\ \{\alpha_{4,1},\alpha_{14,1}\},\
\{\alpha_{5,1},\alpha_{13,1}\},\ \{\alpha_{6,1},\alpha_{12,1}\},\
\linebreak
\{\alpha_{7,1},\alpha_{11,1}\},\ \{\alpha_{8,1},\alpha_{10,1}\}
$.
See \cite[Ch. 16]{CS} and \cite{Nik7}, \cite{Nik8} for details.

\vskip1cm

{\bf Case 3.} For the Niemeier lattice $N_3=N(3E_8)=3E_8$,
the group $A(N_3)$ has order $6$, and it is generated
by
$$
(12)=(\alpha_{1,1}\alpha_{1,2})(\alpha_{2,1}\alpha_{2,2})
(\alpha_{3,1}\alpha_{3,2})(\alpha_{4,1}\alpha_{4,2})
(\alpha_{5,1}\alpha_{5,2})(\alpha_{6,1}\alpha_{6,2})
$$
$$
(\alpha_{7,1}\alpha_{7,2})(\alpha_{8,1}\alpha_{8,2}),
$$
$$
(23)=(\alpha_{1,2}\alpha_{1,3})(\alpha_{2,2}\alpha_{2,3})
(\alpha_{3,2}\alpha_{3,3})(\alpha_{4,2}\alpha_{4,3})
(\alpha_{5,2}\alpha_{5,3})(\alpha_{6,2}\alpha_{6,3})
$$
$$
(\alpha_{7,2}\alpha_{7,3})(\alpha_{8,2}\alpha_{8,3}).
$$

\medskip

\centerline {\bf Classification of KahK3 conjugacy classes for $A(N_{3})$:}

\medskip

{\bf n=1,} $H\cong C_2$ ($|H|=2$, $i=1$):
$\rk N_H=8$ and $(N_H)^\ast/N_H\cong (\bz/2\bz)^8$.
$$
H_{1,1}=[(12)]
$$
with orbits
$
\{\alpha_{1,1},\alpha_{1,2}\},\ \{\alpha_{2,1},\alpha_{2,2}\},\
\{\alpha_{3,1},\alpha_{3,2}\},\ \{\alpha_{4,1},\alpha_{4,2}\},\
\{\alpha_{5,1},\alpha_{5,2}\},\ \{\alpha_{6,1},\alpha_{6,2}\},\
\linebreak
\{\alpha_{7,1}\alpha_{7,2}\},\ \{\alpha_{8,1}\alpha_{8,2}\}
$.
See \cite[Ch. 16]{CS} and \cite{Nik7}, \cite{Nik8} for details.

\vskip1cm

{\bf Cases 1, 2, 4, 5, 10.} For remaining Niemeier
lattices $N_1$, $N_2$, $N_4$, $N_5$, $N_{10}$,
KahK3 conjugacy classes are only trivial.
See \cite[Ch. 16]{CS} and \cite{Nik7}, \cite{Nik8} for details.


\section{Markings by Niemeier lattices of degenerations of codimension $1$ of
non-trivial finite symplectic automorphism groups of K\"ahlerian K3 surfaces}
\label{sec:Markdegen}

Like for the example of Kummer surfaces in Sect. \ref{sec:Kummer}, using Theorem \ref{th:primembb1}
(we use Program 6 in Section \ref{sec:Appendix} based on this Theorem),
we classify all possible markings
by Niemeier lattices $N_j$, $j=1,\dots, 23$, of degenerations of codimension $1$
of K\"ahlerian K3 surfaces. For Niemeier lattices $N_{23}$ and $N_{22}$ we use
classification of KahK3 conjugacy classes from \cite{Nik7} and \cite{Nik8}.
For Niemeier lattices $N_{j}$, $j=21,\dots, 1$, we use similar style
classification from Section \ref{sec:KahK3cl}. It will be more convenient than
our classification of these cases in \cite{Nik7} and \cite{Nik8}.
Here we label by $n=81,\dots, 1$ types of abstract
finite groups of symplectic
automorphisms of  K\"ahlerian K3 surfaces.
See \cite{Nik0}, \cite{Muk}, \cite{Xiao} and \cite{Hash}.

\begin{theorem} Let $N_j$, $j=1,2,\dots,23$, be
one of Niemeier lattices.

Types of degenerations of codimension one
$N_H\subset S\subset N_{j}$
of maximal ($Clos(H)=H$)  non-trivial finite symplectic
automorphism groups $H$
of K\"ahlerian K3 surfaces which can
be marked by the Niemeier lattice $N_{j}$ are shown
in the Table j below. They are given by orbits of $\alpha\in P(N_j)$
of K\"ahlerian K3 conjugation classes $H$ of
$A(N_{j})$ classified in \cite{Nik7} for $j=23,\,22$
(we use notations of Cases 23 and 22 from \cite{Nik7} respectively),
and in Section \ref{sec:KahK3cl}
for $j=21,20,\dots, 11,\,9,\,8,\,7,\,6,\,3$ (we use notations
of the Case j from Section \ref{sec:KahK3cl}).
We give the type of the root system and the
isomorphism class of the discriminant group $A_S$  of
the corresponding lattice
$S=[N_H, \alpha]_{pr}\subset N_{j}$.
\label{theorem:degj}
\end{theorem}

\vskip1cm

{\bf Table 23. All possible markings by $N_{23}=N(24A_1)$ of degenerations of
codimension one of non-trivial finite symplectic automorphism groups
of K\"ahlerian K3 surfaces:}

\medskip

{\bf n=81,} $H\cong M_{20}$ ($|H|=960$, $i=11357$):
$\rk N_H=19$. {\it No degenerations.}

\medskip

{\bf n=80,} $H\cong F_{384}$ ($|H|=384$, $i=18135$):
$\rk N_H=19$. {\it No degenerations.}

\medskip

{\bf n=79,} $H\cong {\mathfrak A}_{6}$ ($|H|=360$, $i=118$):
$\rk N_H=19$. {\it No degenerations.}

\medskip

{\bf n=78,} $H\cong {\mathfrak A}_{4,4}$ ($|H|=288$, $i=1026$):
$\rk N_H=19$. {\it No degenerations.}

\medskip

{\bf n=77,} $H\cong {T}_{192}$ ($|H|=192$, $i=1493$):
$\rk N_H=19$. {\it No degenerations.}

\medskip

{\bf n=76,} $H\cong {H}_{192}$ ($|H|=192$, $i=955$):
$\rk N_H=19$. {\it No degenerations.}

\medskip

{\bf n=75,} $H\cong 4^2{\mathfrak A}_4$ ($|H|=192$, $i=1023$):
$\rk N_H=18$, $(N_H)^\ast/N_H\cong (\bz/8\bz)^2\times
(\bz/2\bz)^2$. {\it Degenerations:}  Type $16\aaa_1$
with $A_S\cong (\bz/8\bz)\times (\bz/2\bz)^2$
($H_{75,1}$ with the orbit of $\alpha_{1}$).

\medskip

{\bf n=74,} $H\cong L_2(7)$ ($|H|=168$, $i=42$):
$\rk N_H=19$. {\it No degenerations.}

\medskip

{\bf n=70,} $H\cong {\mathfrak S}_5$ ($|H|=120$, $i=34$):
$\rk N_H=19$. {\it No degenerations.}

\medskip

{\bf n=65,} $H\cong2^4D_6$ ($|H|=96$, $i=227$):
$\rk N_H=18$,
$(N_H)^\ast/N_H\cong \bz/24\bz \times
\bz/4\bz\times (\bz/2\bz)^2$. {\it Degenerations:}
$4\aaa_1$ with $A_S\cong \bz/12\bz\times (\bz/4\bz)^2$
($H_{65,3}$ with the orbit of $\alpha_3$ or $\alpha_{10}$;
$H_{65,4}$ with the orbit of $\alpha_3$, $\alpha_4$ or $\alpha_7$);
$8\aaa_1$ with $A_S\cong \bz/24\bz\times (\bz/2\bz)^2$
($H_{65,2}$ with the obit of $\alpha_1$; $H_{65,4}$
with the orbit of $\alpha_1$); $12\aaa_{1}$ with
$A_S\cong (\bz/4\bz)^3$ ($H_{65,3}$ with the
orbit of $\alpha_1$);
$16\aaa_1$ with $A_S\cong \bz/12\bz\times (\bz/2\bz)^2$
($H_{65,1}$ with the orbit of $\alpha_1$).

\medskip

{\bf n=63,} $H\cong M_9$ ($|H|=72$, $i=41$):
$\rk N_H=19$. {\it No degenerations.}

\medskip

{\bf n=62,} $H\cong N_{72}$ ($|H|=72$, $i=40$):
$\rk N_H=19$. {\it No degenerations.}

\medskip

{\bf n=61,} $H\cong {\mathfrak A}_{4,3}$ ($|H|=72$, $i=43$):
$\rk N_H=18$ and
$(N_H)^\ast/N_H\cong (\bz/12\bz)^2 \times \bz/3\bz$.
{\it Degenerations:} $3\aaa_1$ with
$A_S\cong (\bz/12\bz)^2 \times \bz/2\bz$ ($H_{61,1}$
with the orbit of $\alpha_7$ or $\alpha_{16}$);
$12\aaa_1$ with $A_S\cong \bz/24\bz\times \bz/3\bz$
($H_{61,1}$ with the orbit of $\alpha_2$).

\medskip

{\bf n=56,} $H\cong \Gamma_{25}a_1$ ($|H|=64$, $i=138$):
$\rk N_H=18$, $(N_H)^\ast/N_H\cong
\bz/8\bz\times (\bz/4\bz)^3$.
{\it Degenerations:} $8\aaa_1$ with
$A_S\cong \bz/8\bz\times (\bz/4\bz)^2$
($H_{56,1}$ with the orbit of $\alpha_1$ or
$\alpha_6$); $16\aaa_1$ with
$A_S\cong (\bz/4\bz)^3$
($H_{56,2}$ with the orbit of $\alpha_2$).

\medskip

{\bf n=55,} $H\cong {\mathfrak A}_5$ ($|H|=60$, $i=5$):
$\rk N_H=18$ and $(N_H)^\ast/N_H\cong \bz/30\bz\times \bz/10\bz$.
{\it Degenerations:}
$\aaa_1$ with $A_S\cong \bz/30\bz\times \bz/10\bz\times \bz/2\bz$
($H_{55,1}$ with the orbit of $\alpha_1$, $\alpha_{10}$ or
$\alpha_{24}$; $H_{55,2}$ with the orbit of $\alpha_1$, $\alpha_5$,
$\alpha_{21}$ or $\alpha_{24}$);
$5\aaa_1$ with
$A_S\cong \bz/30\bz\times (\bz/2\bz)^2$
($H_{55,1}$ with the orbit of $\alpha_2$;
$H_{55,2}$ with the orbit of $\alpha_4$);
$6\aaa_1$ with $A_S\cong \bz/20\bz\times \bz/5\bz$
($H_{55,1}$ with the orbit of $\alpha_4$);
$10\aaa_1$ with $A_S\cong \bz/60\bz$
($H_{55,1}$ with the orbit of $\alpha_3$);
$15\aaa_1$ with $A_S\cong \bz/10\bz \times (\bz/2\bz)^2$
($H_{55,2}$ with the orbit of $\alpha_2$).

\medskip

{\bf n=54,} $H\cong T_{48}$ ($|H|=48$, $i=29$):
$\rk N_H=19$.
{\it No degenerations.}

\medskip

{\bf n=51,} $H\cong C_2\times {\mathfrak S}_4$ ($|H|=48$, $i=48$):
$\rk N_H=18$, $(N_H)^\ast/N_H\cong (\bz/12\bz)^2\times (\bz/2\bz)^2$.
{\it Degenerations:}
$2\aaa_1$ with $A_S\cong (\bz/12\bz)^2\times  \bz/4\bz$
($H_{51,2}$ with the orbit of $\alpha_2$ or $\alpha_{10}$;
$H_{51,3}$ with the orbit of $\alpha_2$; $H_{51,4}$ with the orbit
of $\alpha_3$ or $\alpha_6$);
$4\aaa_1$ with $A_S\cong \bz/24\bz \times \bz/6\bz \times \bz/2\bz$
($H_{51,2}$ with the orbit of $\alpha_1$; $H_{51,4}$ with the orbit of
$\alpha_1$ or $\alpha_5$);
$6\aaa_1$ with $A_S\cong \bz/12\bz \times (\bz/4\bz)^2$
($H_{51,3}$ with the orbit of $\alpha_4$; $H_{51,4}$ with the orbit of
$\alpha_2$ or $\alpha_7$);
$8\aaa_1$ with $A_S\cong \bz/12\bz\times  \bz/6\bz\times \bz/2\bz$
($H_{51,1}$ with the orbit of $\alpha_1$ or $\alpha_4$;
$H_{51,3}$ with the orbit of $\alpha_1$);
$12\aaa_1$ with $A_S\cong \bz/24\bz\times (\bz/2\bz)^2$
($H_{51,2}$ with the orbit of $\alpha_4$).

\medskip

{\bf n=49,} $H\cong 2^4 C_3$ ($|H|=48$, $i=50$):
$\rk N_H=17$, $(N_H)^\ast/N_H \cong
\bz/24\bz\times (\bz/2\bz)^4$. {\it Degenerations of codimension $1$:}
$4\aaa_1$ with $A_S\cong \bz/12\bz \times \bz/4\bz\times (\bz/2\bz)^2$
($H_{49,1}$ with the orbit of $\alpha_2$ or $\alpha_4$;
$H_{49,3}$ with the orbit of $\alpha_2$, $\alpha_4$, $\alpha_8$, $\alpha_{10}$
or $\alpha_{11}$);
$12\aaa_1$ with $A_S\cong (\bz/4\bz)^2 \times (\bz/2\bz)^2$
($H_{49,1}$ with the orbit of $\alpha_1$);
$16\aaa_1$ with $A_S\cong \bz/6\bz \times (\bz/2\bz)^3$
($H_{49,2}$ with the orbit of $\alpha_2$).

\medskip

{\bf n=48,} $H\cong {\mathfrak S}_{3,3}$ ($|H|=36$, $i=10$):
$\rk N_H=18$, $(N_H)^\ast/N_H \cong
\bz/18\bz\times
\bz/6\bz\times (\bz/3\bz)^2$.
{\it Degenerations:}
$3\aaa_1$ with $A_S\cong \bz/18\bz \times (\bz/6\bz)^2$
($H_{48,1}$ with the orbit of $\alpha_5$ or $\alpha_7$);
$6\aaa_1$ with $A_S\cong \bz/36\bz \times (\bz/3\bz)^2$
($H_{48,1}$ with the orbit of $\alpha_2$);
$9\aaa_1$ with $A_S\cong (\bz/6\bz)^3$
($H_{48,1}$ with the orbit of $\alpha_4$).

\medskip

{\bf n=46,} $H\cong 3^2 C_4$ ($|H|=36$, $i=9$):
$\rk N_H=18$ and $(N_H)^\ast/N_H\cong
\bz/18\bz\times \bz/6\bz\times \bz/3\bz$.
{\it Degenerations:}
$6\aaa_1$ with $A_S\cong \bz/36\bz \times \bz/3\bz$
($H_{46,1}$ with the orbit of $\alpha_1$ or $\alpha_3$);
$9\aaa_1$ with $A_S\cong (\bz/6\bz)^2 \times \bz/2\bz$
($H_{46,1}$ with the orbit of $\alpha_2$);

\medskip

{\bf n=40,} $H\cong Q_8*Q_8$ ($|H|=32$, $i=49$):
$\rk N_H=17$, $(N_H)^\ast/N_H\cong
(\bz/4\bz)^5$.
{\it Degenerations of codimension $1$:}
$8\aaa_1$ with $A_S\cong (\bz/4\bz)^4$
($H_{40,1}$ with the orbit of $\alpha_1$ or $\alpha_3$).

\medskip

{\bf n=39,} $H\cong 2^4C_2$ ($|H|=32$, $i=27$):
$\rk N_H=17$, $(N_H)^\ast/N_H\cong
\bz/8\bz\times (\bz/4\bz)^2\times (\bz/2\bz)^2$.
{\it Degenerations of codimension $1$:}
$4\aaa_1$ with $A_S\cong (\bz/4\bz)^4$
($H_{39,2}$ with the orbit of $\alpha_2$, $\alpha_3$
or $\alpha_4$);
$8\aaa_1$ with $A_S\cong \bz/8\bz\times \bz/4\bz\times (\bz/2\bz)^2$
($H_{39,2}$ with the orbit of $\alpha_5$; $H_{39,3}$
with the orbit of $\alpha_3$ or $\alpha_5$);
$16\aaa_1$ with $A_S\cong (\bz/4\bz)^2\times (\bz/2\bz)^2$
($H_{39,1}$ with the orbit of $\alpha_2$).

\medskip

{\bf n=34}, $H\cong {\mathfrak S}_4$ ($|H|=24$, $i=12$):
$\rk N_H=17$ and $(N_H)^\ast/N_H\cong
(\bz/12\bz)^2\times \bz/4\bz$.
{\it Degenerations of codimension $1$:}
$\aaa_1$ with $A_S\cong (\bz/12\bz)^2\times \bz/4\bz\times \bz/2\bz$
($H_{34,1}$ with the orbit of $\alpha_5$, $\alpha_6$,
$\alpha_{17}$ or $\alpha_{21}$;
$H_{34,2}$ with the orbit of $\alpha_6$, $\alpha_{21}$ or $\alpha_{24}$;
$H_{34,3}$ with the orbit of $\alpha_5$ or $\alpha_{21}$;
$H_{34,4}$ with the orbit of $\alpha_{10}$ or $\alpha_{24}$);
$2\aaa_1$ with $A_S\cong (\bz/12\bz)^2\times (\bz/2\bz)^2$
($H_{34,2}$ with the orbit of $\alpha_4$;
$H_{34,3}$ with the orbit of $\alpha_{4}$);
$3\aaa_1$ with $A_S\cong \bz/12\bz\times (\bz/4\bz)^2\times \bz/2\bz$
($H_{34,2}$ with the orbit of $\alpha_2$;
$H_{34,4}$ with the orbit of $\alpha_{2}$ or $\alpha_{3}$);
$4\aaa_1$ with $A_S\cong \bz/24\bz\times \bz/12\bz$
($H_{34,2}$ with the orbit of $\alpha_1$;
$H_{34,3}$ with the orbit of $\alpha_1$ or $\alpha_{10}$;
$H_{34,4}$ with the orbit of $\alpha_{1}$ or $\alpha_{5}$);
$6\aaa_1$ with $A_S\cong \bz/12\bz\times \bz/4\bz\times (\bz/2\bz)^2$
($H_{34,1}$ with the orbit of $\alpha_{1}$ or $\alpha_{8}$;
$H_{34,3}$ with the orbit of $\alpha_2$ or $\alpha_{3}$);
$8\aaa_1$ with $A_S\cong (\bz/12\bz)^2$
($H_{34,1}$ with the orbit of $\alpha_2$;
$H_{34,4}$ with the orbit of $\alpha_{4}$);
$12\aaa_1$ with $A_S\cong \bz/24\bz\times \bz/4\bz$
($H_{34,2}$ with the orbit of $\alpha_3$).

\medskip

{\bf n=33,} $H\cong C_7\rtimes C_3$ ($|H|=21$, $i=1$):
$\rk N_H=18$ and $(N_H)^\ast/N_H\cong (\bz/7\bz)^3$.
{\it Degenerations:}
$7\aaa_1$ with $A_S\cong \bz/14\bz\times \bz/7\bz$
($H_{33,1}$ with the orbit of $\alpha_1$, $\alpha_3$ or $\alpha_{11}$).

\medskip

{\bf n=32,} $H\cong Hol(C_5)$ ($|H|=20$, $i=3$):
$\rk N_H=18$ and
$(N_H)^\ast/N_H \cong (\bz/10\bz)^2 \times \bz/5\bz$.
{\it Degenerations:}
$2\aaa_1$ with $A_S\cong \bz/20\bz\times (\bz/5\bz)^2$
($H_{32,1}$ with the orbit of $\alpha_{19}$);
$5\aaa_1$ with $A_S\cong (\bz/10\bz)^2\times \bz/2\bz$
($H_{32,1}$ with the orbit of $\alpha_3$ or $\alpha_{4}$);
$10\aaa_1$ with $A_S\cong \bz/20\bz\times \bz/5\bz$
($H_{32,1}$ with the orbit of $\alpha_2$).

\medskip

{\bf n=30,} $H\cong {\mathfrak A}_{3,3}$ ($|H|=18$, $i=4$):
$\rk N_H=16$ and
$(N_H)^\ast/N_H\cong \bz/9\bz\times (\bz/3\bz)^4$.
{\it Degenerations of codimension $1$:}
$3\aaa_1$ with $A_S\cong \bz/18\bz\times (\bz/3\bz)^3$
($H_{30,1}$ with the orbit of $\alpha_4$, $\alpha_5$,
$\alpha_{11}$ or $\alpha_{15}$);
$9\aaa_1$ with $A_S\cong \bz/6\bz\times (\bz/3\bz)^3$
($H_{30,1}$ with the orbit of $\alpha_1$).

\medskip

{\bf n=26,} $H\cong SD_{16}$ ($|H|=16$, $i=8$):
$\rk N_H=18$ and
$(N_H)^\ast/N_H\cong (\bz/8\bz)^2\times
\bz/4\bz\times \bz/2\bz$.
{\it Degenerations:}
$8\aaa_1$ with $A_S\cong \bz/16\bz\times \bz/4\bz\times \bz/2\bz$
($H_{26,1}$ with the orbit of $\alpha_1$ or $\alpha_3$).

\medskip

{\bf n=22,} $H\cong C_2\times D_8$ ($|H|=16$, $i=11$):
$\rk N_H=16$, $(N_H)^\ast/N_H\cong
(\bz/4\bz)^4\times (\bz/2\bz)^2$.
{\it Degenerations of codimension $1$:}
$2\aaa_1$ with $A_S\cong (\bz/4\bz)^5$
($H_{22,2}$ with the orbit of $\alpha_{12}$ or $\alpha_{16}$;
$H_{22,3}$ with the orbit of
$\alpha_8$, $\alpha_{11}$, $\alpha_{12}$ or $\alpha_{16}$);
$4\aaa_1$ with
$A_S\cong \bz/8\bz\times (\bz/4\bz)^2\times (\bz/2\bz)^2$
($H_{22,2}$ with the orbit of $\alpha_6$ or $\alpha_8$;
$H_{22,3}$ with the orbit of $\alpha_1$,
$\alpha_2$, $\alpha_4$ or $\alpha_6$).
$8\aaa_1$ with $A_S\cong (\bz/4\bz)^3\times (\bz/2\bz)^2$
($H_{22,1}$ with the orbit of $\alpha_2$ or  $\alpha_4$;
$H_{22,2}$ with the orbit of $\alpha_2$).

\medskip

{\bf n=21,} $H\cong C_2^4$ ($|H|=16$, $i=14$):
$\rk N_H=15$, $(N_H)^\ast/N_H \cong
\bz/8\bz\times (\bz/2\bz)^6$.
{\it Degenerations of codimension $1$:}
$4\aaa_1$ with $A_S\cong (\bz/4\bz)^2\times (\bz/2\bz)^4$
($H_{21,1}$ with the orbit of $\alpha_{1}$, $\alpha_{2}$,
$\alpha_3$, $\alpha_8$ or $\alpha_{12}$);
$16\aaa_1$ with $A_S\cong (\bz/2\bz)^6$
($H_{21,2}$ with the orbit of $\alpha_1$).

\medskip

{\bf n=18,} $H\cong D_{12}$ ($|H|=12$, $i=4$):
$\rk N_H=16$ and $(N_H)^\ast/N_H\cong (\bz/6\bz)^4$.
{\it Degenerations of codimension $1$:}
$\aaa_1$ with $A_S\cong (\bz/6\bz)^4\times \bz/2\bz$
($H_{18,1}$ with the orbit of $\alpha_{1}$ or $\alpha_{6}$);
$2\aaa_1$ with $A_S\cong \bz/12\bz\times (\bz/6\bz)^2\times \bz/3\bz$
($H_{18,1}$ with the orbit of $\alpha_{15}$ or $\alpha_{16}$);
$3\aaa_1$ with $A_S\cong (\bz/6\bz)^3\times (\bz/2\bz)^2$
($H_{18,1}$ with the orbit of $\alpha_{4}$ or  $\alpha_{5}$);
$6\aaa_1$ with $A_S\cong \bz/12\bz\times (\bz/6\bz)^2$
($H_{18,1}$ with the orbit of $\alpha_2$ or $\alpha_8$).

\medskip

{\bf n=17,} $H\cong {\mathfrak A}_4$ ($|H|=12$, $i=3$):
$\rk N_H=16$ and $(N_H)^\ast/N_H\cong (\bz/12\bz)^2\times (\bz/2\bz)^2$.
{\it Degenerations of codimension $1$:}
$\aaa_1$ with $A_S\cong (\bz/12\bz)^2\times (\bz/2\bz)^3$
($H_{17,1}$ with the orbit of $\alpha_{5}$, $\alpha_{15}$,
$\alpha_{17}$, $\alpha_{21}$  or $\alpha_{24}$;
$H_{17,2}$ with the orbit of
$\alpha_{5}$, $\alpha_{6}$,
$\alpha_{10}$ or $\alpha_{17}$;
$H_{17,3}$ with the orbit of
$\alpha_{10}$ or $\alpha_{15}$);
$3\aaa_1$ with $A_S\cong \bz/12\bz\times \bz/4\bz\times (\bz/2\bz)^3$
($H_{17,1}$ with the orbit of $\alpha_{1}$;
$H_{17,3}$ with the orbit of $\alpha_1$ or $\alpha_6$);
$4\aaa_1$ with $A_S\cong \bz/24\bz\times \bz/6\bz\times \bz/2\bz$
($H_{17,1}$ with the orbit of $\alpha_{2}$;
$H_{17,2}$ with the orbit of $\alpha_2$ or $\alpha_4$;
$H_{17,3}$ with the orbit of $\alpha_2$, $\alpha_3$, $\alpha_8$ or
$\alpha_{11}$);
$6\aaa_1$ with $A_S\cong \bz/12\bz\times (\bz/4\bz)^2$
($H_{17,2}$ with the orbit of $\alpha_2$ or $\alpha_4$);
$12\aaa_1$ with $A_S\cong \bz/24\bz\times (\bz/2\bz)^2$
($H_{17,1}$ with the orbit of $\alpha_3$).

\medskip

{\bf n=16,} $H\cong D_{10}$ ($|H|=10$, $i=1$):
$\rk N_H=16$ and $(N_H)^\ast/N_H\cong (\bz/5\bz)^4$.
{\it Degenerations of codimension $1$:}
$\aaa_1$ with $A_S\cong \bz/10\bz\times (\bz/5\bz)^3$
($H_{16,1}$ with the orbit of $\alpha_1$, $\alpha_6$,
$\alpha_{19}$ or $\alpha_{22}$);
$5\aaa_1$ with $A_S\cong \bz/10\bz\times (\bz/5\bz)^2$
($H_{16,1}$ with the orbit of $\alpha_2$, $\alpha_3$, $\alpha_4$
or $\alpha_5$).

\medskip

{\bf n=12,} $H\cong Q_8$ ($|H|=8$, $i=4$):
$\rk N_H=17$, $(N_H)^\ast/N_H \cong
(\bz/8\bz)^2\times (\bz/2\bz)^3$.
{\it Degenerations of codimension $1$:}
$8\aaa_1$ with $A_S\cong \bz/16\bz\times (\bz/2\bz)^3$
($H_{12,1}$ with the orbit of $\alpha_1$ or $\alpha_{4}$).

\medskip

{\bf n=10,} $H\cong D_8$ ($|H|=8$, $i=3$):
$\rk N_H=15$ and $(N_H)^\ast/N_H\cong
(\bz/4\bz)^5$.
{\it Degenerations of codimension $1$:}
$\aaa_1$ with $A_S\cong (\bz/4\bz)^5\times \bz/2\bz$
($H_{10,1}$ with the orbit of $\alpha_{1}$, $\alpha_{5}$,
$\alpha_{6}$ or $\alpha_{24}$;
$H_{10,2}$ with the orbit of
$\alpha_{1}$ or $\alpha_{6}$);
$2\aaa_1$ with $A_S\cong (\bz/4\bz)^4\times (\bz/2\bz)^2$
($H_{10,1}$ with the orbit of $\alpha_{10}$ or $\alpha_{15}$;
$H_{10,2}$ with the orbit of $\alpha_2$, $\alpha_{10}$ or
$\alpha_{16}$);
$4\aaa_1$ with $A_S\cong \bz/8\bz\times (\bz/4\bz)^3$
($H_{10,1}$ with the orbit of $\alpha_{3}$ or $\alpha_7$;
$H_{10,2}$ with the orbit of $\alpha_3$, $\alpha_4$, $\alpha_5$
or $\alpha_7$);
$8\aaa_1$ with $A_S\cong (\bz/4\bz)^4$
($H_{10,1}$ with the orbit of $\alpha_2$).

\medskip

{\bf n=9,} $H\cong C_2^3$ ($|H|=8$, $i=5$):
$\rk N_H=14$, $(N_H)^\ast/N_H \cong
(\bz/4\bz)^2\times (\bz/2\bz)^6$.
{\it Degenerations of codimension $1$}:
$2\aaa_1$ with $A_S\cong (\bz/4\bz)^3\times (\bz/2\bz)^4$
($H_{9,1}$ with the orbit of $\alpha_5$ or $\alpha_{6}$;
$H_{9,3}$ with the orbit of $\alpha_6$,  $\alpha_8$, $\alpha_{10}$,
$\alpha_{11}$, $\alpha_{12}$, $\alpha_{15}$ or $\alpha_{16}$);
$4\aaa_1$ with $A_S\cong \bz/8\bz\times (\bz/2\bz)^6$
($H_{9,1}$ with the orbit of $\alpha_2$,
$\alpha_3$, $\alpha_4$ or $\alpha_7$;
$H_{9,4}$ with the orbit of $\alpha_2$ or $\alpha_4$);
$8\aaa_1$ with $A_S\cong \bz/4\bz\times (\bz/2\bz)^6$
($H_{9,2}$ with the orbit of $\alpha_2$ or $\alpha_4$;
$H_{9,3}$ with the orbit of $\alpha_{2}$).

\medskip

{\bf n=6,} $H\cong D_6$ ($|H|=6$, $i=1$):
$\rk N_H=14$ and $(N_H)^\ast/N_H\cong (\bz/6\bz)^2\times (\bz/3\bz)^3$.
{\it Degenerations of codimension $1$}:
$\aaa_1$ with $A_S\cong (\bz/6\bz)^3\times (\bz/3\bz)^2$
($H_{6,1}$ with the orbit of $\alpha_1$,  $\alpha_5$, $\alpha_{6}$
or $\alpha_{24}$);
$2\aaa_1$ with $A_S\cong \bz/12\bz\times (\bz/3\bz)^4$
($H_{6,1}$ with the orbit of $\alpha_{16}$);
$3\aaa_1$ with $A_S\cong (\bz/6\bz)^3\times \bz/3\bz$
($H_{6,1}$ with the orbit of $\alpha_3$, $\alpha_7$,
$\alpha_{10}$ or $\alpha_{11}$);
$6\aaa_1$ with $A_S\cong \bz/12\bz\times (\bz/3\bz)^3$
($H_{6,1}$ with the orbit of $\alpha_2$).

\medskip

{\bf n=4,} $H\cong C_4$ ($|H|=4$, $i=1$):
$\rk N_H=14$ and $(N_H)^\ast/N_H \cong
(\bz/4\bz)^4\times (\bz/2\bz)^2$.
{\it Degenerations of codimension $1$}:
$\aaa_1$ with $A_S\cong (\bz/4\bz)^4\times (\bz/2\bz)^3$
($H_{4,1}$ with the orbit of $\alpha_1$,  $\alpha_4$, $\alpha_{17}$
or $\alpha_{18}$);
$2\aaa_1$ with $A_S\cong (\bz/4\bz)^5$
($H_{4,1}$ with the orbit of $\alpha_8$ or $\alpha_{12}$);
$4\aaa_1$ with $A_S\cong \bz/8\bz \times (\bz/4\bz)^2\times  (\bz/2\bz)^2$
($H_{4,1}$ with the orbit of $\alpha_2$, $\alpha_3$,
$\alpha_{7}$ or $\alpha_{11}$);

\medskip

{\bf n=3,} $H\cong C_2^2$ ($|H|=4$, $i=2$):
$\rk N_H=12$ and $(N_H)^\ast/N_H\cong (\bz/4\bz)^2\times (\bz/2\bz)^6$.
{\it Degenerations of codimension $1$}:
$\aaa_1$ with $A_S\cong (\bz/4\bz)^2\times (\bz/2\bz)^7$
($H_{3,1}$ with the orbit of $\alpha_1$, $\alpha_5$, $\alpha_{6}$,
$\alpha_{10}$, $\alpha_{15}$, $\alpha_{17}$, $\alpha_{21}$ or
$\alpha_{24}$;
$H_{3,2}$ with the orbit of $\alpha_1$, $\alpha_5$, $\alpha_{6}$
or  $\alpha_{21}$);
$2\aaa_1$ with $A_S\cong (\bz/4\bz)^3\times (\bz/2\bz)^4$
($H_{3,2}$ with the orbit of $\alpha_8$,
$\alpha_{10}$, $\alpha_{11}$, $\alpha_{12}$, $\alpha_{15}$ or
$\alpha_{16}$;
$H_{3,3}$ with the orbit of $\alpha_1$, $\alpha_2$,
$\alpha_3$, $\alpha_4$, $\alpha_5$, $\alpha_6$, $\alpha_7$, $\alpha_8$,
$\alpha_{11}$,  $\alpha_{12}$, $\alpha_{15}$ or  $\alpha_{16}$);
$4\aaa_1$ with $A_S\cong \bz/8\bz \times (\bz/2\bz)^6$
($H_{3,1}$ with the orbit of $\alpha_2$, $\alpha_{3}$, $\alpha_{4}$
or $\alpha_7$;
$H_{3,2}$ with the orbit of $\alpha_{2}$ or $\alpha_4$).

\medskip

{\bf n=2,} $H\cong C_3$ ($|H|=3$, $i=1$):
$\rk N_H=12$ and $(N_H)^\ast/N_H\cong (\bz/3\bz)^6$.
{\it Degenerations of codimension $1$}:
$\aaa_1$ with $A_S\cong \bz/6\bz\times (\bz/3\bz)^5$
($H_{2,1}$ with the orbit of $\alpha_3$,  $\alpha_4$, $\alpha_{14}$,
$\alpha_{17}$, $\alpha_{21}$ or $\alpha_{24}$);
$3\aaa_1$ with $A_S\cong \bz/6\bz \times (\bz/3\bz)^4$
($H_{2,1}$ with the orbit of $\alpha_1$, $\alpha_{2}$,
$\alpha_{5}$, $\alpha_{6}$, $\alpha_{12}$ or $\alpha_{15}$).

\medskip

{\bf n=1,} $H\cong C_2$ ($|H|=2$, $i=1$):
$\rk N_H=8$ and $(N_H)^\ast/N_H\cong (\bz/2\bz)^8$.
{\it Degenerations of codimension $1$}:
$\aaa_1$ with $A_S\cong (\bz/2\bz)^9$
($H_{1,1}$ with the orbit of $\alpha_1$,  $\alpha_4$, $\alpha_{8}$,
$\alpha_{12}$, $\alpha_{15}$,  $\alpha_{17}$, $\alpha_{18}$ or $\alpha_{21}$);
$2\aaa_1$ with $A_S\cong \bz/4\bz\times (\bz/2\bz)^6$
($H_{1,1}$ with the orbit of $\alpha_2$, $\alpha_{3}$,
$\alpha_{5}$, $\alpha_{6}$, $\alpha_{7}$, $\alpha_{10}$,
$\alpha_{11}$  or $\alpha_{14}$).


\vskip 1cm

{\bf Table 22. All possible markings by $N_{22}=N(12A_2)$ of
degenerations of codimension one of
 non-trivial finite symplectic automorphism groups of
 K\"ahlerian K3 surfaces:}

{\bf n=79,} $H\cong {\mathfrak A}_{6}$ ($|H|=360$, $i=118$):
$\rk N_H=19$.
{\it No degenerations}.

\medskip

{\bf n=70,} $H\cong {\mathfrak S}_5$ ($|H|=120$, $i=34$):
$\rk N_H=19$. {\it No degenerations}.

\medskip

{\bf n=63,} $H\cong M_9$ ($|H|=72$, $i=41$):
$\rk N_H=19$. {\it No degenerations}.

\medskip

{\bf n=62,} $H\cong N_{72}$ ($|H|=72$, $i=40$):
$\rk N_H=19$. {\it No degenerations}.

\medskip

{\bf n=55,} $H\cong {\mathfrak A}_5$ ($|H|=60$, $i=5$):
$\rk N_H=18$ and $(N_H)^\ast/N_H \cong\bz/30\bz\times \bz/10\bz$.
{\it Degenerations}:
$\aaa_1$ with $A_S\cong \bz/30\bz \times \bz/10\bz \times \bz/2\bz$
($H_{55,1}$ with the orbit of $\alpha_{1,10}$ or $\alpha_{2,10}$;
$H_{55,2}$ with the orbit of $\alpha_{1,6}$,
$\alpha_{2,6}$, $\alpha_{1,8}$ or $\alpha_{2,8}$);
$5\aaa_1$ with $A_S\cong \bz/30\bz \times (\bz/2\bz)^2$
($H_{55,1}$ with the orbit of $\alpha_{1,6}$ or $\alpha_{2,6}$);
$6\aaa_1$ with $A_S\cong \bz/20\bz \times \bz/5\bz$
($H_{55,1}$ with the orbit of $\alpha_{1,1}$ or $\alpha_{2,1}$);
$10\aaa_1$ with $A_S\cong \bz/60\bz$
($H_{55,2}$ with the orbit of $\alpha_{1,1}$ or $\alpha_{2,1}$).

\medskip

{\bf n=54,} $H\cong T_{48}$ ($|H|=48$, $i=29$):
$\rk N_H=19$.
{\it No degenerations}.

\medskip

{\bf n=48,} $H\cong {\mathfrak S}_{3,3}$ ($|H|=36$, $i=10$):
$\rk N_H=18$, $(N_H)^\ast/N_H \cong
\bz/18\bz\times
\bz/6\bz\times (\bz/3\bz)^2$.
{\it Degenerations}:
$3\aaa_1$ with $A_S\cong \bz/18\bz \times (\bz/6\bz)^2$
($H_{48,2}$ with the orbit of $\alpha_{1,6}$, $\alpha_{2,6}$,
$\alpha_{1,7}$ or $\alpha_{2,7}$);
$6\aaa_1$ with $A_S\cong \bz/36\bz \times (\bz/3\bz)^2$
($H_{48,2}$ with the orbit of $\alpha_{1,1}$ or $\alpha_{1,2}$);
$9\aaa_1$ with $A_S\cong (\bz/6\bz)^3$
($H_{48,1}$ with the orbit of $\alpha_{1,1}$ or $\alpha_{2,1}$).

\medskip

{\bf n=46,} $H\cong 3^2 C_4$ ($|H|=36$, $i=9$):
$\rk N_H=18$ and $(N_H)^\ast/N_H\cong \bz/18\bz\times \bz/6\bz\times \bz/3\bz$.
{\it Degenerations}:
$9\aaa_1$ with $A_S\cong (\bz/6\bz)^2 \times \bz/2\bz$
($H_{46,1}$ with the orbit of $\alpha_{1,1}$ or  $\alpha_{2,1}$);
$9\aaa_2$ with $A_S\cong \bz/6\bz \times \bz/3\bz$
($H_{46,2}$ with the orbit of $\alpha_{1,1}$).

\medskip

{\bf n=34,} $H\cong {\mathfrak S}_4$ ($|H|=24$, $i=12$):
$\rk N_H=17$ and $(N_H)^\ast/N_H\cong (\bz/12\bz)^2\times \bz/4\bz$.
{\it Degenerations of codimension $1$}:
$4\aaa_1$ with $A_S\cong \bz/24\bz \times \bz/12\bz$
($H_{34,1}$ with the orbit of $\alpha_{1,1}$ or  $\alpha_{2,1}$);
$6\aaa_2$ with $A_S\cong \bz/12\bz \times \bz/4\bz$
($H_{34,1}$ with the orbit of $\alpha_{1,2}$).

\medskip

{\bf n=32,} $H\cong Hol(C_5)$ ($|H|=20$, $i=3$):
$\rk N_H=18$ and $(N_H)^\ast/N_H\linebreak
\cong (\bz/10\bz)^2\times \bz/5\bz$.
{\it Degenerations}:
$5\aaa_1$ with $A_S\cong (\bz/10\bz)^2 \times \bz/2\bz$
($H_{32,1}$ with the orbit of $\alpha_{1,3}$ or  $\alpha_{2,3}$);
$5\aaa_2$ with $A_S\cong \bz/10\bz \times \bz/5\bz$
($H_{32,1}$ with the orbit of $\alpha_{1,6}$).

\medskip

{\bf n=30,} $H\cong {\mathfrak A}_{3,3}$ ($|H|=18$, $i=4$):
$\rk N_H=16$ and $(N_H)^\ast/N_H\cong \bz/9\bz\times (\bz/3\bz)^4$.
{\it Degenerations of codimension $1$}:
$3\aaa_1$ with $A_S\cong \bz/18\bz \times (\bz/3\bz)^3$
($H_{30,2}$ with the orbit of
$\alpha_{1,1}$, $\alpha_{2,1}$, $\alpha_{1,2}$,
$\alpha_{2,2}$, $\alpha_{1,6}$, $\alpha_{2,6}$,
$\alpha_{1,7}$ or $\alpha_{2,7}$);
$9\aaa_1$ with $A_S\cong \bz/6\bz \times (\bz/3\bz)^3$
($H_{30,1}$ with the orbit of $\alpha_{1,1}$ or $\alpha_{2,1}$).

\medskip

{\bf n=26,} $H\cong SD_{16}$ ($|H|=16$, $i=8$):
$\rk N_H=18$ and
$(N_H)^\ast/N_H\cong (\bz/8\bz)^2\times \bz/4\bz\times \bz/2\bz$.
{\it Degenerations}:
$2\aaa_2$ with $A_S\cong (\bz/8\bz)^2 \times \bz/2\bz$
($H_{26,1}$ with the orbit of $\alpha_{1,3}$);
$8\aaa_1$ with $A_S\cong \bz/16\bz \times \bz/4\bz\times \bz/2\bz$
($H_{26,1}$ with the orbit of $\alpha_{1,5}$ or $\alpha_{2,5}$).

\medskip

{\bf n=18,} $H\cong D_{12}$ ($|H|=12$, $i=4$):
$\rk N_H=16$ and $(N_H)^\ast/N_H\cong (\bz/6\bz)^4$.
{\it Degenerations of codimension $1$}:
$\aaa_1$ with $A_S\cong (\bz/6\bz)^4 \times \bz/2\bz$
($H_{18,1}$ with the orbit of
$\alpha_{1,1}$ or $\alpha_{2,1}$);
$2\aaa_1$ with $A_S\cong \bz/12\bz \times (\bz/6\bz)^2\times\bz/3\bz$
($H_{18,1}$ with the orbit of $\alpha_{1,7}$ or $\alpha_{2,7}$);
$3\aaa_1$ with $A_S\cong (\bz/6\bz)^3 \times (\bz/2\bz)^2$
($H_{18,1}$ with the orbit of $\alpha_{1,2}$ or $\alpha_{2,2}$);
$6\aaa_1$ with $A_S\cong \bz/12\bz \times (\bz/6\bz)^2$
($H_{18,1}$ with the orbit of $\alpha_{1,5}$ or $\alpha_{2,5}$).

\medskip

{\bf n=17,} $H={\mathfrak A}_4$ ($|H|=12$, $i=3$):
$\rk N_H=16$ and $(N_H)^\ast/N_H\cong (\bz/12\bz)^2\times (\bz/2\bz)^2$.
{\it Degenerations of codimension $1$}:
$\aaa_1$ with $A_S\cong (\bz/12\bz)^2 \times (\bz/2\bz)^3$
($H_{17,1}$ with the orbit of
$\alpha_{1,1}$, $\alpha_{2,1}$, $\alpha_{1,2}$ or $\alpha_{2,2}$);
$4\aaa_1$ with $A_S\cong \bz/24\bz \times \bz/6\bz \times\bz/2\bz$
($H_{17,1}$ with the orbit of $\alpha_{1,3}$ or $\alpha_{2,3}$);
$6\aaa_1$ with $A_S\cong \bz/12\bz \times (\bz/4\bz)^2$
($H_{17,1}$ with the orbit of $\alpha_{1,5}$ or $\alpha_{2,5}$).

\medskip

{\bf n=16,} $H\cong D_{10}$ ($|H|=10$, $i=1$):
$\rk N_H=16$ and $(N_H)^\ast/N_H\cong (\bz/5\bz)^4$.
{\it Degenerations of codimension $1$}:
$\aaa_1$ with $A_S\cong \bz/10\bz \times (\bz/5\bz)^3$
($H_{16,1}$ with the orbit of
$\alpha_{1,1}$, $\alpha_{2,1}$, $\alpha_{1,2}$ or $\alpha_{2,2}$);
$5\aaa_1$ with $A_S\cong \bz/10\bz \times (\bz/5\bz)^2$
($H_{16,1}$ with the orbit of $\alpha_{1,3}$, $\alpha_{2,3}$,
$\alpha_{1,5}$ or $\alpha_{2,5}$).

\medskip

{\bf n=12,} $H\cong Q_8$ ($|H|=8$, $i=4$):
$\rk N_H=17$, $(N_H)^\ast/N_H \cong (\bz/8\bz)^2\times (\bz/2\bz)^3$.
{\it Degenerations of codimension $1$}:
$\aaa_2$ with $A_S\cong (\bz/8\bz)^2 \times (\bz/2\bz)^2$
($H_{12,1}$ with the orbit of
$\alpha_{1,5}$, $\alpha_{1,8}$ or $\alpha_{1,11}$);
$8\aaa_1$ with $A_S\cong \bz/16\bz \times (\bz/2\bz)^3$
($H_{12,1}$ with the orbit of $\alpha_{1,1}$ or $\alpha_{2,1}$).

\medskip

{\bf n=10,} $H\cong D_8$ ($|H|=8$, $i=3$):
$\rk N_H=15$ and $(N_H)^\ast/N_H \cong (\bz/4\bz)^5$.
{\it Degenerations of codimension $1$}:
$\aaa_1$ with $A_S\cong (\bz/4\bz)^5 \times \bz/2\bz$
($H_{10,1}$ with the orbit of
$\alpha_{1,1}$, $\alpha_{2,1}$, $\alpha_{1,2}$ or $\alpha_{2,2}$);
$4\aaa_1$ with $A_S\cong \bz/8\bz \times (\bz/4\bz)^3$
($H_{10,1}$ with the orbit of $\alpha_{1,3}$, $\alpha_{2,3}$,
$\alpha_{1,5}$ or $\alpha_{2,5}$);
$2\aaa_2$ with $A_S\cong (\bz/4\bz)^4$
($H_{10,1}$ with the orbit of $\alpha_{1,9}$).

\medskip

{\bf n=6,} $H\cong D_6$ ($|H|=6$, $i=1$):
$\rk N_H=14$, $(N_H)^\ast/N_H \cong (\bz/6\bz)^2\times (\bz/3\bz)^3$.
{\it Degenerations of codimension $1$}:
$\aaa_1$ with $A_S\cong (\bz/6\bz)^3 \times (\bz/3\bz)^2$
($H_{6,1}$ with the orbit of
$\alpha_{1,8}$ or $\alpha_{2,8}$;
$H_{6,2}$ with the orbit of $\alpha_{1,6}$, $\alpha_{2,6}$,
$\alpha_{1,8}$, $\alpha_{2,8}$, $\alpha_{1,10}$ or $\alpha_{2,10}$);
$2\aaa_1$ with $A_S\cong \bz/12\bz \times (\bz/3\bz)^4$
($H_{6,1}$ with the orbit of $\alpha_{1,6}$ or $\alpha_{2,6}$);
$3\aaa_1$ with $A_S\cong (\bz/6\bz)^3 \times \bz/3\bz$
($H_{6,2}$ with the orbit of $\alpha_{1,1}$ or $\alpha_{2,1}$);
$6\aaa_1$ with $A_S\cong \bz/12\bz \times (\bz/3\bz)^3$
($H_{6,2}$ with the orbit of $\alpha_{1,2}$ or $\alpha_{2,2}$).

\medskip

{\bf n=4,} $H\cong C_4$ ($|H|=4$, $i=1$):
$\rk N_H=14$ and $(N_H)^\ast/N_H\cong (\bz/4\bz)^4\times (\bz/2\bz)^2$.
{\it Degenerations of codimension $1$}:
$\aaa_1$ with $A_S\cong (\bz/4\bz)^4 \times (\bz/2\bz)^3$
($H_{4,1}$ with the orbit of
$\alpha_{1,2}$, $\alpha_{2,2}$, $\alpha_{1,5}$ or $\alpha_{2,5}$);
$\aaa_2$ with $A_S\cong (\bz/4\bz)^4 \times \bz/2\bz$
($H_{4,1}$ with the orbit of $\alpha_{1,1}$ or $\alpha_{1,11}$);
$4\aaa_1$ with $A_S\cong \bz/8\bz\times (\bz/4\bz)^2\times (\bz/2\bz)^2$
($H_{4,1}$ with the orbit of $\alpha_{1,3}$, $\alpha_{2,3}$,
$\alpha_{1,6}$ or $\alpha_{2,6}$).

\medskip

{\bf n=3,} $H\cong C_2^2$ ($|H|=4$, $i=2$):
$\rk N_H=12$ and $(N_H)^\ast/N_H\cong (\bz/4\bz)^2\times (\bz/2\bz)^6$.
{\it Degenerations of codimension $1$}:
$\aaa_1$ with $A_S\cong (\bz/4\bz)^2 \times (\bz/2\bz)^7$
($H_{3,1}$ with the orbit of
$\alpha_{1,1}$, $\alpha_{2,1}$, $\alpha_{1,2}$ or $\alpha_{2,2}$);
$2\aaa_1$ with $A_S\cong (\bz/4\bz)^3 \times (\bz/2\bz)^4$
($H_{1,3}$ with the orbit of
$\alpha_{1,3}$, $\alpha_{2,3}$, $\alpha_{1,5}$, $\alpha_{2,5}$,
$\alpha_{1,7}$ or $\alpha_{2,7}$);
$4\aaa_1$ with $A_S\cong \bz/8\bz\times (\bz/2\bz)^6$
($H_{3,1}$ with the orbit of $\alpha_{1,6}$ or $\alpha_{2,6}$).

{\bf n=2,} $H\cong C_3$ ($|H|=3$, $i=1$):
$\rk N_H=12$ and $(N_H)^\ast/N_H\cong (\bz/3\bz)^6$.
{\it Degenerations of codimension $1$}:
$\aaa_1$ with $A_S\cong \bz/6\bz \times (\bz/3\bz)^5$
($H_{2,1}$ with the orbit of
$\alpha_{1,6}$, $\alpha_{2,6}$, $\alpha_{1,8}$, $\alpha_{2,8}$,
$\alpha_{1,10}$ or $\alpha_{2,10}$);
$3\aaa_1$ with $A_S\cong \bz/6\bz \times (\bz/3\bz)^4$
($H_{2,1}$ with the orbit of
$\alpha_{1,1}$, $\alpha_{2,1}$, $\alpha_{1,2}$, $\alpha_{2,2}$,
$\alpha_{1,7}$ or $\alpha_{2,7}$).

\medskip

{\bf n=1,} $H\cong C_2$ ($|H|=2$, $i=1$):
$\rk N_H=8$ and $(N_H)^\ast/N_H\cong (\bz/2\bz)^8$.
{\it Degenerations of codimension $1$}:
$\aaa_1$ with $A_S\cong (\bz/2\bz)^9$
($H_{1,1}$ with the orbit of
$\alpha_{1,1}$, $\alpha_{2,1}$, $\alpha_{1,2}$, $\alpha_{2,2}$,
$\alpha_{1,5}$, $\alpha_{2,5}$, $\alpha_{1,11}$ or $\alpha_{2,11}$);
$2\aaa_1$ with $A_S\cong \bz/4\bz \times (\bz/2\bz)^6$
($H_{1,1}$ with the orbit of
$\alpha_{1,3}$, $\alpha_{2,3}$, $\alpha_{1,6}$, $\alpha_{2,6}$,
$\alpha_{1,7}$, $\alpha_{2,7}$, $\alpha_{1,9}$ or $\alpha_{2,9}$).


\vskip1cm

{\bf Table 21. All possible markings by $N_{21}=N(8A_3)$ of degenerations of
codimension one of  non-trivial finite symplectic automorphism
groups of K\"ahlerian K3 surfaces:}

\medskip

{\bf n=74,} $H\cong L_2(7)$ ($|H|=168$, $i=42$):
$\rk N_H=19$. {\it No degenerations.}

\medskip

{\bf n=51,} $H\cong C_2\times {\mathfrak S}_4$ ($|H|=48$, $i=48$):
$\rk N_H=18$, $(N_H)^\ast/N_H\cong (\bz/12\bz)^2\times (\bz/2\bz)^2$.
{\it Degenerations:}
$2\aaa_1$ with $A_S\cong (\bz/12\bz)^2\times  \bz/4\bz$
($H_{51,1}$ with the orbit of $\alpha_{1,1}$ or $\alpha_{1,2}$;
$H_{51,2}$ with the orbit of $\alpha_{1,1}$);
$4\aaa_1$ with $A_S\cong \bz/24\bz \times \bz/6\bz \times \bz/2\bz$
($H_{51,2}$ with the orbit of $\alpha_{2,4}$);
$6\aaa_1$ with $A_S\cong \bz/12\bz \times (\bz/4\bz)^2$
($H_{51,2}$ with the orbit of $\alpha_{1,2}$);
$8\aaa_1$ with $A_S\cong \bz/12\bz\times  \bz/6\bz\times \bz/2\bz$
($H_{51,2}$ with the orbit of $\alpha_{1,4}$);
$12\aaa_1$ with $A_S\cong \bz/24\bz\times (\bz/2\bz)^2$
($H_{51,1}$ with the orbit of $\alpha_{1,3}$).

\medskip

{\bf n=34}, $H\cong {\mathfrak S}_4$ ($|H|=24$, $i=12$):
$\rk N_H=17$ and $(N_H)^\ast/N_H\cong
(\bz/12\bz)^2\times \bz/4\bz$.
{\it Degenerations of codimension $1$:}
$\aaa_1$ with $A_S\cong (\bz/12\bz)^2\times \bz/4\bz\times \bz/2\bz$
($H_{34,1}$ with the orbit of $\alpha_{2,1}$, $\alpha_{1,2}$,
$\alpha_{2,2}$ or $\alpha_{3,2}$;
$H_{34,2}$ with the orbit of $\alpha_{2,1}$;
$H_{34,3}$ with the orbit of $\alpha_{1,2}$, $\alpha_{2,2}$ or
$\alpha_{3,2}$);
$2\aaa_1$ with $A_S\cong (\bz/12\bz)^2\times (\bz/2\bz)^2$
($H_{34,1}$ with the orbit of $\alpha_{1,1}$;
$H_{34,2}$ with the orbit of $\alpha_{1,1}$);
$3\aaa_1$ with $A_S\cong \bz/12\bz\times (\bz/4\bz)^2\times \bz/2\bz$
($H_{34,2}$ with the orbit of $\alpha_{2,2}$;
$H_{34,3}$ with the orbit of $\alpha_{2,1}$);
$4\aaa_1$ with $A_S\cong \bz/24\bz\times \bz/12\bz$
($H_{34,2}$ with the orbit of $\alpha_{1,4}$, $\alpha_{2,4}$ or
$\alpha_{3,4}$;
$H_{34,3}$ with the orbit of $\alpha_{2,4}$);
$6\aaa_1$ with $A_S\cong \bz/12\bz\times \bz/4\bz\times (\bz/2\bz)^2$
($H_{34,1}$ with the orbit of $\alpha_{2,3}$;
$H_{34,2}$ with the orbit of $\alpha_{1,2}$;
$H_{34,3}$ with the orbit of $\alpha_{1,1}$);
$8\aaa_1$ with $A_S\cong (\bz/12\bz)^2$
($H_{34,3}$ with the orbit of $\alpha_{1,4}$);
$12\aaa_1$ with $A_S\cong \bz/24\bz\times \bz/4\bz$
($H_{34,1}$ with the orbit of $\alpha_{1,3}$).

\medskip

{\bf n=33,} $H\cong C_7\rtimes C_3$ ($|H|=21$, $i=1$):
$\rk N_H=18$ and $(N_H)^\ast/N_H\cong (\bz/7\bz)^3$.
{\it Degenerations:}
$7\aaa_1$ with $A_S\cong \bz/14\bz\times \bz/7\bz$
($H_{33,1}$ with the orbit of $\alpha_{1,2}$, $\alpha_{2,2}$
or $\alpha_{3,2}$).

\medskip

{\bf n=22,} $H\cong C_2\times D_8$ ($|H|=16$, $i=11$):
$\rk N_H=16$, $(N_H)^\ast/N_H\cong
(\bz/4\bz)^4\times (\bz/2\bz)^2$.
{\it Degenerations of codimension $1$:}
$2\aaa_1$ with $A_S\cong (\bz/4\bz)^5$
($H_{22,1}$ with the orbit of $\alpha_{1,1}$ or $\alpha_{1,2}$);
$4\aaa_1$ with $A_S\cong \bz/8\bz\times (\bz/4\bz)^2\times (\bz/2\bz)^2$
($H_{22,1}$ with the orbit of $\alpha_{2,3}$ or $\alpha_{1,4}$);
$8\aaa_1$ with $A_S\cong (\bz/4\bz)^3\times (\bz/2\bz)^2$
($H_{22,1}$ with the orbit of $\alpha_{1,3}$).

\medskip

{\bf n=18,} $H\cong D_{12}$ ($|H|=12$, $i=4$):
$\rk N_H=16$ and $(N_H)^\ast/N_H\cong (\bz/6\bz)^4$.
{\it Degenerations of codimension $1$:}
$\aaa_1$ with $A_S\cong (\bz/6\bz)^4\times \bz/2\bz$
($H_{18,1}$ with the orbit of $\alpha_{2,1}$ or $\alpha_{2,2}$);
$2\aaa_1$ with $A_S\cong \bz/12\bz\times (\bz/6\bz)^2\times \bz/3\bz$
($H_{18,1}$ with the orbit of $\alpha_{1,1}$ or $\alpha_{1,2}$);
$3\aaa_1$ with $A_S\cong (\bz/6\bz)^3\times (\bz/2\bz)^2$
($H_{18,1}$ with the orbit of $\alpha_{2,3}$ or  $\alpha_{2,5}$);
$6\aaa_1$ with $A_S\cong \bz/12\bz\times (\bz/6\bz)^2$
($H_{18,1}$ with the orbit of $\alpha_{1,3}$ or $\alpha_{1,5}$).

\medskip

{\bf n=17,} $H\cong {\mathfrak A}_4$ ($|H|=12$, $i=3$):
$\rk N_H=16$ and $(N_H)^\ast/N_H\cong (\bz/12\bz)^2\times (\bz/2\bz)^2$.
{\it Degenerations of codimension $1$:}
$\aaa_1$ with $A_S\cong (\bz/12\bz)^2\times (\bz/2\bz)^3$
($H_{17,1}$ with the orbit of $\alpha_{1,1}$, $\alpha_{2,1}$,
$\alpha_{3,1}$, $\alpha_{1,5}$, $\alpha_{2,5}$ or $\alpha_{3,5}$;
$H_{17,2}$ with the orbit of $\alpha_{1,1}$, $\alpha_{2,1}$
or $\alpha_{3,1}$);
$3\aaa_1$ with $A_S\cong \bz/12\bz\times \bz/4\bz\times (\bz/2\bz)^3$
($H_{17,2}$ with the orbit of $\alpha_{2,2}$);
$4\aaa_1$ with $A_S\cong \bz/24\bz\times \bz/6\bz\times \bz/2\bz$
($H_{17,2}$ with the orbit of $\alpha_{1,3}$, $\alpha_{2,3}$ or
$\alpha_{3,3}$);
$6\aaa_1$ with $A_S\cong \bz/12\bz\times (\bz/4\bz)^2$
($H_{17,1}$ with the orbit of $\alpha_{2,2}$;
$H_{17,2}$ with the orbit of $\alpha_{1,2}$);
$12\aaa_1$ with $A_S\cong \bz/24\bz\times (\bz/2\bz)^2$
($H_{17,1}$ with the orbit of $\alpha_{1,2}$).

\medskip

{\bf n=10,} $H\cong D_8$ ($|H|=8$, $i=3$):
$\rk N_H=15$ and $(N_H)^\ast/N_H\cong
(\bz/4\bz)^5$.
{\it Degenerations of codimension $1$:}
$\aaa_1$ with $A_S\cong (\bz/4\bz)^5\times \bz/2\bz$
($H_{10,1}$ with the orbit of $\alpha_{1,1}$, $\alpha_{2,1}$,
$\alpha_{3,1}$ or $\alpha_{2,2}$;
$H_{10,2}$ with the orbit of
$\alpha_{2,1}$ or $\alpha_{2,2}$);
$2\aaa_1$ with $A_S\cong (\bz/4\bz)^4\times (\bz/2\bz)^2$
($H_{10,1}$ with the orbit of $\alpha_{1,2}$ or $\alpha_{2,4}$;
$H_{10,2}$ with the orbit of $\alpha_{1,1}$, $\alpha_{1,2}$ or
$\alpha_{2,4}$);
$4\aaa_1$ with $A_S\cong \bz/8\bz\times (\bz/4\bz)^3$
($H_{10,1}$ with the orbit of $\alpha_{2,3}$ or $\alpha_{1,4}$;
$H_{10,2}$ with the orbit of $\alpha_{1,3}$, $\alpha_{2,3}$,
$\alpha_{3,3}$
or $\alpha_{1,4}$);
$8\aaa_1$ with $A_S\cong (\bz/4\bz)^4$
($H_{10,1}$ with the orbit of $\alpha_{1,3}$).

\medskip

{\bf n=9,} $H\cong C_2^3$ ($|H|=8$, $i=5$):
$\rk N_H=14$, $(N_H)^\ast/N_H \cong
(\bz/4\bz)^2\times (\bz/2\bz)^6$.
{\it Degenerations of codimension $1$}:
$2\aaa_1$ with $A_S\cong (\bz/4\bz)^3\times (\bz/2\bz)^4$
($H_{9,1}$ with the orbit of $\alpha_{1,1}$, $\alpha_{1,2}$,
$\alpha_{1,4}$
or $\alpha_{1,5}$;
$H_{9,2}$ with the orbit of $\alpha_{1,1}$,  $\alpha_{2,2}$,
$\alpha_{2,3}$,
$\alpha_{1,5}$ or $\alpha_{2,7}$);
$4\aaa_1$ with $A_S\cong \bz/8\bz\times (\bz/2\bz)^6$
($H_{9,1}$ with the orbit of $\alpha_{2,3}$;
$H_{9,2}$ with the orbit of $\alpha_{1,2}$, $\alpha_{1,3}$
or $\alpha_{1,7}$);
$8\aaa_1$ with $A_S\cong \bz/4\bz\times (\bz/2\bz)^6$
($H_{9,1}$ with the orbit of $\alpha_{1,3}$).

\medskip

{\bf n=6,} $H\cong D_6$ ($|H|=6$, $i=1$):
$\rk N_H=14$ and $(N_H)^\ast/N_H\cong (\bz/6\bz)^2\times (\bz/3\bz)^3$.
{\it Degenerations of codimension $1$}:
$\aaa_1$ with $A_S\cong (\bz/6\bz)^3\times (\bz/3\bz)^2$
($H_{6,1}$ with the orbit of $\alpha_{1,1}$,  $\alpha_{2,1}$,
$\alpha_{3,1}$
or $\alpha_{2,2}$);
$2\aaa_1$ with with $A_S\cong \bz/12\bz\times (\bz/3\bz)^4$
($H_{6,1}$ with the orbit of $\alpha_{1,2}$);
$3\aaa_1$ with $A_S\cong (\bz/6\bz)^3\times \bz/3\bz$
($H_{6,1}$ with the orbit of $\alpha_{2,3}$, $\alpha_{1,5}$,
$\alpha_{2,5}$ or $\alpha_{3,5}$);
$6\aaa_1$ with $A_S\cong \bz/12\bz\times (\bz/3\bz)^3$
($H_{6,1}$ with the orbit of $\alpha_{1,3}$).

\medskip

{\bf n=4,} $H\cong C_4$ ($|H|=4$, $i=1$):
$\rk N_H=14$ and $(N_H)^\ast/N_H \cong
(\bz/4\bz)^4\times (\bz/2\bz)^2$.
{\it Degenerations of codimension $1$}:
$\aaa_1$ with $A_S\cong (\bz/4\bz)^4\times (\bz/2\bz)^3$
($H_{4,1}$ with the orbit of $\alpha_{2,1}$,  $\alpha_{1,2}$, $\alpha_{2,2}$
or $\alpha_{3,2}$);
$2\aaa_1$ with $A_S\cong (\bz/4\bz)^5$
($H_{4,1}$ with the orbit of $\alpha_{1,1}$ or $\alpha_{2,3}$);
$4\aaa_1$ with $A_S\cong \bz/8\bz \times (\bz/4\bz)^2\times  (\bz/2\bz)^2$
($H_{4,1}$ with the orbit of $\alpha_{1,3}$, $\alpha_{1,4}$,
$\alpha_{2,4}$ or $\alpha_{3,4}$).

\medskip

{\bf n=3,} $H\cong C_2^2$ ($|H|=4$, $i=2$):
$\rk N_H=12$ and $(N_H)^\ast/N_H\cong (\bz/4\bz)^2\times (\bz/2\bz)^6$.
{\it Degenerations of codimension $1$}:
$\aaa_1$ with $A_S\cong (\bz/4\bz)^2\times (\bz/2\bz)^7$
($H_{3,1}$ with the orbit of $\alpha_{2,1}$, $\alpha_{2,2}$, $\alpha_{2,4}$
or  $\alpha_{2,5}$;
$H_{3,2}$ with the orbit of $\alpha_{1,1}$, $\alpha_{2,1}$, $\alpha_{3,1}$,
$\alpha_{2,2}$, $\alpha_{2,4}$
or  $\alpha_{2,5}$;
$H_{3,3}$ with the orbit of $\alpha_{1,1}$, $\alpha_{2,1}$, $\alpha_{3,1}$,
$\alpha_{1,5}$, $\alpha_{2,5}$ or  $\alpha_{3,5}$;
$H_{3,4}$ with the orbit of $\alpha_{2,1}$ or  $\alpha_{2,5}$);
$2\aaa_1$ with $A_S\cong (\bz/4\bz)^3\times (\bz/2\bz)^4$
($H_{3,1}$ with the orbit of $\alpha_{1,1}$, $\alpha_{1,2}$,
$\alpha_{2,3}$, $\alpha_{1,4}$, $\alpha_{1,5}$ or $\alpha_{2,7}$;
$H_{3,2}$ with the orbit of $\alpha_{1,2}$, $\alpha_{1,4}$
or $\alpha_{1,5}$;
$H_{3,3}$ with the orbit of $\alpha_{2,2}$, $\alpha_{2,3}$,
or $\alpha_{2,7}$;
$H_{3,4}$ with the orbit of $\alpha_{1,1}$, $\alpha_{2,2}$,
$\alpha_{1,3}$, $\alpha_{2,3}$, $\alpha_{3,3}$,
$\alpha_{1,5}$, $\alpha_{1,7}$, $\alpha_{2,7}$ or $\alpha_{3,7}$);
$4\aaa_1$ with $A_S\cong \bz/8\bz \times (\bz/2\bz)^6$
($H_{3,1}$ with the orbit of $\alpha_{1,3}$ or $\alpha_{1,7}$;
$H_{3,2}$ with the orbit of $\alpha_{1,3}$, $\alpha_{2,3}$
or $\alpha_{3,3}$;
$H_{3,3}$ with the orbit of $\alpha_{1,2}$, $\alpha_{1,3}$
or $\alpha_{1,7}$;
$H_{3,4}$ with the orbit of $\alpha_{1,2}$).

\medskip

{\bf n=2,} $H\cong C_3$ ($|H|=3$, $i=1$):
$\rk N_H=12$ and $(N_H)^\ast/N_H\cong (\bz/3\bz)^6$.
{\it Degenerations of codimension $1$}:
$\aaa_1$ with $A_S\cong \bz/6\bz\times (\bz/3\bz)^5$
($H_{2,1}$ with the orbit of $\alpha_{1,1}$,  $\alpha_{2,1}$, $\alpha_{3,1}$,
$\alpha_{1,2}$, $\alpha_{2,2}$ or $\alpha_{3,2}$);
$3\aaa_1$ with $A_S\cong \bz/6\bz \times (\bz/3\bz)^4$
($H_{2,1}$ with the orbit of $\alpha_{1,3}$, $\alpha_{2,3}$,
$\alpha_{3,3}$, $\alpha_{1,5}$, $\alpha_{2,5}$ or $\alpha_{3,5}$).

\medskip

{\bf n=1,} $H\cong C_2$ ($|H|=2$, $i=1$):
$\rk N_H=8$ and $(N_H)^\ast/N_H\cong (\bz/2\bz)^8$.
{\it Degenerations of codimension $1$}:
$\aaa_1$ with $A_S\cong (\bz/2\bz)^9$
($H_{1,1}$ with the orbit of $\alpha_{2,1}$,  $\alpha_{2,2}$, $\alpha_{2,3}$,
$\alpha_{2,4}$, $\alpha_{2,5}$,  $\alpha_{2,6}$, $\alpha_{2,7}$ or $\alpha_{2,8}$;
$H_{1,2}$ with the orbit of $\alpha_{1,1}$,  $\alpha_{2,1}$, $\alpha_{3,1}$,
$\alpha_{2,2}$, $\alpha_{2,4}$,  $\alpha_{1,5}$, $\alpha_{2,5}$ or $\alpha_{3,5}$);
$2\aaa_1$ with $A_S\cong \bz/4\bz\times (\bz/2\bz)^6$
($H_{1,1}$ with the orbit of $\alpha_{1,1}$,  $\alpha_{1,2}$, $\alpha_{1,3}$,
$\alpha_{1,4}$, $\alpha_{1,5}$,  $\alpha_{1,6}$, $\alpha_{1,7}$ or $\alpha_{1,8}$;
$H_{1,2}$ with the orbit of $\alpha_{1,2}$,  $\alpha_{1,3}$, $\alpha_{2,3}$,
$\alpha_{3,3}$, $\alpha_{1,4}$,  $\alpha_{1,7}$, $\alpha_{2,7}$ or $\alpha_{3,7}$).


\vskip1cm

{\bf Table 20. All possible markings by $N_{20}=N(6A_4)$ of degenerations of
codimension one of non-trivial finite symplectic automorphism
groups of K\"ahlerian K3 surfaces:}

\medskip

{\bf n=32,} $H\cong Hol(C_5)$ ($|H|=20$, $i=3$):
$\rk N_H=18$ and
$(N_H)^\ast/N_H \cong (\bz/10\bz)^2 \times \bz/5\bz$.
{\it Degenerations:}
$2\aaa_1$ with $A_S\cong \bz/20\bz\times (\bz/5\bz)^2$
($H_{32,2}$ with the orbit of $\alpha_{1,1}$);
$5\aaa_1$ with $A_S\cong (\bz/10\bz)^2\times \bz/2\bz$
($H_{32,2}$ with the orbit of $\alpha_{1,2}$, $\alpha_{2,2}$,
$\alpha_{3,2}$ or $\alpha_{4,2}$);
$10\aaa_1$ with $A_S\cong \bz/20\bz\times \bz/5\bz$
($H_{32,1}$ with the orbit of $\alpha_{1,2}$);
$5\aaa_2$ with $A_S\cong \bz/10\bz\times \bz/5\bz$
($H_{32,1}$ with the orbit of $\alpha_{2,2}$).

\medskip

{\bf n=16,} $H\cong D_{10}$ ($|H|=10$, $i=1$):
$\rk N_H=16$ and $(N_H)^\ast/N_H\cong (\bz/5\bz)^4$.
{\it Degenerations of codimension $1$:}
$\aaa_1$ with $A_S\cong \bz/10\bz\times (\bz/5\bz)^3$
($H_{16,1}$ with the orbit of $\alpha_{1,1}$, $\alpha_{2,1}$,
$\alpha_{3,1}$ or $\alpha_{4,1}$);
$5\aaa_1$ with $A_S\cong \bz/10\bz\times (\bz/5\bz)^2$
($H_{16,1}$ with the orbit of $\alpha_{1,2}$, $\alpha_{2,2}$,
$\alpha_{3,2}$
or $\alpha_{4,2}$).

\medskip

{\bf n=4,} $H\cong C_4$ ($|H|=4$, $i=1$):
$\rk N_H=14$ and $(N_H)^\ast/N_H \cong
(\bz/4\bz)^4\times (\bz/2\bz)^2$.
{\it Degenerations of codimension $1$}:
$\aaa_1$ with $A_S\cong (\bz/4\bz)^4\times (\bz/2\bz)^3$
($H_{4,1}$ with the orbit of $\alpha_{1,1}$,
$\alpha_{2,1}$, $\alpha_{3,1}$
or $\alpha_{4,1}$);
$2\aaa_1$ with $A_S\cong (\bz/4\bz)^5$
($H_{4,1}$ with the orbit of $\alpha_{1,2}$);
$\aaa_2$ with $A_S\cong (\bz/4\bz)^4\times  \bz/2\bz$
($H_{4,1}$ with the orbit of $\alpha_{2,2}$);
$4\aaa_1$ with
$A_S\cong \bz/8\bz \times (\bz/4\bz)^2\times  (\bz/2\bz)^2$
($H_{4,1}$ with the orbit of $\alpha_{1,3}$, $\alpha_{2,3}$,
$\alpha_{3,3}$ or $\alpha_{4,3}$).

\medskip

{\bf n=1,} $H\cong C_2$ ($|H|=2$, $i=1$):
$\rk N_H=8$ and $(N_H)^\ast/N_H\cong (\bz/2\bz)^8$.
{\it Degenerations of codimension $1$}:
$\aaa_1$ with $A_S\cong (\bz/2\bz)^9$
($H_{1,1}$ with the orbit of $\alpha_{1,1}$,
$\alpha_{2,1}$, $\alpha_{3,1}$,
$\alpha_{4,1}$, $\alpha_{1,2}$,  $\alpha_{2,2}$, $\alpha_{3,2}$ or
$\alpha_{4,2}$);
$2\aaa_1$ with $A_S\cong \bz/4\bz\times (\bz/2\bz)^6$
($H_{1,1}$ with the orbit of $\alpha_{1,3}$,  $\alpha_{2,3}$,
$\alpha_{3,3}$,
$\alpha_{4,3}$, $\alpha_{1,4}$,  $\alpha_{2,4}$, $\alpha_{3,4}$
or $\alpha_{4,4}$).


\vskip1cm

{\bf Table 19. All possible markings by $N_{19}=N(6D_4)$ of degenerations of
codimension one of non-trivial  finite symplectic automorphism groups
of K\"ahlerian K3 surfaces:}

\medskip

{\bf n=70,} $H\cong {\mathfrak S}_5$ ($|H|=120$, $i=34$):
$\rk N_H=19$. {\it No degenerations.}

\medskip

{\bf n=61,} $H\cong {\mathfrak A}_{4,3}$ ($|H|=72$, $i=43$):
$\rk N_H=18$ and
$(N_H)^\ast/N_H\cong (\bz/12\bz)^2 \times \bz/3\bz$.
{\it Degenerations:}
$3\aaa_1$ with
$A_S\cong (\bz/12\bz)^2 \times \bz/2\bz$
($H_{61,1}$ with the orbit of $\alpha_{1,5}$ or $\alpha_{1,6}$);
$12\aaa_1$ with $A_S\cong \bz/24\bz\times \bz/3\bz$
($H_{61,1}$ with the orbit of $\alpha_{1,1}$).

\medskip

{\bf n=55,} $H\cong {\mathfrak A}_5$ ($|H|=60$, $i=5$):
$\rk N_H=18$ and $(N_H)^\ast/N_H\cong \bz/30\bz\times \bz/10\bz$.
{\it Degenerations:}
$\aaa_1$ with
$A_S\cong \bz/30\bz\times \bz/10\bz\times \bz/2\bz$
($H_{55,1}$ with the orbit of $\alpha_{1,4}$, $\alpha_{2,4}$
$\alpha_{3,4}$ or $\alpha_{4,4}$);
$5\aaa_1$ with
$A_S\cong \bz/30\bz\times (\bz/2\bz)^2$
($H_{55,1}$ with the orbit of $\alpha_{2,1}$);
$15\aaa_1$ with $A_S\cong \bz/10\bz \times (\bz/2\bz)^2$
($H_{55,1}$ with the orbit of $\alpha_{1,1}$).

\medskip

{\bf n=48,} $H\cong {\mathfrak S}_{3,3}$ ($|H|=36$, $i=10$):
$\rk N_H=18$, $(N_H)^\ast/N_H \cong
\bz/18\bz\times
\bz/6\bz\times (\bz/3\bz)^2$.
{\it Degenerations:}
$3\aaa_1$ with $A_S\cong \bz/18\bz \times (\bz/6\bz)^2$
($H_{48,1}$ with the orbit of $\alpha_{1,3}$ or $\alpha_{2,4}$);
$6\aaa_1$ with $A_S\cong \bz/36\bz \times (\bz/3\bz)^2$
($H_{48,1}$ with the orbit of $\alpha_{1,1}$);
$9\aaa_1$ with $A_S\cong (\bz/6\bz)^3$
($H_{48,1}$ with the orbit of $\alpha_{1,4}$).

\medskip

{\bf n=34}, $H\cong {\mathfrak S}_4$ ($|H|=24$, $i=12$):
$\rk N_H=17$ and $(N_H)^\ast/N_H\cong
(\bz/12\bz)^2\times \bz/4\bz$.
{\it Degenerations of codimension $1$:}
$\aaa_1$ with $A_S\cong (\bz/12\bz)^2\times \bz/4\bz\times \bz/2\bz$
($H_{34,1}$ with the orbit of $\alpha_{2,3}$ or $\alpha_{2,4}$;
$H_{34,2}$ with the orbit of $\alpha_{2,3}$, $\alpha_{1,4}$ or $\alpha_{2,4}$);
$2\aaa_1$ with $A_S\cong (\bz/12\bz)^2\times (\bz/2\bz)^2$
($H_{34,2}$ with the orbit of $\alpha_{3,4}$);
$3\aaa_1$ with $A_S\cong \bz/12\bz\times (\bz/4\bz)^2\times \bz/2\bz$
($H_{34,1}$ with the orbit of $\alpha_{1,3}$ or $\alpha_{1,4}$;
$H_{34,2}$ with the orbit of $\alpha_{1,3}$);
$4\aaa_1$ with $A_S\cong \bz/24\bz\times \bz/12\bz$
($H_{34,1}$ with the orbit of $\alpha_{2,1}$ or $\alpha_{4,1}$;
$H_{34,2}$ with the orbit of $\alpha_{2,1}$);
$8\aaa_1$ with $A_S\cong (\bz/12\bz)^2$
($H_{34,1}$ with the orbit of $\alpha_{1,1}$);
$12\aaa_1$ with $A_S\cong \bz/24\bz\times \bz/4\bz$
($H_{34,2}$ with the orbit of $\alpha_{1,1}$).

\medskip

{\bf n=32,} $H\cong Hol(C_5)$ ($|H|=20$, $i=3$):
$\rk N_H=18$ and
$(N_H)^\ast/N_H \cong (\bz/10\bz)^2 \times \bz/5\bz$.
{\it Degenerations:}
$2\aaa_1$ with $A_S\cong \bz/20\bz\times (\bz/5\bz)^2$
($H_{32,1}$ with the orbit of $\alpha_{3,3}$);
$5\aaa_1$ with $A_S\cong (\bz/10\bz)^2\times \bz/2\bz$
($H_{32,1}$ with the orbit of $\alpha_{1,1}$ or $\alpha_{2,1}$);
$10\aaa_1$ with $A_S\cong \bz/20\bz\times \bz/5\bz$
($H_{32,1}$ with the orbit of $\alpha_{3,1}$).

\medskip

{\bf n=30,} $H\cong {\mathfrak A}_{3,3}$ ($|H|=18$, $i=4$):
$\rk N_H=16$ and
$(N_H)^\ast/N_H\cong \bz/9\bz\times (\bz/3\bz)^4$.
{\it Degenerations of codimension $1$:}
$3\aaa_1$ with $A_S\cong \bz/18\bz\times (\bz/3\bz)^3$
($H_{30,1}$ with the orbit of $\alpha_{2,1}$, $\alpha_{1,4}$,
$\alpha_{1,5}$ or $\alpha_{1,6}$);
$9\aaa_1$ with $A_S\cong \bz/6\bz\times (\bz/3\bz)^3$
($H_{30,1}$ with the orbit of $\alpha_{1,1}$).

\medskip

{\bf n=18,} $H\cong D_{12}$ ($|H|=12$, $i=4$):
$\rk N_H=16$ and $(N_H)^\ast/N_H\cong (\bz/6\bz)^4$.
{\it Degenerations of codimension $1$:}
$\aaa_1$ with $A_S\cong (\bz/6\bz)^4\times \bz/2\bz$
($H_{18,1}$ with the orbit of $\alpha_{2,3}$ or $\alpha_{2,4}$;
$H_{18,2}$ with the orbit of $\alpha_{2,3}$ or $\alpha_{4,3}$);
$2\aaa_1$ with $A_S\cong \bz/12\bz\times (\bz/6\bz)^2\times \bz/3\bz$
($H_{18,1}$ with the orbit of $\alpha_{2,1}$ or $\alpha_{2,5}$;
$H_{18,2}$ with the orbit of $\alpha_{1,3}$ or $\alpha_{2,4}$);
$3\aaa_1$ with $A_S\cong (\bz/6\bz)^3\times (\bz/2\bz)^2$
($H_{18,1}$ with the orbit of $\alpha_{1,3}$ or  $\alpha_{1,4}$;
$H_{18,2}$ with the orbit of $\alpha_{2,1}$ or $\alpha_{3,1}$);
$6\aaa_1$ with $A_S\cong \bz/12\bz\times (\bz/6\bz)^2$
($H_{18,1}$ with the orbit of $\alpha_{1,1}$ or $\alpha_{1,5}$;
$H_{18,2}$ with the orbit of $\alpha_{1,1}$ or $\alpha_{1,4}$).

\medskip

{\bf n=17,} $H\cong {\mathfrak A}_4$ ($|H|=12$, $i=3$):
$\rk N_H=16$ and $(N_H)^\ast/N_H\cong (\bz/12\bz)^2\times (\bz/2\bz)^2$.
{\it Degenerations of codimension $1$:}
$\aaa_1$ with $A_S\cong (\bz/12\bz)^2\times (\bz/2\bz)^3$
($H_{17,1}$ with the orbit of $\alpha_{2,3}$  or $\alpha_{2,4}$;
$H_{17,2}$ with the orbit of $\alpha_{2,3}$, $\alpha_{1,4}$,
$\alpha_{2,4}$, $\alpha_{3,4}$ or $\alpha_{4,4}$);
$3\aaa_1$ with $A_S\cong \bz/12\bz\times \bz/4\bz\times (\bz/2\bz)^3$
($H_{17,1}$ with the orbit of $\alpha_{1,3}$ or $\alpha_{1,4}$;
$H_{17,2}$ with the orbit of $\alpha_{1,3}$);
$4\aaa_1$ with $A_S\cong \bz/24\bz\times \bz/6\bz\times \bz/2\bz$
($H_{17,1}$ with the orbit of $\alpha_{1,1}$, $\alpha_{2,1}$, $\alpha_{3,1}$
or $\alpha_{4,1}$;
$H_{17,2}$ with the orbit of $\alpha_{2,1}$);
$12\aaa_1$ with $A_S\cong \bz/24\bz\times (\bz/2\bz)^2$
($H_{17,2}$ with the orbit of $\alpha_{1,1}$).

\medskip

{\bf n=16,} $H\cong D_{10}$ ($|H|=10$, $i=1$):
$\rk N_H=16$ and $(N_H)^\ast/N_H\cong (\bz/5\bz)^4$.
{\it Degenerations of codimension $1$:}
$\aaa_1$ with $A_S\cong \bz/10\bz\times (\bz/5\bz)^3$
($H_{16,1}$ with the orbit of $\alpha_{1,3}$, $\alpha_{2,3}$,
$\alpha_{3,3}$ or $\alpha_{4,3}$);
$5\aaa_1$ with $A_S\cong \bz/10\bz\times (\bz/5\bz)^2$
($H_{16,1}$ with the orbit of $\alpha_{1,1}$, $\alpha_{2,1}$,
$\alpha_{3,1}$ or $\alpha_{4,1}$).

\medskip

{\bf n=10,} $H\cong D_8$ ($|H|=8$, $i=3$):
$\rk N_H=15$ and $(N_H)^\ast/N_H\cong
(\bz/4\bz)^5$.
{\it Degenerations of codimension $1$:}
$\aaa_1$ with $A_S\cong (\bz/4\bz)^5\times \bz/2\bz$
($H_{10,1}$ with the orbit of $\alpha_{1,3}$, $\alpha_{2,3}$,
$\alpha_{2,4}$ or $\alpha_{3,4}$);
$2\aaa_1$ with $A_S\cong (\bz/4\bz)^4\times (\bz/2\bz)^2$
($H_{10,1}$ with the orbit of $\alpha_{3,3}$ or $\alpha_{1,4}$);
$4\aaa_1$ with $A_S\cong \bz/8\bz\times (\bz/4\bz)^3$
($H_{10,1}$ with the orbit of $\alpha_{2,1}$ or $\alpha_{4,1}$);
$8\aaa_1$ with $A_S\cong (\bz/4\bz)^4$
($H_{10,1}$ with the orbit of $\alpha_{1,1}$).

\medskip

{\bf n=6,} $H\cong D_6$ ($|H|=6$, $i=1$):
$\rk N_H=14$ and $(N_H)^\ast/N_H\cong (\bz/6\bz)^2\times (\bz/3\bz)^3$.
{\it Degenerations of codimension $1$}:
$\aaa_1$ with $A_S\cong (\bz/6\bz)^3\times (\bz/3\bz)^2$
($H_{6,1}$ with the orbit of $\alpha_{2,3}$,  $\alpha_{2,4}$, $\alpha_{2,5}$
or $\alpha_{2,6}$;
$H_{6,2}$ with the orbit of $\alpha_{2,4}$,  $\alpha_{4,4}$, $\alpha_{2,5}$
or $\alpha_{2,6}$;
$H_{6,3}$ with the orbit of $\alpha_{1,3}$,  $\alpha_{2,3}$, $\alpha_{3,3}$
or $\alpha_{4,3}$);
$2\aaa_1$ with with $A_S\cong \bz/12\bz\times (\bz/3\bz)^4$
($H_{6,1}$ with the orbit of $\alpha_{2,1}$;
$H_{6,2}$ with the orbit of $\alpha_{1,4}$;
$H_{6,3}$ with the orbit of $\alpha_{2,1}$);
$3\aaa_1$ with $A_S\cong (\bz/6\bz)^3\times \bz/3\bz$
($H_{6,1}$ with the orbit of $\alpha_{1,3}$, $\alpha_{1,4}$,
$\alpha_{1,5}$ or $\alpha_{1,6}$;
$H_{6,2}$ with the orbit of $\alpha_{2,1}$, $\alpha_{4,1}$,
$\alpha_{1,5}$ or $\alpha_{1,6}$;
$H_{6,3}$ with the orbit of $\alpha_{1,4}$, $\alpha_{2,4}$,
$\alpha_{3,4}$ or $\alpha_{4,4}$);
$6\aaa_1$ with $A_S\cong \bz/12\bz\times (\bz/3\bz)^3$
($H_{6,1}$ with the orbit of $\alpha_{1,1}$;
$H_{6,2}$ with the orbit of $\alpha_{1,1}$;
$H_{6,3}$ with the orbit of $\alpha_{1,1}$).

\medskip

{\bf n=4,} $H\cong C_4$ ($|H|=4$, $i=1$):
$\rk N_H=14$ and $(N_H)^\ast/N_H \cong
(\bz/4\bz)^4\times (\bz/2\bz)^2$.
{\it Degenerations of codimension $1$}:
$\aaa_1$ with $A_S\cong (\bz/4\bz)^4\times (\bz/2\bz)^3$
($H_{4,1}$ with the orbit of $\alpha_{2,3}$,  $\alpha_{3,3}$
$\alpha_{1,4}$ or $\alpha_{2,4}$);
$2\aaa_1$ with with $A_S\cong (\bz/4\bz)^5$
($H_{4,1}$ with the orbit of $\alpha_{1,3}$ or $\alpha_{3,4}$);
$4\aaa_1$ with $A_S\cong \bz/8\bz \times (\bz/4\bz)^2\times  (\bz/2\bz)^2$
($H_{4,1}$ with the orbit of $\alpha_{1,1}$, $\alpha_{2,1}$,
$\alpha_{3,1}$ or $\alpha_{4,1}$).

\medskip

{\bf n=3,} $H\cong C_2^2$ ($|H|=4$, $i=2$):
$\rk N_H=12$ and $(N_H)^\ast/N_H\cong (\bz/4\bz)^2\times (\bz/2\bz)^6$.
{\it Degenerations of codimension $1$}:
$\aaa_1$ with $A_S\cong (\bz/4\bz)^2\times (\bz/2\bz)^7$
($H_{3,1}$ with the orbit of $\alpha_{2,5}$, $\alpha_{4,5}$, $\alpha_{2,6}$
or  $\alpha_{4,6}$;
$H_{3,2}$ with the orbit of $\alpha_{1,3}$, $\alpha_{2,3}$, $\alpha_{3,3}$,
$\alpha_{4,3}$, $\alpha_{1,4}$, $\alpha_{3,4}$
or  $\alpha_{4,4}$);
$2\aaa_1$ with $A_S\cong (\bz/4\bz)^3\times (\bz/2\bz)^4$
($H_{3,1}$ with the orbit of $\alpha_{2,1}$,
$\alpha_{4,1}$, $\alpha_{2,3}$, $\alpha_{4,3}$, $\alpha_{1,5}$ or $\alpha_{1,6}$;
$H_{3,3}$ with the orbit of $\alpha_{1,1}$, $\alpha_{2,1}$,
$\alpha_{3,1}$, $\alpha_{4,1}$, $\alpha_{1,3}$,
$\alpha_{2,3}$, $\alpha_{3,3}$, $\alpha_{4,3}$,
$\alpha_{1,5}$,  $\alpha_{2,5}$, $\alpha_{3,5}$ or $\alpha_{4,5}$);
$4\aaa_1$ with $A_S\cong \bz/8\bz \times (\bz/2\bz)^6$
($H_{3,1}$ with the orbit of $\alpha_{1,1}$ or $\alpha_{1,3}$;
$H_{3,2}$ with the orbit of $\alpha_{1,1}$ $\alpha_{2,1}$, $\alpha_{3,1}$
or $\alpha_{4,1}$).

\medskip

{\bf n=2,} $H\cong C_3$ ($|H|=3$, $i=1$):
$\rk N_H=12$ and $(N_H)^\ast/N_H\cong (\bz/3\bz)^6$.
{\it Degenerations of codimension $1$}:
$\aaa_1$ with $A_S\cong \bz/6\bz\times (\bz/3\bz)^5$
($H_{2,1}$ with the orbit of $\alpha_{2,1}$,  $\alpha_{2,2}$, $\alpha_{2,3}$,
$\alpha_{2,4}$, $\alpha_{2,5}$ or $\alpha_{2,6}$;
$H_{2,2}$ with the orbit of $\alpha_{2,3}$,  $\alpha_{1,4}$, $\alpha_{2,4}$,
$\alpha_{3,4}$, $\alpha_{4,4}$ or $\alpha_{2,5}$);
$3\aaa_1$ with with $A_S\cong \bz/6\bz \times (\bz/3\bz)^4$
($H_{2,1}$ with the orbit of $\alpha_{1,1}$, $\alpha_{1,2}$,
$\alpha_{1,3}$, $\alpha_{1,4}$, $\alpha_{1,5}$ or $\alpha_{1,6}$;
$H_{2,2}$ with the orbit of $\alpha_{1,1}$, $\alpha_{2,1}$,
$\alpha_{3,1}$, $\alpha_{4,1}$, $\alpha_{1,3}$ or $\alpha_{1,5}$).

\medskip

{\bf n=1,} $H\cong C_2$ ($|H|=2$, $i=1$):
$\rk N_H=8$ and $(N_H)^\ast/N_H\cong (\bz/2\bz)^8$.
{\it Degenerations of codimension $1$}:
$\aaa_1$ with $A_S\cong (\bz/2\bz)^9$
($H_{1,1}$ with the orbit of $\alpha_{2,3}$,  $\alpha_{4,3}$, $\alpha_{2,4}$,
$\alpha_{4,4}$, $\alpha_{2,5}$,  $\alpha_{4,5}$, $\alpha_{2,6}$ or $\alpha_{4,6}$;
$H_{1,2}$ with the orbit of $\alpha_{1,3}$,  $\alpha_{2,3}$, $\alpha_{3,3}$,
$\alpha_{4,3}$, $\alpha_{1,4}$,  $\alpha_{2,4}$, $\alpha_{3,4}$ or $\alpha_{4,4}$);
$2\aaa_1$ with $A_S\cong \bz/4\bz\times (\bz/2\bz)^6$
($H_{1,1}$ with the orbit of $\alpha_{1,1}$, $\alpha_{2,1}$,
$\alpha_{3,1}$, $\alpha_{4,1}$, $\alpha_{1,3}$, $\alpha_{1,4}$,
$\alpha_{1,5}$  or $\alpha_{1,6}$;
$H_{1,2}$ with the orbit of $\alpha_{1,1}$, $\alpha_{2,1}$,
$\alpha_{3,1}$, $\alpha_{4,1}$, $\alpha_{1,5}$, $\alpha_{2,5}$,
$\alpha_{3,5}$  or $\alpha_{4,5}$).


\vskip1cm

{\bf Table 18. All possible markings by $N_{18}=N(4A_5\oplus D_4)$ of degenerations of
codimension one of non-trivial  finite symplectic automorphism
groups of K\"ahlerian K3 surfaces:}

\medskip

{\bf n=54,} $H\cong T_{48}$ ($|H|=48$, $i=29$):
$\rk N_H=19$.
{\it No degenerations.}

\medskip

{\bf n=26,} $H\cong SD_{16}$ ($|H|=16$, $i=8$):
$\rk N_H=18$ and
$(N_H)^\ast/N_H\cong (\bz/8\bz)^2\times
\bz/4\bz\times \bz/2\bz$.
{\it Degenerations:}
$8\aaa_1$ with $A_S\cong \bz/16\bz\times \bz/4\bz\times \bz/2\bz$
($H_{26,1}$ with the orbit of $\alpha_{1,1}$ or $\alpha_{2,1}$).

\medskip

{\bf n=18,} $H\cong D_{12}$ ($|H|=12$, $i=4$):
$\rk N_H=16$ and $(N_H)^\ast/N_H\cong (\bz/6\bz)^4$.
{\it Degenerations of codimension $1$:}
$\aaa_1$ with $A_S\cong (\bz/6\bz)^4\times \bz/2\bz$
($H_{18,1}$ with the orbit of $\alpha_{3,1}$ or $\alpha_{2,5}$);
$2\aaa_1$ with $A_S\cong \bz/12\bz\times (\bz/6\bz)^2\times \bz/3\bz$
($H_{18,1}$ with the orbit of $\alpha_{1,1}$ or $\alpha_{2,1}$);
$3\aaa_1$ with $A_S\cong (\bz/6\bz)^3\times (\bz/2\bz)^2$
($H_{18,1}$ with the orbit of $\alpha_{3,2}$ or  $\alpha_{1,5}$);
$6\aaa_1$ with $A_S\cong \bz/12\bz\times (\bz/6\bz)^2$
($H_{18,1}$ with the orbit of $\alpha_{1,2}$ or $\alpha_{2,2}$).

\medskip

{\bf n=12,} $H\cong Q_8$ ($|H|=8$, $i=4$):
$\rk N_H=17$, $(N_H)^\ast/N_H \cong
(\bz/8\bz)^2\times (\bz/2\bz)^3$.
{\it Degenerations of codimension $1$:}
$8\aaa_1$ with $A_S\cong \bz/16\bz\times (\bz/2\bz)^3$
($H_{12,1}$ with the orbit of $\alpha_{1,1}$ or $\alpha_{2,1}$).

\medskip

{\bf n=10,} $H\cong D_8$ ($|H|=8$, $i=3$):
$\rk N_H=15$ and $(N_H)^\ast/N_H\cong
(\bz/4\bz)^5$.
{\it Degenerations of codimension $1$:}
$\aaa_1$ with $A_S\cong (\bz/4\bz)^5\times \bz/2\bz$
($H_{10,1}$ with the orbit of $\alpha_{2,5}$ or $\alpha_{4,5}$);
$2\aaa_1$ with $A_S\cong (\bz/4\bz)^4\times (\bz/2\bz)^2$
($H_{10,1}$ with the orbit of $\alpha_{3,1}$, $\alpha_{3,3}$ or
$\alpha_{1,5}$);
$4\aaa_1$ with $A_S\cong \bz/8\bz\times (\bz/4\bz)^3$
($H_{10,1}$ with the orbit of $\alpha_{1,1}$, $\alpha_{2,1}$,
$\alpha_{1,3}$ or $\alpha_{2,3}$).

\medskip

{\bf n=6,} $H\cong D_6$ ($|H|=6$, $i=1$):
$\rk N_H=14$ and $(N_H)^\ast/N_H\cong (\bz/6\bz)^2\times (\bz/3\bz)^3$.
{\it Degenerations of codimension $1$}:
$\aaa_1$ with $A_S\cong (\bz/6\bz)^3\times (\bz/3\bz)^2$
($H_{6,1}$ with the orbit of $\alpha_{1,1}$,  $\alpha_{2,1}$,
$\alpha_{3,1}$, $\alpha_{4,1}$, $\alpha_{5,1}$
or $\alpha_{2,5}$;
$H_{6,2}$ with the orbit of $\alpha_{3,1}$ or $\alpha_{2,5}$);
$2\aaa_1$ with $A_S\cong \bz/12\bz\times (\bz/3\bz)^4$
($H_{6,2}$ with the orbit of $\alpha_{1,1}$ or $\alpha_{2,1}$);
$3\aaa_1$ with $A_S\cong (\bz/6\bz)^3\times \bz/3\bz$
($H_{6,1}$ with the orbit of $\alpha_{3,2}$ or $\alpha_{1,5}$;
$H_{6,2}$ with the orbit of $\alpha_{1,2}$, $\alpha_{2,2}$,
$\alpha_{3,2}$, $\alpha_{4,2}$, $\alpha_{5,2}$ or $\alpha_{1,5}$);
$6\aaa_1$ with $A_S\cong \bz/12\bz\times (\bz/3\bz)^3$
($H_{6,1}$ with the orbit of $\alpha_{1,2}$ or $\alpha_{2,2}$).

\medskip

{\bf n=4,} $H\cong C_4$ ($|H|=4$, $i=1$):
$\rk N_H=14$ and $(N_H)^\ast/N_H \cong
(\bz/4\bz)^4\times (\bz/2\bz)^2$.
{\it Degenerations of codimension $1$}:
$\aaa_1$ with $A_S\cong (\bz/4\bz)^4\times (\bz/2\bz)^3$
($H_{4,1}$ with the orbit of $\alpha_{1,5}$, $\alpha_{2,5}$,
$\alpha_{3,5}$ or $\alpha_{4,5}$);
$2\aaa_1$ with $A_S\cong (\bz/4\bz)^5$
($H_{4,1}$ with the orbit of $\alpha_{3,1}$ or $\alpha_{3,3}$);
$4\aaa_1$ with $A_S\cong \bz/8\bz \times (\bz/4\bz)^2\times (\bz/2\bz)^2$
($H_{4,1}$ with the orbit of $\alpha_{1,1}$, $\alpha_{2,1}$,
$\alpha_{1,3}$ or $\alpha_{2,3}$).

\medskip

{\bf n=3,} $H\cong C_2^2$ ($|H|=4$, $i=2$):
$\rk N_H=12$ and $(N_H)^\ast/N_H\cong (\bz/4\bz)^2\times (\bz/2\bz)^6$.
{\it Degenerations of codimension $1$}:
$\aaa_1$ with $A_S\cong (\bz/4\bz)^2\times (\bz/2\bz)^7$
($H_{3,1}$ with the orbit of $\alpha_{3,1}$, $\alpha_{3,2}$,
$\alpha_{2,5}$ or  $\alpha_{4,5}$);
$2\aaa_1$ with $A_S\cong (\bz/4\bz)^3\times (\bz/2\bz)^4$
($H_{3,1}$ with the orbit of $\alpha_{1,1}$,
$\alpha_{2,1}$, $\alpha_{1,2}$, $\alpha_{2,2}$, $\alpha_{3,3}$
or $\alpha_{1,5}$);
$4\aaa_1$ with $A_S\cong \bz/8\bz \times (\bz/2\bz)^6$
($H_{3,1}$ with the orbit of $\alpha_{1,3}$ or $\alpha_{1,4}$).

\medskip

{\bf n=2,} $H\cong C_3$ ($|H|=3$, $i=1$):
$\rk N_H=12$ and $(N_H)^\ast/N_H\cong (\bz/3\bz)^6$.
{\it Degenerations of codimension $1$}:
$\aaa_1$ with $A_S\cong \bz/6\bz\times (\bz/3\bz)^5$
($H_{2,1}$ with the orbit of $\alpha_{1,1}$,  $\alpha_{2,1}$,
$\alpha_{3,1}$, $\alpha_{4,1}$, $\alpha_{5,1}$ or $\alpha_{2,5}$);
$3\aaa_1$ with $A_S\cong \bz/6\bz \times (\bz/3\bz)^4$
($H_{2,1}$ with the orbit of $\alpha_{1,2}$, $\alpha_{2,2}$,
$\alpha_{3,2}$, $\alpha_{4,2}$, $\alpha_{5,2}$ or $\alpha_{1,5}$).

\medskip

{\bf n=1,} $H\cong C_2$ ($|H|=2$, $i=1$):
$\rk N_H=8$ and $(N_H)^\ast/N_H\cong (\bz/2\bz)^8$.
{\it Degenerations of codimension $1$}:
$\aaa_1$ with $A_S\cong (\bz/2\bz)^9$
($H_{1,1}$ with the orbit of $\alpha_{1,1}$,
$\alpha_{2,1}$, $\alpha_{3,1}$,
$\alpha_{4,1}$, $\alpha_{5,1}$, $\alpha_{3,2}$, $\alpha_{2,5}$
or $\alpha_{4,5}$;
$H_{1,2}$ with the orbit of $\alpha_{3,1}$, $\alpha_{3,2}$,
$\alpha_{3,3}$,
$\alpha_{3,4}$, $\alpha_{1,5}$, $\alpha_{2,5}$, $\alpha_{3,5}$ or
$\alpha_{4,5}$);
$2\aaa_1$ with $A_S\cong \bz/4\bz\times (\bz/2\bz)^6$
($H_{1,1}$ with the orbit of $\alpha_{1,2}$, $\alpha_{2,2}$,
$\alpha_{1,3}$, $\alpha_{2,3}$, $\alpha_{3,3}$, $\alpha_{4,3}$,
$\alpha_{5,3}$  or $\alpha_{1,5}$;
$H_{1,2}$ with the orbit of $\alpha_{1,1}$, $\alpha_{2,1}$,
$\alpha_{1,2}$, $\alpha_{2,2}$, $\alpha_{1,3}$, $\alpha_{2,3}$,
$\alpha_{1,4}$  or $\alpha_{2,4}$).


\vskip1cm

{\bf Table 17. All possible markings by $N_{17}=N(4A_6)$ of degenerations of
codimension one of non-trivial  finite symplectic automorphism groups
of K\"ahlerian K3 surfaces:}

\medskip

{\bf n=2,} $H\cong C_3$ ($|H|=3$, $i=1$):
$\rk N_H=12$ and $(N_H)^\ast/N_H\cong (\bz/3\bz)^6$.
{\it Degenerations of codimension $1$}:
$\aaa_1$ with $A_S\cong \bz/6\bz\times (\bz/3\bz)^5$
($H_{2,1}$ with the orbit of $\alpha_{1,1}$,  $\alpha_{2,1}$,
$\alpha_{3,1}$, $\alpha_{4,1}$, $\alpha_{5,1}$ or $\alpha_{6,1}$);
$3\aaa_1$ with $A_S\cong \bz/6\bz \times (\bz/3\bz)^4$
($H_{2,1}$ with the orbit of $\alpha_{1,2}$, $\alpha_{2,2}$,
$\alpha_{3,2}$, $\alpha_{4,2}$, $\alpha_{5,2}$ or $\alpha_{6,2}$).


\vskip1cm

{\bf Table 16. All possible markings by $N_{16}=N(2A_7\oplus 2D_5)$ of degenerations of
codimension one of  non-trivial finite symplectic automorphism groups
of K\"ahlerian K3 surfaces:}

\medskip

{\bf n=3,} $H\cong C_2^2$ ($|H|=4$, $i=2$):
$\rk N_H=12$ and $(N_H)^\ast/N_H\cong (\bz/4\bz)^2\times (\bz/2\bz)^6$.
{\it Degenerations of codimension $1$}:
$\aaa_1$ with $A_S\cong (\bz/4\bz)^2\times (\bz/2\bz)^7$
($H_{3,1}$ with the orbit of $\alpha_{4,1}$ or $\alpha_{4,2}$;
$H_{3,2}$ with the orbit of $\alpha_{1,3}$, $\alpha_{2,3}$,
$\alpha_{3,3}$, $\alpha_{1,4}$, $\alpha_{2,4}$ or $\alpha_{3,4}$);
$2\aaa_1$ with $A_S\cong (\bz/4\bz)^3\times (\bz/2\bz)^4$
($H_{3,1}$ with the orbit of $\alpha_{1,1}$, $\alpha_{2,1}$,
$\alpha_{3,1}$, $\alpha_{1,2}$, $\alpha_{2,2}$, $\alpha_{3,2}$,
$\alpha_{1,3}$, $\alpha_{2,3}$ or $\alpha_{3,3}$;
$H_{3,2}$ with the orbit of $\alpha_{4,1}$, $\alpha_{4,3}$ or
$\alpha_{4,4}$);
$4\aaa_1$ with $A_S\cong \bz/8\bz \times (\bz/2\bz)^6$
($H_{3,1}$ with the orbit of $\alpha_{4,3}$;
$H_{3,2}$ with the orbit of $\alpha_{1,1}$, $\alpha_{2,1}$ or
$\alpha_{3,1}$).

\medskip

{\bf n=1,} $H\cong C_2$ ($|H|=2$, $i=1$):
$\rk N_H=8$ and $(N_H)^\ast/N_H\cong (\bz/2\bz)^8$.
{\it Degenerations of codimension $1$}:
$\aaa_1$ with $A_S\cong (\bz/2\bz)^9$
($H_{1,1}$ with the orbit of $\alpha_{1,1}$, $\alpha_{2,1}$,
$\alpha_{3,1}$,
$\alpha_{4,1}$, $\alpha_{5,1}$, $\alpha_{6,1}$, $\alpha_{7,1}$ or
$\alpha_{4,2}$;
$H_{1,2}$ with the orbit of $\alpha_{4,1}$, $\alpha_{4,2}$,
$\alpha_{1,3}$,
$\alpha_{2,3}$, $\alpha_{3,3}$, $\alpha_{1,4}$, $\alpha_{2,4}$ or
$\alpha_{3,4}$;
$H_{1,3}$ with the orbit of $\alpha_{1,3}$, $\alpha_{2,3}$,
$\alpha_{3,3}$,
$\alpha_{4,3}$, $\alpha_{5,3}$, $\alpha_{1,4}$, $\alpha_{2,4}$ or
$\alpha_{3,4}$);
$2\aaa_1$ with $A_S\cong \bz/4\bz\times (\bz/2\bz)^6$
($H_{1,1}$ with the orbit of $\alpha_{1,2}$, $\alpha_{2,2}$,
$\alpha_{3,2}$, $\alpha_{1,3}$, $\alpha_{2,3}$, $\alpha_{3,3}$,
$\alpha_{4,3}$  or $\alpha_{5,3}$;
$H_{1,2}$ with the orbit of $\alpha_{1,1}$, $\alpha_{2,1}$,
$\alpha_{3,1}$, $\alpha_{1,2}$, $\alpha_{2,2}$, $\alpha_{3,2}$,
$\alpha_{4,3}$  or $\alpha_{4,4}$;
$H_{1,3}$ with the orbit of $\alpha_{1,1}$, $\alpha_{2,1}$,
$\alpha_{3,1}$, $\alpha_{4,1}$, $\alpha_{5,1}$, $\alpha_{6,1}$,
$\alpha_{7,1}$  or $\alpha_{4,4}$).


\vskip1cm

{\bf Table 15. All possible markings by $N_{15}=N(3A_8)$ of degenerations of
codimension one of non-trivial finite symplectic automorphism groups
of K\"ahlerian K3 surfaces:}

\medskip

{\bf n=1,} $H\cong C_2$ ($|H|=2$, $i=1$):
$\rk N_H=8$ and $(N_H)^\ast/N_H\cong (\bz/2\bz)^8$.
{\it Degenerations of codimension $1$}:
$\aaa_1$ with $A_S\cong (\bz/2\bz)^9$
($H_{1,1}$ with the orbit of $\alpha_{1,1}$, $\alpha_{2,1}$,
$\alpha_{3,1}$,
$\alpha_{4,1}$, $\alpha_{5,1}$, $\alpha_{6,1}$, $\alpha_{7,1}$
or $\alpha_{8,1}$);
$2\aaa_1$ with $A_S\cong \bz/4\bz\times (\bz/2\bz)^6$
($H_{1,1}$ with the orbit of $\alpha_{1,2}$, $\alpha_{2,2}$,
$\alpha_{3,2}$, $\alpha_{4,2}$, $\alpha_{5,2}$, $\alpha_{6,2}$,
$\alpha_{7,2}$  or $\alpha_{8,2}$).


\vskip1cm

{\bf Table 14. All possible markings by $N_{14}=N(4D_6)$ of degenerations of
codimension one of  non-trivial finite symplectic automorphism groups
of K\"ahlerian K3 surfaces:}

\medskip

{\bf n=6,} $H\cong D_6$ ($|H|=6$, $i=1$):
$\rk N_H=14$ and $(N_H)^\ast/N_H\cong (\bz/6\bz)^2\times (\bz/3\bz)^3$.
{\it Degenerations of codimension $1$}:
$\aaa_1$ with $A_S\cong (\bz/6\bz)^3\times (\bz/3\bz)^2$
($H_{6,1}$ with the orbit of $\alpha_{1,1}$,  $\alpha_{2,1}$,
$\alpha_{3,1}$ or $\alpha_{4,1}$);
$2\aaa_1$ with $A_S\cong \bz/12\bz\times (\bz/3\bz)^4$
($H_{6,1}$ with the orbit of $\alpha_{5,1}$);
$3\aaa_1$ with $A_S\cong (\bz/6\bz)^3\times \bz/3\bz$
($H_{6,1}$ with the orbit of $\alpha_{1,2}$, $\alpha_{2,2}$,
$\alpha_{3,2}$ or $\alpha_{4,2}$);
$6\aaa_1$ with $A_S\cong \bz/12\bz\times (\bz/3\bz)^3$
($H_{6,1}$ with the orbit of $\alpha_{5,2}$).

\medskip

{\bf n=2,} $H\cong C_3$ ($|H|=3$, $i=1$):
$\rk N_H=12$ and $(N_H)^\ast/N_H\cong (\bz/3\bz)^6$.
{\it Degenerations of codimension $1$}:
$\aaa_1$ with $A_S\cong \bz/6\bz\times (\bz/3\bz)^5$
($H_{2,1}$ with the orbit of $\alpha_{1,1}$,  $\alpha_{2,1}$,
$\alpha_{3,1}$, $\alpha_{4,1}$, $\alpha_{5,1}$ or $\alpha_{6,1}$);
$3\aaa_1$ with $A_S\cong \bz/6\bz \times (\bz/3\bz)^4$
($H_{2,1}$ with the orbit of $\alpha_{1,2}$, $\alpha_{2,2}$,
$\alpha_{3,2}$, $\alpha_{4,2}$, $\alpha_{5,2}$ or $\alpha_{6,2}$).

\medskip

{\bf n=1,} $H\cong C_2$ ($|H|=2$, $i=1$):
$\rk N_H=8$ and $(N_H)^\ast/N_H\cong (\bz/2\bz)^8$.
{\it Degenerations of codimension $1$}:
$\aaa_1$ with $A_S\cong (\bz/2\bz)^9$
($H_{1,1}$ with the orbit of $\alpha_{1,1}$, $\alpha_{2,1}$,
$\alpha_{3,1}$,
$\alpha_{4,1}$, $\alpha_{1,2}$, $\alpha_{2,2}$, $\alpha_{3,2}$ or
$\alpha_{4,2}$);
$2\aaa_1$ with $A_S\cong \bz/4\bz\times (\bz/2\bz)^6$
($H_{1,1}$ with the orbit of $\alpha_{5,1}$, $\alpha_{5,2}$,
$\alpha_{1,3}$, $\alpha_{2,3}$, $\alpha_{3,3}$, $\alpha_{4,3}$,
$\alpha_{5,3}$  or $\alpha_{6,3}$).


\vskip1cm

{\bf Table 13. All possible markings by $N_{13}=N(2A_9\oplus D_6)$ of degenerations of
codimension one of  non-trivial finite symplectic automorphism groups
of K\"ahlerian K3 surfaces:}

\medskip

{\bf n=4,} $H\cong C_4$ ($|H|=4$, $i=1$):
$\rk N_H=14$ and $(N_H)^\ast/N_H \cong
(\bz/4\bz)^4\times (\bz/2\bz)^2$.
{\it Degenerations of codimension $1$}:
$\aaa_1$ with $A_S\cong (\bz/4\bz)^4\times (\bz/2\bz)^3$
($H_{4,1}$ with the orbit of $\alpha_{1,3}$, $\alpha_{2,3}$,
$\alpha_{3,3}$ or $\alpha_{4,3}$);
$2\aaa_1$ with $A_S\cong (\bz/4\bz)^5$
($H_{4,1}$ with the orbit of $\alpha_{2,1}$ or $\alpha_{5,3}$);
$4\aaa_1$ with $A_S\cong \bz/8\bz \times (\bz/4\bz)^2\times (\bz/2\bz)^2$
($H_{4,1}$ with the orbit of $\alpha_{1,1}$, $\alpha_{2,1}$,
$\alpha_{3,1}$ or $\alpha_{4,1}$).

\medskip

{\bf n=1,} $H\cong C_2$ ($|H|=2$, $i=1$):
$\rk N_H=8$ and $(N_H)^\ast/N_H\cong (\bz/2\bz)^8$.
{\it Degenerations of codimension $1$}:
$\aaa_1$ with $A_S\cong (\bz/2\bz)^9$
($H_{1,1}$ with the orbit of $\alpha_{5,1}$, $\alpha_{5,2}$,
$\alpha_{1,3}$,
$\alpha_{2,3}$, $\alpha_{3,3}$, $\alpha_{4,3}$, $\alpha_{5,3}$
or $\alpha_{6,3}$);
$2\aaa_1$ with $A_S\cong \bz/4\bz\times (\bz/2\bz)^6$
($H_{1,1}$ with the orbit of $\alpha_{1,1}$, $\alpha_{2,1}$,
$\alpha_{3,1}$, $\alpha_{4,1}$, $\alpha_{1,2}$, $\alpha_{2,2}$,
$\alpha_{3,2}$  or $\alpha_{4,2}$).


\vskip1cm

{\bf Table 12. All possible markings by $N_{12}=N(4E_6)$ of degenerations of
codimension one of non-trivial  finite symplectic automorphism groups
of K\"ahlerian K3 surfaces:}

\medskip

{\bf n=18,} $H\cong D_{12}$ ($|H|=12$, $i=4$):
$\rk N_H=16$ and $(N_H)^\ast/N_H\cong (\bz/6\bz)^4$.
{\it Degenerations of codimension $1$:}
$\aaa_1$ with $A_S\cong (\bz/6\bz)^4\times \bz/2\bz$
($H_{18,1}$ with the orbit of $\alpha_{2,1}$ or $\alpha_{4,1}$);
$2\aaa_1$ with $A_S\cong \bz/12\bz\times (\bz/6\bz)^2\times \bz/3\bz$
($H_{18,1}$ with the orbit of $\alpha_{1,1}$ or $\alpha_{3,1}$);
$3\aaa_1$ with $A_S\cong (\bz/6\bz)^3\times (\bz/2\bz)^2$
($H_{18,1}$ with the orbit of $\alpha_{2,2}$ or  $\alpha_{4,2}$);
$6\aaa_1$ with $A_S\cong \bz/12\bz\times (\bz/6\bz)^2$
($H_{18,1}$ with the orbit of $\alpha_{1,2}$ or $\alpha_{3,2}$).

\medskip

{\bf n=6,} $H\cong D_6$ ($|H|=6$, $i=1$):
$\rk N_H=14$ and $(N_H)^\ast/N_H\cong (\bz/6\bz)^2\times (\bz/3\bz)^3$.
{\it Degenerations of codimension $1$}:
$\aaa_1$ with $A_S\cong (\bz/6\bz)^3\times (\bz/3\bz)^2$
($H_{6,1}$ with the orbit of $\alpha_{1,1}$,  $\alpha_{2,1}$,
$\alpha_{3,1}$, $\alpha_{4,1}$, $\alpha_{5,1}$ or $\alpha_{6,1}$;
$H_{6,2}$ with the orbit of $\alpha_{2,1}$
or $\alpha_{4,1}$);
$2\aaa_1$ with $A_S\cong \bz/12\bz\times (\bz/3\bz)^4$
($H_{6,2}$ with the orbit of $\alpha_{1,1}$ or
$\alpha_{3,1}$);
$3\aaa_1$ with $A_S\cong (\bz/6\bz)^3\times \bz/3\bz$
($H_{6,1}$ with the orbit of $\alpha_{2,2}$ or $\alpha_{4,2}$;
$H_{6,2}$ with the orbit of $\alpha_{1,2}$, $\alpha_{2,2}$,
$\alpha_{3,2}$, $\alpha_{4,2}$, $\alpha_{5,2}$ or $\alpha_{6,2}$);
$6\aaa_1$ with $A_S\cong \bz/12\bz\times (\bz/3\bz)^3$
($H_{6,1}$ with the orbit of $\alpha_{1,2}$ or $\alpha_{3,2}$).

\medskip

{\bf n=3,} $H\cong C_2^2$ ($|H|=4$, $i=2$):
$\rk N_H=12$ and $(N_H)^\ast/N_H\cong (\bz/4\bz)^2\times (\bz/2\bz)^6$.
{\it Degenerations of codimension $1$:}
$\aaa_1$ with $A_S\cong (\bz/4\bz)^2\times (\bz/2\bz)^7$
($H_{3,1}$ with the orbit of $\alpha_{2,1}$, $\alpha_{4,1}$,
$\alpha_{2,2}$ or $\alpha_{4,2}$);
$2\aaa_1$ with $A_S\cong (\bz/4\bz)^3\times (\bz/2\bz)^4$
($H_{3,1}$ with the orbit of $\alpha_{1,1}$, $\alpha_{3,1}$,
$\alpha_{1,2}$, $\alpha_{3,2}$, $\alpha_{2,3}$ or $\alpha_{4,3}$);
$4\aaa_1$ with $A_S\cong \bz/8\bz \times (\bz/2\bz)^6$
($H_{3,1}$ with the orbit of $\alpha_{1,3}$ or $\alpha_{3,3}$).

\medskip

{\bf n=2,} $H\cong C_3$ ($|H|=3$, $i=1$):
$\rk N_H=12$ and $(N_H)^\ast/N_H\cong (\bz/3\bz)^6$.
{\it Degenerations of codimension $1$}:
$\aaa_1$ with $A_S\cong \bz/6\bz\times (\bz/3\bz)^5$
($H_{2,1}$ with the orbit of $\alpha_{1,1}$,  $\alpha_{2,1}$,
$\alpha_{3,1}$, $\alpha_{4,1}$, $\alpha_{5,1}$ or $\alpha_{6,1}$);
$3\aaa_1$ with $A_S\cong \bz/6\bz \times (\bz/3\bz)^4$
($H_{2,1}$ with the orbit of $\alpha_{1,2}$, $\alpha_{2,2}$,
$\alpha_{3,2}$, $\alpha_{4,2}$, $\alpha_{5,2}$ or $\alpha_{6,2}$).

\medskip

{\bf n=1,} $H\cong C_2$ ($|H|=2$, $i=1$):
$\rk N_H=8$ and $(N_H)^\ast/N_H\cong (\bz/2\bz)^8$.
{\it Degenerations of codimension $1$}:
$\aaa_1$ with $A_S\cong (\bz/2\bz)^9$
($H_{1,1}$ with the orbit of $\alpha_{2,1}$, $\alpha_{4,1}$,
$\alpha_{2,2}$,
$\alpha_{4,2}$, $\alpha_{2,3}$, $\alpha_{4,3}$, $\alpha_{2,4}$
or $\alpha_{4,4}$;
$H_{1,2}$ with the orbit of $\alpha_{2,1}$, $\alpha_{4,1}$,
$\alpha_{1,2}$,
$\alpha_{2,2}$, $\alpha_{3,2}$, $\alpha_{4,2}$, $\alpha_{5,2}$ or
$\alpha_{6,2}$);
$2\aaa_1$ with $A_S\cong \bz/4\bz\times (\bz/2\bz)^6$
($H_{1,1}$ with the orbit of $\alpha_{1,1}$, $\alpha_{1,2}$,
$\alpha_{1,3}$, $\alpha_{1,4}$, $\alpha_{3,1}$, $\alpha_{3,2}$,
$\alpha_{3,3}$  or $\alpha_{3,4}$;
$H_{1,2}$ with the orbit of $\alpha_{1,1}$, $\alpha_{1,3}$,
$\alpha_{6,3}$, $\alpha_{3,1}$, $\alpha_{2,3}$, $\alpha_{3,3}$,
$\alpha_{4,3}$ or $\alpha_{5,3}$).


\vskip1cm

{\bf Table 11. All possible markings by $N_{11}=N(A_{11}\oplus D_7\oplus E_6)$
of degenerations of
codimension one of  non-trivial finite symplectic automorphism groups
of K\"ahlerian K3 surfaces:}

\medskip

{\bf n=1,} $H\cong C_2$ ($|H|=2$, $i=1$):
$\rk N_H=8$ and $(N_H)^\ast/N_H\cong (\bz/2\bz)^8$.
{\it Degenerations of codimension $1$}:
$\aaa_1$ with $A_S\cong (\bz/2\bz)^9$
($H_{1,1}$ with the orbit of $\alpha_{6,1}$, $\alpha_{1,2}$,
$\alpha_{2,2}$, $\alpha_{3,2}$, $\alpha_{4,2}$, $\alpha_{5,2}$,
$\alpha_{2,3}$  or $\alpha_{4,3}$);
$2\aaa_1$ with $A_S\cong \bz/4\bz\times (\bz/2\bz)^6$
($H_{1,1}$ with the orbit of $\alpha_{1,1}$, $\alpha_{2,1}$,
$\alpha_{3,1}$,
$\alpha_{4,1}$, $\alpha_{5,1}$, $\alpha_{6,2}$, $\alpha_{1,3}$
or $\alpha_{3,3}$).


\vskip1cm

{\bf Table 9. All possible markings by $N_{9}=N(3D_8)$ of degenerations of
codimension one of  non-trivial finite symplectic automorphism groups
of K\"ahlerian K3 surfaces:}

\medskip

{\bf n=1,} $H\cong C_2$ ($|H|=2$, $i=1$):
$\rk N_H=8$ and $(N_H)^\ast/N_H\cong (\bz/2\bz)^8$.
{\it Degenerations of codimension $1$}:
$\aaa_1$ with $A_S\cong (\bz/2\bz)^9$
($H_{1,1}$ with the orbit of $\alpha_{1,3}$, $\alpha_{2,3}$,
$\alpha_{3,3}$, $\alpha_{4,3}$, $\alpha_{5,3}$, $\alpha_{6,3}$,
$\alpha_{7,3}$  or $\alpha_{8,3}$);
$2\aaa_1$ with $A_S\cong \bz/4\bz\times (\bz/2\bz)^6$
($H_{1,1}$ with the orbit of $\alpha_{1,1}$, $\alpha_{2,1}$,
$\alpha_{3,1}$,
$\alpha_{4,1}$, $\alpha_{5,1}$, $\alpha_{6,1}$, $\alpha_{7,1}$
or $\alpha_{8,1}$).


\vskip1cm

{\bf Table 8. All possible markings by $N_{8}=N(A_{15}\oplus D_9)$ of degenerations of
codimension one of non-trivial  finite symplectic automorphism groups
of K\"ahlerian K3 surfaces:}

\medskip

{\bf n=1,} $H\cong C_2$ ($|H|=2$, $i=1$):
$\rk N_H=8$ and $(N_H)^\ast/N_H\cong (\bz/2\bz)^8$.
{\it Degenerations of codimension $1$}:
$\aaa_1$ with $A_S\cong (\bz/2\bz)^9$
($H_{1,1}$ with the orbit of $\alpha_{8,1}$, $\alpha_{1,2}$,
$\alpha_{2,2}$, $\alpha_{3,2}$, $\alpha_{4,2}$, $\alpha_{5,2}$,
$\alpha_{6,2}$  or $\alpha_{7,2}$);
$2\aaa_1$ with $A_S\cong \bz/4\bz\times (\bz/2\bz)^6$
($H_{1,1}$ with the orbit of $\alpha_{1,1}$, $\alpha_{2,1}$,
$\alpha_{3,1}$,
$\alpha_{4,1}$, $\alpha_{5,1}$, $\alpha_{6,1}$, $\alpha_{7,1}$
or $\alpha_{8,2}$).


\vskip1cm

{\bf Table 7. All possible markings by $N_{7}=N(D_{10}\oplus 2E_7)$ of degenerations of
codimension one of  non-trivial finite symplectic automorphism groups
of K\"ahlerian K3 surfaces:}

\medskip

{\bf n=1,} $H\cong C_2$ ($|H|=2$, $i=1$):
$\rk N_H=8$ and $(N_H)^\ast/N_H\cong (\bz/2\bz)^8$.
{\it Degenerations of codimension $1$}:
$\aaa_1$ with $A_S\cong (\bz/2\bz)^9$
($H_{1,1}$ with the orbit of $\alpha_{1,1}$, $\alpha_{2,1}$,
$\alpha_{3,1}$, $\alpha_{4,1}$, $\alpha_{5,1}$, $\alpha_{6,1}$,
$\alpha_{7,1}$  or $\alpha_{8,1}$);
$2\aaa_1$ with $A_S\cong \bz/4\bz\times (\bz/2\bz)^6$
($H_{1,1}$ with the orbit of $\alpha_{9,1}$, $\alpha_{1,2}$,
$\alpha_{2,2}$,
$\alpha_{3,2}$, $\alpha_{4,2}$, $\alpha_{5,2}$, $\alpha_{6,2}$
or $\alpha_{7,2}$).


\vskip1cm

{\bf Table 6. All possible markings by $N_{6}=N(A_{17}\oplus E_7)$ of degenerations of
codimension one of  non-trivial finite symplectic automorphism groups
of K\"ahlerian K3 surfaces:}

\medskip

{\bf n=1,} $H\cong C_2$ ($|H|=2$, $i=1$):
$\rk N_H=8$ and $(N_H)^\ast/N_H\cong (\bz/2\bz)^8$.
{\it Degenerations of codimension $1$}:
$\aaa_1$ with $A_S\cong (\bz/2\bz)^9$
($H_{1,1}$ with the orbit of $\alpha_{9,1}$, $\alpha_{1,2}$,
$\alpha_{2,2}$, $\alpha_{3,2}$, $\alpha_{4,2}$, $\alpha_{5,2}$,
$\alpha_{6,2}$  or $\alpha_{7,2}$);
$2\aaa_1$ with $A_S\cong \bz/4\bz\times (\bz/2\bz)^6$
($H_{1,1}$ with the orbit of $\alpha_{1,1}$, $\alpha_{2,1}$,
$\alpha_{3,1}$,
$\alpha_{4,1}$, $\alpha_{5,1}$, $\alpha_{6,1}$, $\alpha_{7,1}$
or $\alpha_{8,1}$).


\vskip1cm

{\bf Table 3. All possible markings by $N_{3}=N(3E_8)$ of degenerations of
codimension one of non-trivial finite symplectic automorphism groups
of K\"ahlerian K3 surfaces:}

\medskip

{\bf n=1,} $H\cong C_2$ ($|H|=2$, $i=1$):
$\rk N_H=8$ and $(N_H)^\ast/N_H\cong (\bz/2\bz)^8$.
{\it Degenerations of codimension $1$}:
$\aaa_1$ with $A_S\cong (\bz/2\bz)^9$
($H_{1,1}$ with the orbit of $\alpha_{1,3}$, $\alpha_{2,3}$,
$\alpha_{3,3}$, $\alpha_{4,3}$, $\alpha_{5,3}$, $\alpha_{6,3}$,
$\alpha_{7,3}$  or $\alpha_{8,3}$);
$2\aaa_1$ with $A_S\cong \bz/4\bz\times (\bz/2\bz)^6$
($H_{1,1}$ with the orbit of $\alpha_{1,1}$, $\alpha_{2,1}$,
$\alpha_{3,1}$,
$\alpha_{4,1}$, $\alpha_{5,1}$, $\alpha_{6,1}$, $\alpha_{7,1}$
or $\alpha_{8,1}$).



\newpage

\section{Classification of types of degenerations\\
of finite symplectic automorphism groups\\
of K\"ahlerian K3 surfaces}
\label{sec:degALL}

Using Theorem \ref{theorem:degj}, we obtain the following classification of
types of degenerations of codimension $1$ of K\"ahlerian K3 surfaces.

\begin{theorem}
For non-trivial finite symplectic
automorphism groups $H$ of K\"ahlerian K3 surfaces $X$,
types of degenerations of codimension one $N_H\subset S$
(that is $\rk S=\rk N_H+1$ and $S<0$) are
the following and only the following which are given in
Table GEN below.
We give the type of the set of non-singular rational curves
of $X$ and the
isomorphism class of the discriminant group
$A_S=S^\ast/S$ of $S=S_X$.
The description of these degenerations using
markings by Niemeier lattices
$N_j$, $j=3,6,7,8,9,11\ -\ 23$, are given in
Theorem \ref{theorem:degj} above and its Tables j.
\label{theorem:degengeneral}
\end{theorem}

\vskip1cm

{\bf Table GEN. Classification of types of degenerations
of codimension one of  non-trivial finite symplectic
automorphism groups
of K\"ahlerian K3 surfaces:}

\medskip

{\bf n=81,} $H\cong M_{20}$ ($|H|=960$, $i=11357$):
$\rk N_H=19$. {\it No degenerations.}

\medskip

{\bf n=80,} $H\cong F_{384}$ ($|H|=384$, $i=18135$):
$\rk N_H=19$. {\it No degenerations.}

\medskip

{\bf n=79,} $H\cong {\mathfrak A}_{6}$ ($|H|=360$, $i=118$):
$\rk N_H=19$. {\it No degenerations.}

\medskip

{\bf n=78,} $H\cong {\mathfrak A}_{4,4}$ ($|H|=288$, $i=1026$):
$\rk N_H=19$. {\it No degenerations.}

\medskip

{\bf n=77,} $H\cong {T}_{192}$ ($|H|=192$, $i=1493$):
$\rk N_H=19$. {\it No degenerations.}

\medskip

{\bf n=76,} $H\cong {H}_{192}$ ($|H|=192$, $i=955$):
$\rk N_H=19$. {\it No degenerations.}

\medskip

{\bf n=75,} $H\cong 4^2{\mathfrak A}_4$ ($|H|=192$, $i=1023$):
$\rk N_H=18$, $(N_H)^\ast/N_H\cong (\bz/8\bz)^2\times
(\bz/2\bz)^2$. {\it Degenerations:}
$16\aaa_1$
with $A_S\cong (\bz/8\bz)\times (\bz/2\bz)^2$.

\medskip

{\bf n=74,} $H\cong L_2(7)$ ($|H|=168$, $i=42$):
$\rk N_H=19$. {\it No degenerations.}

\medskip


{\bf n=73,} $H\cong 2^4D_{10}$ ($|H|=160$, $i=234$):
The same as for $n=81$ above.

\medskip

{\bf n=72,} $H\cong {\mathfrak A}_4^2$ ($|H|=144$, $i=184$):
The same as for $n=78$ above.

\medskip

{\bf n=71,} $H\cong F_{128}$ ($|H|=128$, $i=931$):
The same as for $n=80$ above.

\medskip

{\bf n=70,} $H\cong {\mathfrak S}_5$ ($|H|=120$, $i=34$):
$\rk N_H=19$. {\it No degenerations.}

\medskip

{\bf n=69,} $H\cong (Q_8 * Q_8)\rtimes C_3$ ($|H|=96$, $i=204$):
The same as for $n=77$ above.

\medskip

{\bf n=68,} $H\cong 2^3D_{12}$, ($|H|=96$, $i=195$):
The same as for $n=78$ above.

\medskip

{\bf n=67,} $H\cong 4^2D_6$, ($|H|=96$, $i=64$):
The same as for $n=80$ above.

\medskip

{\bf n=66,} $H\cong 2^4 C_6$ ($|H|=96$, $i=70$):
The same as for $n=76$ above.

\medskip

{\bf n=65,} $H\cong2^4D_6$ ($|H|=96$, $i=227$):
$\rk N_H=18$,
$(N_H)^\ast/N_H\cong \bz/24\bz \times
\bz/4\bz\times (\bz/2\bz)^2$. {\it Degenerations:}
$4\aaa_1$ with $A_S\cong \bz/12\bz\times (\bz/4\bz)^2$;
$8\aaa_1$ with $A_S\cong \bz/24\bz\times (\bz/2\bz)^2$;
$12\aaa_{1}$ with $A_S\cong (\bz/4\bz)^3$.
$16\aaa_1$ with $A_S\cong \bz/12\bz\times (\bz/2\bz)^2$.

\medskip

{\bf n=64,} $H\cong 2^4 C_5$ ($|H|=80$, $i=49$):
The same as for $n=81$ above.

\medskip

{\bf n=63,} $H\cong M_9$ ($|H|=72$, $i=41$):
$\rk N_H=19$. {\it No degenerations.}

\medskip

{\bf n=62,} $H\cong N_{72}$ ($|H|=72$, $i=40$):
$\rk N_H=19$. {\it No degenerations.}

\medskip

{\bf n=61,} $H\cong {\mathfrak A}_{4,3}$ ($|H|=72$, $i=43$):
$\rk N_H=18$ and $(N_H)^\ast/N_H\cong (\bz/12\bz)^2 \times \bz/3\bz$.
{\it Degenerations:}
$3\aaa_1$ with $A_S\cong (\bz/12\bz)^2 \times \bz/2\bz$;
$12\aaa_1$ with $A_S\cong \bz/24\bz\times \bz/3\bz$.

\medskip

{\bf n=60,} $H\cong \Gamma_{26}a_2$ ($|H|=64$, $i=136$):
The same as for $n=80$ above.

\medskip

{\bf n=59,} $H\cong \Gamma_{23}a_2$ ($|H|=64$, $i=35$):
The same as for $n=80$ above.

\medskip

{\bf n=58,} $H\cong \Gamma_{22}a_1$ ($|H|=64$, $i=32$):
The same as for $n=80$ above.

\medskip

{\bf n=57,} $H\cong \Gamma_{13}a_1$ ($|H|=64$, $i=242$):
The same as for $n=75$ above.

\medskip

{\bf n=56,} $H\cong \Gamma_{25}a_1$ ($|H|=64$, $i=138$):
$\rk N_H=18$, $(N_H)^\ast/N_H\cong
\bz/8\bz\times (\bz/4\bz)^3$.
{\it Degenerations:}
$8\aaa_1$ with $A_S\cong \bz/8\bz\times (\bz/4\bz)^2$;
$16\aaa_1$ with $A_S\cong (\bz/4\bz)^3$.

\medskip

{\bf n=55,} $H\cong {\mathfrak A}_5$ ($|H|=60$, $i=5$):
$\rk N_H=18$ and $(N_H)^\ast/N_H\cong \bz/30\bz\times \bz/10\bz$.
{\it Degenerations:}
$\aaa_1$ with $A_S\cong \bz/30\bz\times \bz/10\bz\times \bz/2\bz$;
$5\aaa_1$ with $A_S\cong \bz/30\bz\times (\bz/2\bz)^2$;
$6\aaa_1$ with $A_S\cong \bz/20\bz\times \bz/5\bz$;
$10\aaa_1$ with $A_S\cong \bz/60\bz$;
$15\aaa_1$ with $A_S\cong \bz/10\bz \times (\bz/2\bz)^2$.

\medskip

{\bf n=54,} $H\cong T_{48}$ ($|H|=48$, $i=29$):
$\rk N_H=19$ and $(N_H)^\ast/N_H\cong
\bz/24\bz\times \bz/8\bz\times \bz/2\bz$.
{\it No degenerations.}

\medskip

{\bf n=53,} $H\cong 2^2Q_{12}$ ($|H|=48$, $i=30$):
The same as for $n=78$ above.

\medskip

 {\bf n=52,} $H\cong 2^2(C_2\times C_6)$ ($|H|=48$, $i=49$):
The same as for $n=78$ above.

\medskip

{\bf n=51,} $H\cong C_2\times {\mathfrak S}_4$ ($|H|=48$, $i=48$):
$\rk N_H=18$, $(N_H)^\ast/N_H\cong (\bz/12\bz)^2\times (\bz/2\bz)^2$.
{\it Degenerations:}
$2\aaa_1$ with $A_S\cong (\bz/12\bz)^2\times  \bz/4\bz$;
$4\aaa_1$ with $A_S\cong \bz/24\bz \times \bz/6\bz \times \bz/2\bz$;
$6\aaa_1$ with $A_S\cong \bz/12\bz \times (\bz/4\bz)^2$;
$8\aaa_1$ with $A_S\cong \bz/12\bz\times  \bz/6\bz\times \bz/2\bz$;
$12\aaa_1$ with $A_S\cong \bz/24\bz\times (\bz/2\bz)^2$.

\medskip

{\bf n=50,} $H\cong 4^2C_3$ ($|H|=48$, $i=3$):
The same as for $n=75$ above.

\medskip

{\bf n=49,} $H\cong 2^4 C_3$ ($|H|=48$, $i=50$):
$\rk N_H=17$, $(N_H)^\ast/N_H \cong
\bz/24\bz\times (\bz/2\bz)^4$. {\it Degenerations of codimension $1$:}
$4\aaa_1$ with $A_S\cong \bz/12\bz \times \bz/4\bz\times (\bz/2\bz)^2$;
$12\aaa_1$ with $A_S\cong (\bz/4\bz)^2 \times (\bz/2\bz)^2$;
$16\aaa_1$ with $A_S\cong \bz/6\bz \times (\bz/2\bz)^3$.

\medskip

{\bf n=48,} $H\cong {\mathfrak S}_{3,3}$ ($|H|=36$, $i=10$):
$\rk N_H=18$, $(N_H)^\ast/N_H \cong
\bz/18\bz\times \bz/6\bz\times (\bz/3\bz)^2$.
{\it Degenerations:}
$3\aaa_1$ with $A_S\cong \bz/18\bz \times (\bz/6\bz)^2$;
$6\aaa_1$ with $A_S\cong \bz/36\bz \times (\bz/3\bz)^2$;
$9\aaa_1$ with $A_S\cong (\bz/6\bz)^3$.

\medskip

{\bf n=47,} $H\cong C_3\times {\mathfrak A}_4$ ($|H|=36$, $i=11$):
The same as for $n=61$ above.

\medskip

\medskip

{\bf n=46,} $H\cong 3^2 C_4$ ($|H|=36$, $i=9$):
$\rk N_H=18$ and $(N_H)^\ast/N_H\cong \bz/18\bz\times
\bz/6\bz\times \bz/3\bz$.
{\it Degenerations:}
$6\aaa_1$ with $A_S\cong \bz/36\bz \times \bz/3\bz$;
$9\aaa_1$ with $A_S\cong (\bz/6\bz)^2 \times \bz/2\bz$;
$9\aaa_2$ with $A_S\cong \bz/6\bz \times \bz/3\bz$.

\medskip

{\bf n=45,} $H\cong \Gamma_{6}a_2$ ($|H|=32$, $i=44$):
The same as for $n=80$ above.

\medskip

{\bf n=44,} $H\cong \Gamma_{3}e$  ($|H|=32$, $i=11$):
The same as for $n=80$ above.

\medskip

{\bf n=43,} $H\cong \Gamma_{7}a_2$ ($|H|=32$, $i=7$):
The same as for $n=80$ above.

\medskip

{\bf n=42,} $H\cong \Gamma_4c_2$ {$|H|=32$, $i=31$):
The same as for $n=75$ above.

\medskip

{\bf n=41,} $H\cong \Gamma_7 a_1$, ($|H|=32$, $i=6$):
The same as for $n=56$ above.

\medskip

{\bf n=40,} $H\cong Q_8*Q_8$ ($|H|=32$, $i=49$):
$\rk N_H=17$, $(N_H)^\ast/N_H\cong (\bz/4\bz)^5$.
{\it Degenerations of codimension $1$:}
$8\aaa_1$ with $A_S\cong (\bz/4\bz)^4$.

\medskip

{\bf n=39,} $H\cong 2^4C_2$ ($|H|=32$, $i=27$):
$\rk N_H=17$, $(N_H)^\ast/N_H\cong \bz/8\bz\times
(\bz/4\bz)^2\times (\bz/2\bz)^2$.
{\it Degenerations of codimension $1$:}
$4\aaa_1$ with $A_S\cong (\bz/4\bz)^4$;
$8\aaa_1$ with $A_S\cong \bz/8\bz\times \bz/4\bz\times (\bz/2\bz)^2$;
$16\aaa_1$ with $A_S\cong (\bz/4\bz)^2\times (\bz/2\bz)^2$.

\medskip

{\bf n=38,} $H\cong T_{24}$ ($|H|=24$, $i=3$):
The same as for $n=54$ above.

\medskip

{\bf n=37,} $H\cong T_{24}$, ($|H|=24$, $i=3$):
The same as for $n=77$ above.

\medskip

{\bf n=36,} $H\cong C_3\rtimes D_8$ ($|H|=24$, $i=8$):
The same as for $n=61$ above.

\medskip

{\bf n=35,} $H\cong C_2\times {\mathfrak A}_4$, ($|H|=24$, $i=13$):
The same as for $n=51$ above.

\medskip

{\bf n=34}, $H\cong {\mathfrak S}_4$ ($|H|=24$, $i=12$):
$\rk N_H=17$ and $(N_H)^\ast/N_H\cong
(\bz/12\bz)^2\times \bz/4\bz$.
{\it Degenerations of codimension $1$:}
$\aaa_1$ with $A_S\cong (\bz/12\bz)^2\times \bz/4\bz\times \bz/2\bz$;
$2\aaa_1$ with $A_S\cong (\bz/12\bz)^2\times (\bz/2\bz)^2$;
$3\aaa_1$ with $A_S\cong \bz/12\bz\times (\bz/4\bz)^2\times \bz/2\bz$;
$4\aaa_1$ with $A_S\cong \bz/24\bz\times \bz/12\bz$;
$6\aaa_1$ with $A_S\cong \bz/12\bz\times \bz/4\bz\times (\bz/2\bz)^2$;
$8\aaa_1$ with $A_S\cong (\bz/12\bz)^2$;
$12\aaa_1$ with $A_S\cong \bz/24\bz\times \bz/4\bz$;
$6\aaa_2$ with $A_S\cong \bz/12\bz \times \bz/4\bz$.

\medskip

{\bf n=33,} $H\cong C_7\rtimes C_3$ ($|H|=21$, $i=1$):
$\rk N_H=18$ and $(N_H)^\ast/N_H\cong (\bz/7\bz)^3$.
{\it Degenerations:}
$7\aaa_1$ with $A_S\cong \bz/14\bz\times \bz/7\bz$.

\medskip

{\bf n=32,} $H\cong Hol(C_5)$ ($|H|=20$, $i=3$):
$\rk N_H=18$ and $(N_H)^\ast/N_H \cong (\bz/10\bz)^2 \times \bz/5\bz$.
{\it Degenerations:}
$2\aaa_1$ with $A_S\cong \bz/20\bz\times (\bz/5\bz)^2$;
$5\aaa_1$ with $A_S\cong (\bz/10\bz)^2\times \bz/2\bz$;
$10\aaa_1$ with $A_S\cong \bz/20\bz\times \bz/5\bz$;
$5\aaa_2$ with $A_S\cong \bz/10\bz \times \bz/5\bz$.

\medskip

{\bf n=31,} $H\cong C_3\times D_6$ ($|H|=18$, $i=3$):
The same as for $n=48$ above.

\medskip

{\bf n=30,} $H\cong {\mathfrak A}_{3,3}$ ($|H|=18$, $i=4$):
$\rk N_H=16$ and $(N_H)^\ast/N_H\cong \bz/9\bz\times (\bz/3\bz)^4$.
{\it Degenerations of codimension $1$:}
$3\aaa_1$ with $A_S\cong \bz/18\bz\times (\bz/3\bz)^3$;
$9\aaa_1$ with $A_S\cong \bz/6\bz\times (\bz/3\bz)^3$.

\medskip

{\bf n=29,} $H\cong Q_{16}$ ($|H|=16$, $i=9$):
The same as for $n=80$ above.

\medskip

{\bf n=28,} $H\cong \Gamma_2 d$ ($|H|=16$, $i=6$):
The same as for $n=80$ above.

\medskip

{\bf n=27,} $H\cong C_2\times Q_8$ ($|H|=16$, $i=12$):
The same as for $n=75$ above.

\medskip

{\bf n=26,} $H\cong SD_{16}$ ($|H|=16$, $i=8$):
$\rk N_H=18$ and $(N_H)^\ast/N_H\cong (\bz/8\bz)^2\times
\bz/4\bz\times \bz/2\bz$.
{\it Degenerations:}
$8\aaa_1$ with $A_S\cong \bz/16\bz\times \bz/4\bz\times \bz/2\bz$;
$2\aaa_2$ with $A_S\cong (\bz/8\bz)^2 \times \bz/2\bz$.

\medskip

{\bf n=25,} $H\cong C_4^2$ ($|H|=16$, $i=2$):
The same as for $n=75$ above.

\medskip

{\bf n=24,} $H\cong Q_8 * C_4$ ($|H|=16$, $i=13$):
The same as for $n=40$ above.

\medskip

{\bf n=23,} $H\cong \Gamma_2 c_1$ ($|H|=16$, $i=3$):
The same as for $n=39$ above.

\medskip

{\bf n=22,} $H\cong C_2\times D_8$ ($|H|=16$, $i=11$):
$\rk N_H=16$, $(N_H)^\ast/N_H\cong (\bz/4\bz)^4\times (\bz/2\bz)^2$.
{\it Degenerations of codimension $1$:}
$2\aaa_1$ with $A_S\cong (\bz/4\bz)^5$;
$4\aaa_1$ with $A_S\cong \bz/8\bz\times (\bz/4\bz)^2\times (\bz/2\bz)^2$;
$8\aaa_1$ with $A_S\cong (\bz/4\bz)^3\times (\bz/2\bz)^2$.

\medskip

{\bf n=21,} $H\cong C_2^4$ ($|H|=16$, $i=14$):
$\rk N_H=15$, $(N_H)^\ast/N_H \cong \bz/8\bz\times (\bz/2\bz)^6$.
{\it Degenerations of codimension $1$:}
$4\aaa_1$ with $A_S\cong (\bz/4\bz)^2\times (\bz/2\bz)^4$;
$16\aaa_1$ with $A_S\cong (\bz/2\bz)^6$.

\medskip

{\bf n=20,} $H\cong Q_{12}$ ($|H|=12$, $i=1$):
The same as for $n=61$ above.

\medskip

{\bf n=19,} $H\cong C_2\times C_6$ ($|H|=12$, $i=5$):
The same as for $n=61$ above.

\medskip

{\bf n=18,} $H\cong D_{12}$ ($|H|=12$, $i=4$):
$\rk N_H=16$ and $(N_H)^\ast/N_H\cong (\bz/6\bz)^4$.
{\it Degenerations of codimension $1$:}
$\aaa_1$ with $A_S\cong (\bz/6\bz)^4\times \bz/2\bz$;
$2\aaa_1$ with $A_S\cong \bz/12\bz\times (\bz/6\bz)^2\times \bz/3\bz$;
$3\aaa_1$ with $A_S\cong (\bz/6\bz)^3\times (\bz/2\bz)^2$;
$6\aaa_1$ with $A_S\cong \bz/12\bz\times (\bz/6\bz)^2$.

\medskip

{\bf n=17,} $H\cong {\mathfrak A}_4$ ($|H|=12$, $i=3$):
$\rk N_H=16$ and $(N_H)^\ast/N_H\cong (\bz/12\bz)^2\times (\bz/2\bz)^2$.
{\it Degenerations of codimension $1$:}
$\aaa_1$ with $A_S\cong (\bz/12\bz)^2\times (\bz/2\bz)^3$;
$3\aaa_1$ with $A_S\cong \bz/12\bz\times \bz/4\bz\times (\bz/2\bz)^3$;
$4\aaa_1$ with $A_S\cong \bz/24\bz\times \bz/6\bz\times \bz/2\bz$;
$6\aaa_1$ with $A_S\cong \bz/12\bz\times (\bz/4\bz)^2$;
$12\aaa_1$ with $A_S\cong \bz/24\bz\times (\bz/2\bz)^2$.

\medskip

{\bf n=16,} $H\cong D_{10}$ ($|H|=10$, $i=1$):
$\rk N_H=16$ and $(N_H)^\ast/N_H\cong (\bz/5\bz)^4$.
{\it Degenerations of codimension $1$:}
$\aaa_1$ with $A_S\cong \bz/10\bz\times (\bz/5\bz)^3$;
$5\aaa_1$ with $A_S\cong \bz/10\bz\times (\bz/5\bz)^2$.

\medskip

{\bf n=15,} $H\cong C_3^2$ ($|H|=9$, $i=2$):
The same as for $n=30$ above.

\medskip

{\bf n=14,} $H\cong C_8$ ($|H|=8$, $i=1$):
The same as for $n=26$ above.

\medskip

{\bf n=13,} $H\cong Q_8$ ($|H|=8$, $i=4$):
The same as for $n=40$ above.

\medskip

{\bf n=12,} $H\cong Q_8$ ($|H|=8$, $i=4$):
$\rk N_H=17$, $(N_H)^\ast/N_H \cong
(\bz/8\bz)^2\times (\bz/2\bz)^3$.
{\it Degenerations of codimension $1$:}
$8\aaa_1$ with $A_S\cong \bz/16\bz\times (\bz/2\bz)^3$;
$\aaa_2$ with $A_S\cong (\bz/8\bz)^2 \times (\bz/2\bz)^2$.

\medskip

{\bf n=11,} $H\cong C_2\times C_4$ ($|H|=8$, $i=2$):
The same as for $n=22$ above.

\medskip

{\bf n=10,} $H\cong D_8$ ($|H|=8$, $i=3$):
$\rk N_H=15$ and $(N_H)^\ast/N_H\cong
(\bz/4\bz)^5$.
{\it Degenerations of codimension $1$:}
$\aaa_1$ with $A_S\cong (\bz/4\bz)^5\times \bz/2\bz$;
$2\aaa_1$ with $A_S\cong (\bz/4\bz)^4\times (\bz/2\bz)^2$;
$4\aaa_1$ with $A_S\cong \bz/8\bz\times (\bz/4\bz)^3$;
$8\aaa_1$ with $A_S\cong (\bz/4\bz)^4$;
$2\aaa_2$ with $A_S\cong (\bz/4\bz)^4$.

\medskip

{\bf n=9,} $H\cong C_2^3$ ($|H|=8$, $i=5$):
$\rk N_H=14$, $(N_H)^\ast/N_H \cong (\bz/4\bz)^2\times (\bz/2\bz)^6$.
{\it Degenerations of codimension $1$}:
$2\aaa_1$ with $A_S\cong (\bz/4\bz)^3\times (\bz/2\bz)^4$;
$4\aaa_1$ with $A_S\cong \bz/8\bz\times (\bz/2\bz)^6$;
$8\aaa_1$ with $A_S\cong \bz/4\bz\times (\bz/2\bz)^6$.

\medskip

{\bf n=8,} $H\cong C_7$ ($|H|=7$, $i=1$):
The same as for $n=33$ above.

\medskip

{\bf n=7,} $H\cong C_6$ ($|H|=6$, $i=2$):
The same as for $n=18$ above.

\medskip

{\bf n=6,} $H\cong D_6$ ($|H|=6$, $i=1$):
$\rk N_H=14$ and $(N_H)^\ast/N_H\cong (\bz/6\bz)^2\times (\bz/3\bz)^3$.
{\it Degenerations of codimension $1$}:
$\aaa_1$ with $A_S\cong (\bz/6\bz)^3\times (\bz/3\bz)^2$;
$2\aaa_1$ with $A_S\cong \bz/12\bz\times (\bz/3\bz)^4$;
$3\aaa_1$ with $A_S\cong (\bz/6\bz)^3 \times \bz/3\bz$;
$6\aaa_1$ with $A_S\cong \bz/12\bz\times (\bz/3\bz)^3$.

\medskip

{\bf n=5,} $H\cong C_5$ ($|H|=5$, $i=1$):
The same as for $n=16$ above.

\medskip

{\bf n=4,} $H\cong C_4$ ($|H|=4$, $i=1$):
$\rk N_H=14$ and $(N_H)^\ast/N_H \cong
(\bz/4\bz)^4\times (\bz/2\bz)^2$.
{\it Degenerations of codimension $1$}:
$\aaa_1$ with $A_S\cong (\bz/4\bz)^4 \times (\bz/2\bz)^3$;
$2\aaa_1$ with $A_S\cong (\bz/4\bz)^5$;
$4\aaa_1$ with $A_S\cong \bz/8\bz \times (\bz/4\bz)^2\times  (\bz/2\bz)^2$;
$\aaa_2$ with $A_S\cong (\bz/4\bz)^4 \times \bz/2\bz$;

\medskip

{\bf n=3,} $H\cong C_2^2$ ($|H|=4$, $i=2$):
$\rk N_H=12$ and $(N_H)^\ast/N_H\cong (\bz/4\bz)^2\times (\bz/2\bz)^6$.
{\it Degenerations of codimension $1$}:
$\aaa_1$ with $A_S\cong (\bz/4\bz)^2\times (\bz/2\bz)^7$;
$2\aaa_1$ with $A_S\cong (\bz/4\bz)^3\times (\bz/2\bz)^4$;
$4\aaa_1$ with $A_S\cong \bz/8\bz\times (\bz/2\bz)^6$.

\medskip

{\bf n=2,} $H\cong C_3$ ($|H|=3$, $i=1$):
$\rk N_H=12$ and $(N_H)^\ast/N_H\cong (\bz/3\bz)^6$.
{\it Degenerations of codimension $1$}:
$\aaa_1$ with $A_S\cong \bz/6\bz\times (\bz/3\bz)^5$;
$3\aaa_1$ with $A_S\cong \bz/6\bz \times (\bz/3\bz)^4$.

\medskip

{\bf n=1,} $H\cong C_2$ ($|H|=2$, $i=1$):
$\rk N_H=8$ and $(N_H)^\ast/N_H\cong (\bz/2\bz)^8$.
{\it Degenerations of codimension $1$}:
$\aaa_1$ with $A_S\cong (\bz/2\bz)^9$;
$2\aaa_1$ with $A_S\cong \bz/4\bz\times (\bz/2\bz)^6$.


\section{Conjecture and some remarks.}
\label{sec:Conj}

Theorem \ref{theorem:degengeneral} shows that for
a fixed type (defined by $n$)
of abstract finite symplectic group of automorphisms of K3 and for
a fixed type of degeneration
of codimension $1$ (type of Dynkin diagram), the
discriminant group $A_S$
of the corresponding lattice $S$ is always the same.
It is natural to suppose
that the more strong statement is valid: that the isomorphism class of
the lattice $S$ is defined uniquely.
But, exact calculations and considerations show that it is valid only with
some few exceptions given in the Conjecture below.

\begin{conjecture}
For a fixed type (defined by $n$)
of abstract finite symplectic group $H$ of automorphisms of K\"ahlerian K3 surfaces
and for a fixed type of degeneration 
of codimension $1$ (type of Dynkin diagram) of K\"ahlerian
K3 surfaces, the corresponding lattice $S=[N_H,\alpha_{i,j}]_{pr}$
is unique up to isomorphisms of lattices.

The same can be formulated geometrically. Let $X$ be a K\"ahlerian
K3 surface with $S_X<0$. Let $H=Aut^0 X$ be the group of symplectic
automorphisms of $X$ and $P(X)$ be a the set
of classes of non-singular rational curves on $X$. Let
$N_H=((S_X)^H)^\perp_{S_X}$ be the coinvariant sublattice. Assume that
$S_X$ is generated by $N_H$ and $P(X)$ up to finite index, and
$\rk S_X=\rk N_H+1$. Then the isomorphism class of the lattice $S=S_X$
is defined uniquely by the type of $H$ as abstract group (equivalent to $n$)
and the type of the Dynkin diagram of $P(X)$.

But, it is valid with the following and only the following two exceptions:

(1) {\bf n=34} (equivalently, $H\cong {\mathfrak S}_4$) and the degeneration of
the type $6\aaa_1$;

(2) {\bf n=10} (equivalently, $H\cong {D_8}$) and the degeneration of the type $2\aaa_1$.

In both these cases there are exactly two isomorphism classes of the
lattices $S$.
\label{conj}
\end{conjecture}

Obviously, the conjecture is valid for a trivial $H$.
Theorem  \ref{theorem:degengeneral} shows that the conjecture
is valid for the discriminant group $A_S$ of the lattice $S$.
We hope to present proofs of Conjecture \ref{conj} and of more delicate statements
in further publications. They will give classification of degenerations 
of codimension one of K\"ahlerian K3 surfaces with finite symplectic automorphism groups
(not of types of degenerations only).

\medskip

From Theorems \ref{theorem:degj} and \ref{theorem:degengeneral} 
we also obtain:

\begin{corollary} All types of degenerations
of codimension $1$ of K\"ahlerian K3 surfaces with finite
symplectic automorphism
groups defined by the invariant $n$ of the group and by the
type of Dynkin diagram
can be obtained from markings by Niemeier lattices
$N_{23}=N(24A_1)$ and $N_{22}=N(12A_2)$.
Degenerations of codimension $1$ marked by
other Niemeier lattices $N_j$, $j=1,\dots 21$,
give the same types
of degenerations.
\end{corollary}


\section{Appendix: Programs}
\label{sec:Appendix}

Here we give Programs 5 and 6 for GP/PARI Calculator, Version 2.7.0 which were used
for calculations above. They also include Programs 1 - 4 from
\cite{Nik7} and \cite{Nik8}.

\medskip

Program5:  niemeier$\backslash$general2a.txt

\medskip

\noindent
$\backslash\backslash$for a Niemeier lattice N\_{}i given by root matrix r

\noindent
$\backslash\backslash$and cord matrix cord\hfill

\noindent
$\backslash\backslash$and subgroup H$\backslash$subset A(N\_{}i)\hfill

\noindent
$\backslash\backslash$and its orbits ORB matrix, each line gives\hfill

\noindent
$\backslash\backslash$orbit of length\  $>$\ 1 with 0 remaining elements\hfill

\noindent
$\backslash\backslash$it calculates all additional\hfill

\noindent
$\backslash\backslash$1-elements orbits to matrix ORBF and prints them\hfill

\noindent
$\backslash\backslash$(1-elements orbits the last)\hfill

\noindent
$\backslash\backslash$it calculates coinvariant sublattice N\_{}H\hfill

\noindent
$\backslash\backslash$by its rational generators, and checks if\hfill

\noindent
$\backslash\backslash$N\_{}H has primitive embedding to L\_{}K3\hfill

\noindent
$\backslash\backslash$Then it prints invariants of the discriminant group of N\_{}H\hfill

\noindent
sORB=matsize(ORB);m1=0;\hfill

\noindent
for(k1=1,sORB[1],for(k2=1,sORB[2],$\backslash$\hfill

\noindent
if(ORB[k1,k2]==0,,m1=m1+1)));\hfill

\noindent
ORBF=matrix(sORB[1]+(24-m1),sORB[2]);\hfill

\noindent
for(k=1,sORB[1],ORBF[k,]=ORB[k,]);l=sORB[1];\hfill

\noindent
for(t=1,24,mu=0;for(k1=1,sORB[1],for(k2=1,sORB[2],$\backslash$\hfill

\noindent
if(ORB[k1,k2]!=t,,mu=1)));if(mu==1,,l=l+1;ORBF[l,1]=t));$\backslash$\hfill

\noindent
print("ORBF=",ORBF);$\backslash$\hfill

\noindent
SUBL0=matrix(24,24);alpha=0;$\backslash$\hfill

\noindent
for(k1=1,sORB[1],for(k2=1,sORB[2]-1,$\backslash$\hfill

\noindent
if(ORB[k1,k2+1]$>$0,alpha=alpha+1;$\backslash$\hfill

\noindent
SUBL0[,alpha]=R[,ORB[k1,k2]]-R[,ORB[k1,k2+1]])));$\backslash$\hfill

\noindent
sORBF=matsize(ORBF);$\backslash$\hfill

\noindent
SUBL=SUBL0;a=matrix(24,24+matsize(cord)[1]);$\backslash$\hfill

\noindent
for(i=1,24,a[i,i]=1);for(i=1,matsize(cord)[1],a[,24+i]=cord[i,]$\widetilde{\ }\,$);$\backslash$\hfill

\noindent
L=a;N=SUBL;ggg=gcd(N);N1=N/ggg;M=L;gg=gcd(M);M1=M/gg;$\backslash$\hfill

\noindent
ww=matsnf(M1,1);uu=ww[1];vv=ww[2];dd=ww[3];$\backslash$\hfill

\noindent
mm=matsize(dd)[1];nn=matsize(dd)[2];$\backslash$\hfill

\noindent
nnn=nn;for(i=1,nn,if(dd[,i]==0,nnn=nnn-1));VV=matrix(nn,nnn);$\backslash$\hfill

\noindent
nnnn=0;for(i=1,nn,if(dd[,i]==0,,nnnn=nnnn+1;VV[,nnnn]=vv[,i]));$\backslash$\hfill

\noindent
M2=M1*VV;MM=M2*gg;$\backslash$\hfill

\noindent
kill(gg);kill(M1);kill(ww);kill(uu);kill(vv);kill(dd);kill(mm);$\backslash$\hfill

\noindent
kill(nn);kill(nnn);kill(nnnn);kill(M2);$\backslash$\hfill

\noindent
L1=MM;kill(VV);N2=L1\^{}-1*N1;\

\noindent
ww=matsnf(N2,1);uu=ww[1];vv=ww[2];dd=ww[3];$\backslash$\hfill

\noindent
N3=N2*vv;mm=matsize(dd)[1];nn=matsize(dd)[2];$\backslash$\hfill

\noindent
nnn=nn;for(i=1,nn,if(dd[,i]==0,nnn=nnn-1));$\backslash$\hfill

\noindent
N4=matrix(mm,nnn);nnnn=0;$\backslash$\hfill

\noindent
for(i=1,nn,if(dd[,i]==0,,nnnn=nnnn+1;$\backslash$\hfill

\noindent
ddd=gcd(dd[,i]);N4[,nnnn]=N3[,i]/ddd));$\backslash$\hfill

\noindent
Npr=L1*N4;$\backslash$\hfill

\noindent
kill(ggg);kill(N1);kill(M);kill(L1);kill(MM);$\backslash$\hfill

\noindent
kill(N2);kill(ww);kill(uu);kill(vv);kill(dd);kill(N3);$\backslash$\hfill

\noindent
kill(mm);kill(nn);kill(nnn);kill(nnnn);kill(ddd);kill(N4);$\backslash$\hfill

\noindent
SUBLpr1=Npr;$\backslash$\hfill

\noindent
R=r;B=SUBLpr1;$\backslash$\hfill

\noindent
l=B$\widetilde{\ }\,$*R*B;$\backslash$\hfill

\noindent
ww=matsnf(l,1);uu=ww[1];vv=ww[2];dd=ww[3];$\backslash$\hfill

\noindent
nn=matsize(l)[1];nnn=nn;for(i=1,nn,if(dd[i,i]==0,nnn=nnn-1));$\backslash$\hfill

\noindent
b=matrix(nn,nnn,X,Y,vv[X,Y+nn-nnn]);$\backslash$\hfill

\noindent
ll=b$\widetilde{\ }\,$*l*b;$\backslash$\hfill

\noindent
d=vector(nnn,X,dd[X+nn-nnn,X+nn-nnn]);$\backslash$\hfill

\noindent
kill(ww);kill(uu);kill(vv);kill(dd);$\backslash$\hfill

\noindent
kill(nn);kill(nnn);$\backslash$\hfill

\noindent
BB=B*b;G=BB$\widetilde{\ }\,$*R*BB;D=d;$\backslash$\hfill

\noindent
SUBLpr=BB;DSUBLpr=D;rSUBLpr=G;$\backslash$\hfill

\noindent
RANK=matsize(DSUBLpr)[2];$\backslash$\hfill

\noindent
RANKA=0;for(k1=1,RANK,if(DSUBLpr[k1]$>$1,RANKA=RANKA+1));$\backslash$\hfill

\noindent
\&if(RANK$>$19$||$RANK+RANKA$>$22,,$\backslash$\hfill

\noindent
DISCGRS=vector(RANKA,X,DSUBLpr[X]);$\backslash$\hfill

\noindent
if(RANK+RANKA$<$22,print("DISCGRS=",DISCGRS),$\backslash$\hfill

\noindent
pp=factor(DISCGRS[RANKA]);$\backslash$\hfill

\noindent
u1=2;if(pp[1,1]$>$2,u1=1);if(pp[1,1]==2\&\&pp[1,2]$>$1,u1=1);$\backslash$\hfill

\noindent
if(pp[1,1]==2\&\&pp[1,2]==1,delta=0;$\backslash$\hfill

\noindent
for(k=1,RANKA,if(type(rSUBLpr[k,k]/4)!="t\_{}INT",delta=1)));$\backslash$\hfill

\noindent
if(delta==0,u1=1);gamm=0;$\backslash$\hfill

\noindent
BBB=matrix(RANK-RANKA,RANK-RANKA,X,Y,rSUBLpr[X+RANKA,Y+RANKA]);$\backslash$\hfill

\noindent
dBBB=matdet(BBB);$\backslash$\hfill

\noindent
for(t1=u1,matsize(pp)[1],if(pp[t1,1]==2\&\&$\backslash$\hfill

\noindent
Mod(dBBB,8)!=Mod(1,8)\&\&Mod(dBBB,8)!=Mod(-1,8),gamm=1);$\backslash$\hfill

\noindent
if(pp[t1,1]!=2\&\&kronecker((-1)\^{}(RANK+1)*dBBB,pp[t1,1])!=1,gamm=1));
$\backslash$\hfill

\noindent
if(gamm==0,print("DISCGRS=",DISCGRS))));\hfill


\medskip

Program 6: niemeier$\backslash$degen1a.txt

\medskip

\noindent
$\backslash\backslash$for a Niemeier lattice N\_{}i given by root matrix r\hfill

\noindent
$\backslash\backslash$and cord matrix cord

\noindent
$\backslash\backslash$and sugroup H$\backslash$subset A(N\_{}i)\hfill

\noindent
$\backslash\backslash$and its orbits ORB matrix, each line gives orbit of length $>$ 1,\hfill

\noindent
$\backslash\backslash$with 0 remaining elements\hfill

\noindent
$\backslash\backslash$it calculates all additional\hfill

\noindent
$\backslash\backslash$1-elements orbits to matrix ORBF and prints them\hfill

\noindent
$\backslash\backslash$(1-elements orbits the last)\hfill

\noindent
$\backslash\backslash$it calculates coinvariant sublattice N\_{}H by
its rational generators,\hfill

\noindent
$\backslash\backslash$for each first element k\_{}i of i line of ORBF\hfill

\noindent
$\backslash\backslash$it calculates invariants of the degeneration\hfill

\noindent
$\backslash\backslash$sublattice N(H,i)=[N\_{}H,i] using program general1.txt\hfill

\noindent
$\backslash\backslash$and checks if\hfill

\noindent
$\backslash\backslash$N(H,i) has primitive embedding to L\_{}K3\hfill

\noindent
$\backslash\backslash$Then it prints k\_{}i and invariants of the discriminant\hfill

\noindent
$\backslash\backslash$group of N(H,i)\hfill

\noindent
$\backslash\backslash$it prints: k\_{}i,\hfill

\noindent
$\backslash\backslash$and invariants of N(H\_{}i)$\backslash$subset N(H\_i)\^{}$
\backslash$ast\hfill

\noindent
sORB=matsize(ORB);m1=0;\hfill

\noindent
for(k1=1,sORB[1],for(k2=1,sORB[2],$\backslash$\hfill

\noindent
if(ORB[k1,k2]==0,,m1=m1+1)));\hfill

\noindent
ORBF=matrix(sORB[1]+(24-m1),sORB[2]);\hfill

\noindent
for(k=1,sORB[1],ORBF[k,]=ORB[k,]);l=sORB[1];\hfill

\noindent
for(t=1,24,mu=0;for(k1=1,sORB[1],for(k2=1,sORB[2],$\backslash$\hfill

\noindent
if(ORB[k1,k2]!=t,,mu=1)));if(mu==1,,l=l+1;ORBF[l,1]=t));$\backslash$\hfill

\noindent
print("ORBF=",ORBF);$\backslash$\hfill

\noindent
SUBL0=matrix(24,24);alpha=0;$\backslash$\hfill

\noindent
for(k1=1,sORB[1],for(k2=1,sORB[2]-1,$\backslash$\hfill

\noindent
if(ORB[k1,k2+1]$>$0,alpha=alpha+1;$\backslash$\hfill

\noindent
SUBL0[,alpha]=R[,ORB[k1,k2]]-R[,ORB[k1,k2+1]])));$\backslash$\hfill

\noindent
sORBF=matsize(ORBF);$\backslash$\hfill

\noindent
for(n=1,sORBF[1],t=ORBF[n,1];SUBL=SUBL0;SUBL[t,alpha+1]=1;
$\backslash$\hfill

\noindent
a=matrix(24,24+matsize(cord)[1]);$\backslash$\hfill

\noindent
for(i=1,24,a[i,i]=1);for(i=1,matsize(cord)[1],a[,24+i]=cord[i,]$\widetilde{\ }\,$);
$\backslash$\hfill

\noindent
L=a;N=SUBL;$\backslash$\hfill

\noindent
ggg=gcd(N);N1=N/ggg;M=L;$\backslash$\hfill

\noindent
gg=gcd(M);M1=M/gg;$\backslash$\hfill

\noindent
ww=matsnf(M1,1);uu=ww[1];vv=ww[2];dd=ww[3];$\backslash$\hfill

\noindent
mm=matsize(dd)[1];nn=matsize(dd)[2];$\backslash$\hfill

\noindent
nnn=nn;for(i=1,nn,if(dd[,i]==0,nnn=nnn-1));$\backslash$\hfill

\noindent
VV=matrix(nn,nnn);$\backslash$\hfill

\noindent
nnnn=0;for(i=1,nn,if(dd[,i]==0,,nnnn=nnnn+1;VV[,nnnn]=vv[,i]));$\backslash$\hfill

\noindent
M2=M1*VV;MM=M2*gg;$\backslash$\hfill

\noindent
kill(gg);kill(M1);kill(ww);kill(uu);kill(vv);kill(dd);$\backslash$\hfill

\noindent
kill(mm);kill(nn);kill(nnn);kill(nnnn);kill(M2);$\backslash$\hfill

\noindent
L1=MM;kill(VV);N2=L1\^{}-1*N1;$\backslash$\hfill

\noindent
ww=matsnf(N2,1);uu=ww[1];vv=ww[2];dd=ww[3];$\backslash$\hfill

\noindent
N3=N2*vv;mm=matsize(dd)[1];nn=matsize(dd)[2];$\backslash$\hfill

\noindent
nnn=nn;for(i=1,nn,if(dd[,i]==0,nnn=nnn-1));$\backslash$\hfill

\noindent
N4=matrix(mm,nnn);nnnn=0;$\backslash$\hfill

\noindent
for(i=1,nn,if(dd[,i]==0,,nnnn=nnnn+1;$\backslash$\hfill

\noindent
ddd=gcd(dd[,i]);N4[,nnnn]=N3[,i]/ddd));Npr=L1*N4;$\backslash$\hfill

\noindent
kill(ggg);kill(N1);kill(M);kill(L1);kill(MM);$\backslash$\hfill

\noindent
kill(N2);kill(ww);kill(uu);kill(vv);kill(dd);kill(N3);$\backslash$\hfill

\noindent
kill(mm);kill(nn);kill(nnn);kill(nnnn);kill(ddd);kill(N4);$\backslash$\hfill

\noindent
SUBLpr1=Npr;R=r;B=SUBLpr1;l=B$\widetilde{\ }\,$*R*B;$\backslash$\hfill

\noindent
ww=matsnf(l,1);uu=ww[1];vv=ww[2];dd=ww[3];$\backslash$\hfill

\noindent
nn=matsize(l)[1];nnn=nn;for(i=1,nn,if(dd[i,i]==0,nnn=nnn-1));$\backslash$\hfill

\noindent
b=matrix(nn,nnn,X,Y,vv[X,Y+nn-nnn]);$\backslash$\hfill

\noindent
ll=b$\widetilde{\ }\,$*l*b;d=vector(nnn,X,dd[X+nn-nnn,X+nn-nnn]);$\backslash$\hfill

\noindent
kill(ww);kill(uu);kill(vv);kill(dd);kill(nn);kill(nnn);$\backslash$\hfill

\noindent
BB=B*b;G=BB$\widetilde{\ }\,$*R*BB;D=d;$\backslash$\hfill

\noindent
SUBLpr=BB;DSUBLpr=D;rSUBLpr=G;$\backslash$\hfill

\noindent
RANK=matsize(DSUBLpr)[2];$\backslash$\hfill

\noindent
RANKA=0;for(k1=1,RANK,if(DSUBLpr[k1]$>$1,RANKA=RANKA+1));$\backslash$\hfill

\noindent
if(RANK$>$19$|\,|$RANK+RANKA$>$22,,$\backslash$\hfill

\noindent
DISCGRS=vector(RANKA,X,DSUBLpr[X]);$\backslash$\hfill

\noindent
if(RANK+RANKA$<$22,print("t=",t,"DISCGRS=",DISCGRS),$\backslash$\hfill

\noindent
pp=factor(DISCGRS[RANKA]);$\backslash$\hfill

\noindent
u1=2;if(pp[1,1]$>$2,u1=1);if(pp[1,1]==2\&\&pp[1,2]$>$1,u1=1);$\backslash$\hfill

\noindent
if(pp[1,1]==2\&\&pp[1,2]==1,delta=0;$\backslash$\hfill

\noindent
for(k=1,RANKA,if(type(rSUBLpr[k,k]/4)!="t\_{}INT",delta=1)));$\backslash$\hfill

\noindent
if(delta==0,u1=1);gamm=0;$\backslash$\hfill

\noindent
BBB=matrix(RANK-RANKA,RANK-RANKA,X,Y,rSUBLpr[X+RANKA,Y+RANKA]);$\backslash$\hfill

\noindent
dBBB=matdet(BBB);$\backslash$\hfill

\noindent
for(t1=u1,matsize(pp)[1],if(pp[t1,1]==2\&\&$\backslash$\hfill

\noindent
Mod(dBBB,8)!=Mod(1,8)\&\&Mod(dBBB,8)!=Mod(-1,8),gamm=1);$\backslash$\hfill

\noindent
if(pp[t1,1]!=2\&\&kronecker((-1)\^{}(RANK+1)*dBBB,pp[t1,1])!=1,gamm=1));$\backslash$\hfill

\noindent
if(gamm==0,print("t=",t);print("DISCGRS=",DISCGRS)))));\hfill



V.V. Nikulin \par Deptm. of Pure Mathem. The University of
Liverpool, Liverpool\par L69 3BX, UK; \vskip1pt Steklov
Mathematical Institute,\par ul. Gubkina 8, Moscow 117966, GSP-1,
Russia

vnikulin@liv.ac.uk \, \ \ \  vvnikulin@list.ru

Personal page: http://vnikulin.com

\end{document}